\documentclass[12pt,a4paper,reqno]{amsart}

\usepackage[english]{babel}
\usepackage[T1]{fontenc}
\pdfoutput=1

\usepackage[dvipsnames]{xcolor}
\definecolor{MyBlue}{cmyk}{1,0.13,0,0.63}
\definecolor{MyGreen}{cmyk}{0.91,0,0.88,0.52}
\definecolor{MyRed}{rgb}{.6,0,0}

\usepackage{hyperref}
\usepackage{xurl}
\hypersetup{breaklinks=true,colorlinks = true,
linkcolor = MyBlue,
urlcolor  = MyRed,
citecolor = MyGreen}
\usepackage[nameinlink]{cleveref}

\usepackage{amsmath,amsfonts,amsthm,amssymb,mathrsfs,mathtools}
\usepackage{tikz-cd}

\usepackage[margin=3cm]{geometry}
\usepackage[shortlabels]{enumitem}
\setenumerate[1]{label=(\arabic*),font=\upshape}

\usepackage[noindentafter]{titlesec}
\titleformat{name=\section}{}{\bf\thetitle.}{0.5em}{\centering\bf}
\titlespacing{\section}{0em}{2.5em}{1em}
\titleformat{name=\subsection}[block]{}{\bfseries\thetitle.}{0.5em}{\itshape\bf}
\titleformat{name=\subsubsection}[runin]{}{\itshape\thetitle.}{0.5em}{\itshape}[.]

\numberwithin{equation}{section}
\newtheorem{theorem}{Theorem}[section]
\newtheorem*{theorem*}{Theorem}
\newtheorem{lemma}[theorem]{Lemma}
\newtheorem{corollary}[theorem]{Corollary}
\newtheorem{proposition}[theorem]{Proposition}
\theoremstyle{definition}
\newtheorem{definition}[theorem]{Definition}
\newtheorem{remark}[theorem]{Remark}
\newtheorem{example}[theorem]{Example}

\newcommand{\alg}[1]{\mathcal{#1}} 
\newcommand{\adj}[2][\alg{A}]{\mathcal L_{#1}\left(#2\right)} 
\newcommand{\comp}[2][\alg{A}]{\mathcal K_{#1}(#2)} 
\newcommand{\frank}[2][\alg{A}]{\mathcal F_{#1}(#2)} 
\newcommand{\reg}[2][\alg{A}]{\mathcal R_{#1}{(#2)}} 
\newcommand{\norm}[2][]{\|#2\|_{#1}} 
\newcommand{\Kclass}[1]{[#1]_0} 
\newcommand{\K}[1][\alg{A}]{K_0(#1)} 
\newcommand{\ind}[1]{\mathrm{ind}\left(#1\right)} 
\newcommand{\weakind}[1]{\widehat{\mathrm{ind}}\left(#1\right)} 
\newcommand{\complex}[2][\md E]{(#1,#2)} 
\newcommand{\gapmet}[2]{\mathrm{d}_{\mathrm{gap}}(#1,#2)} 
\newcommand{\tilgapmet}[2]{\mathrm{gap}(#1,#2)} 
\newcommand{\Rieszmet}[2]{\mathrm{d}_{\mathrm{R}}(#1,#2)} 
\newcommand{\md}[1]{\mathscr{#1}} 
\newcommand{\dom}[1]{\mathrm{dom}{(#1)}} 
\newcommand{\ran}[1]{\mathrm{ran}{(#1)}} 
\renewcommand{\ker}[1]{\mathrm{ker}{(#1)}}
\newcommand{\coH}[2]{\mathcal H_{#1}(#2)} 
\newcommand{\ev}[3][]{\langle #2 \mid #3 \rangle_{#1}} 
\newcommand{\graph}[1]{\Gamma(#1)} 
\newcommand {\airquotes}[1]{``\emph{#1}''}
\newcommand{\evDirac}[1]{\mathcal D^+_{#1}}
\newcommand{\oddDirac}[1]{\mathcal D^-_{#1}}
\newcommand{\Dirac}[1]{\mathcal D_{#1}}
\newcommand{\Laplace}[2][\md E]{\Delta^{\complex[#1]{#2}}}
\renewcommand{\sharp}{*}
	
\title{Fredholm complexes of Hilbert C*-modules}
\author{Brian Villegas-Villalpando}
\author{Koen van den Dungen}
\date{}
\thanks{This paper is based on the first named author’s master’s thesis \cite{Villegas-Villalpando2024} at the University of Bonn written under the supervision of the second named author.}

\begin{document}
	\begin{abstract}
		We investigate \emph{complexes of Hilbert C*-modules}, which are cochain complexes with (unbounded) regular operators between Hilbert C*-modules as differential maps. In particular, we provide various equivalent characterizations of the Fredholm property for such complexes of Hilbert C*-modules, and we define the Fredholm index taking values in the $K$-theory group of the C*-algebra. Among other properties of this index, we prove the stability under small or relatively compact perturbations, and we obtain alternative expressions for the index under the existence of a (weak or strong) Hodge decomposition.
	\end{abstract}	
	\maketitle
    
    \addcontentsline{toc}{section}{Introduction}
    \section*{Introduction}\label{sec:introduction}\noindent
                Segal gave in \cite{Segal1970} one of the earliest treatments on \emph{Fredholm complexes} whose \emph{index} takes values in an abstract $K$-group. His work shows that any \emph{Fredholm complex of vector bundles} on a compact space $X$ can be associated with an \emph{Euler characteristic} in $K(X)$. We now introduce the \emph{noncommutative counterpart} of this theory.
					
                Cochain complexes, or complexes for short, are sequences of morphisms $\{d_k\}_{k\in \mathbb Z}$, known as \emph{differential maps}, such that $d_{k+1}d_k=0$ for every $k\in \mathbb Z$. The study of complexes with an analytic structure goes back at least as far as Atiyah and Bott \cite{AtiyahBott1967} in the context of elliptic operators, and, by using Hilbert-space methods, an important class of complexes emerges: \emph{Hilbert complexes}. These are rich structures that prove their usefulness in global analysis, combining functional analysis with homological algebra, see, e.g., \cite{BruningLesch1992}.
					
				If one is concerned with continuous families of Hilbert complexes, then the natural setting for such families is (by Gelfand duality) that of Hilbert C*-modules over an algebra of complex-valued continuous functions, which are given by the sections of continuous fields of Hilbert spaces \cite{DixmierDouady1963}. In this paper, we provide a systematic treatment of a more general situation, namely we investigate \emph{complexes of Hilbert C*-modules} over an arbitrary (possibly noncommutative) C*-algebra. In short, such a complex takes the form
				\[\begin{tikzcd}
					0 & {\md E_0} & {\md E_1} & \cdots & {\md E_{N}} & {\md{E}_{N+1}} & 0,
					\arrow[from=1-1, to=1-2]
					\arrow["{t_0}", from=1-2, to=1-3]
					\arrow["{t_1}", from=1-3, to=1-4]
					\arrow["{t_{N-1}}", from=1-4, to=1-5]
					\arrow["{t_{N}}", from=1-5, to=1-6]
					\arrow[from=1-6, to=1-7]
				\end{tikzcd}\]
				where each $\md E_k$ is a Hilbert C*-module (over the same C*-algebra), and each (densely defined) differential map $t_{k}\colon \dom{t_k}\subset \md E_k\to \md E_{k+1}$ is a \emph{regular} operator such that $\ran{t_k}\subset \ker{t_{k+1}}$ for all $0\leq k\leq N-1$.
															
				Our goal is to define and study the \emph{Fredholm property} of complexes of Hilbert C*-modules, and to show how to associate a \emph{Fredholm index} that takes values in the $K$-group of the C*-algebra. To this end, we shall prove that it suffices to look at a single operator, which we will call the \emph{even Dirac operator}. 
				Although our discussion of complexes of Hilbert C*-modules is largely in parallel to Hilbert complexes, the lack of a general Hilbert projection theorem (not every closed submodule of a Hilbert C*-module is orthogonally complemented) leads to several difficulties and requires us to formulate alternative proofs.
					
				This paper is organized as follows. 
                We start in \Cref{sec:main-results} with our main definitions and a summary of our main results. 
				In \Cref{sec:notation}, we collect some useful facts regarding regular operators on Hilbert C*-modules. \Cref{sec:complexes_hilbert_modules,sec:weak_hodge_decomposition,sec:maps_complexes} provide a systematic treatment of complexes of Hilbert C*-modules. First, in \Cref{sec:complexes_hilbert_modules}, given a complex of Hilbert C*-modules, we will construct the bounded transform complex, the adjoint complex, the graph-norm complex, and the \emph{Dirac} and \emph{Laplace} operators. Afterwards, in \Cref{sec:weak_hodge_decomposition}, we characterize complexes that have weak and strong Hodge decomposition. \Cref{sec:maps_complexes} briefly discusses maps of complexes. 
				In \Cref{sec:fredholm_property}, we study Fredholm complexes in the context of Hilbert C*-modules, and in \Cref{sec:fredholm-equivalences} we provide various equivalent characterizations of the Fredholm property. Afterwards, in \Cref{sec:fredholm_index_map}, we study the index map that assigns to any Fredholm complex an element in the $K$-group of the C*-algebra, and provide alternative expressions for the index under the assumption of weak or strong Hodge decomposition. 
				We describe some stability properties of the index in \Cref{sec:stability}, and in \Cref{sec:sums-seq-tensor} we consider direct sums, sequences, and tensor products of Fredholm complexes. 
				Lastly, in \Cref{sec:applications}, we consider two special cases of the general theory of C*-Fredholm complexes.

    \section{Main definitions and main results}\label{sec:main-results}
				Throughout this paper, let us fix a complex C*-algebra $\alg A$ with norm $\norm{\cdot}$, and let us denote by $\K$ the $K$-group of the C*-algebra $\alg A$. 
				In this section, we present our main definitions and constructions, and provide a summary of our main results (all proofs are postponed to later sections). 
							
		\subsection{Complexes of Hilbert C*-modules}\label{subsec:C*-hilbert-complexes}
								
				We can define a complex of Hilbert C*-modules in analogy with the Hilbert-space case, see, e.g., \cite[Section 7]{AtiyahSinger1968} and \cite[Section 2]{BruningLesch1992}, which we recall first. 
				A \emph{Hilbert complex} is given by a (finite) collection of Hilbert spaces $\md H_{0}, \md H_1, \cdots, \md H_{N+1}$ along with a collection of closed, densely defined, linear operators 
				\begin{equation*}
					\{d_k\colon \dom{d_k}\subset \md H_k\to \md H_{k+1}\}_{k=0}^N
				\end{equation*} 
				such that $\ran{d_{k}}\subset \ker{d_{k+1}}$ for every $0\leq k\leq N-1$. 
				This notion has been present in the literature for a long time, and it has been studied in particular for its relation with elliptic complexes. 
													
				We shall generalize Hilbert complexes by considering instead a sequence of Hilbert $\alg A$-modules $\{\md E_k\}_{k\in \mathbb Z}$ over a fixed C*-algebra $\alg A$. 
																
				\begin{definition}\label{def:C*-hilbert_complexes}
					A \emph{complex of Hilbert $\alg A$-modules}, or an \emph{$\alg A$-Hilbert complex} for short, is a sequence of regular operators $\{t_k\in \reg{\md E_k,\md E_{k+1}}\}_{k\in \mathbb Z}$ that satisfies the complex property 
					\begin{equation}\label{eq1:subsec:C*-hilbert-complexes}		
				    	\ran{t_k}\subset \ker{t_{k+1}}, \quad \forall k\in \mathbb Z.
					\end{equation}
					We will use the short notation
					\begin{equation*}
					   \complex{t}:=\{t_k\in \reg{\md E_k,\md E_{k+1}}\}_{k\in \mathbb Z} 
					\end{equation*}
					to denote an $\alg A$-Hilbert complex, and we refer to the maps $\{t_k\}_{k\in \mathbb Z}$ as the \emph{differential maps} of $\complex{t}$. 
				\end{definition}
				Observe that the sequence $\{t_k\colon  \dom{t_k}\to \dom{t_{k+1}}\}_{k\in\mathbb Z}$ is a (cochain) complex of $\alg A$-modules in the sense of homological algebra, meaning that each operator $t_k\colon \dom{t_k}\to \dom{t_{k+1}}$ is an $\alg A$-module map, and that the composition $t_{k+1}t_{k}=0$ is well defined and vanishes identically for every $k\in \mathbb Z$. (Observe that we do not follow the standard notation of superscripts to denote cochain complexes.)
																
				On Hilbert spaces, any closed and densely defined operator is regular, see, e.g.,  \cite[Theorem 5.1.9]{Pedersen1989}. Thus, with our terminology, Hilbert complexes are $\mathbb C$-Hilbert complexes, but we will continue referring to them as Hilbert complexes. 
																
				While many of our results are valid for infinitely long complexes, we will mostly be interested in complexes with only finitely many differential maps. 
				We will call $\complex{t}$ a \emph{finite-length} $\alg A$-Hilbert complex if $\md E_k=0$ when $|k|$ is large enough. In this case, we will assume (without loss of generality) that $\complex{t}$ has the form
				\[\begin{tikzcd}
					0 & {\dom{t_0}} & {\dom{t_1}} & \cdots & {\dom{t_N}} & {\md{E}_{N+1}} & 0.
					\arrow[from=1-1, to=1-2]
					\arrow["{t_0}", from=1-2, to=1-3]
				    \arrow["{t_1}", from=1-3, to=1-4]
					\arrow["{t_{N-1}}", from=1-4, to=1-5]
					\arrow["{t_{N}}", from=1-5, to=1-6]
					\arrow[from=1-6, to=1-7]
				\end{tikzcd}\]																
												
		\subsection{The Dirac operator}\label{sec:main:Dirac}
															
				The following construction, which mimics what is done in the case of elliptic complexes \cite[Section 7]{AtiyahSinger1968} and in the case of Hilbert complexes \cite[Section 2]{BruningLesch1992}, allows us to describe both the Fredholm property as well as the Fredholm index of a finite-length complex of Hilbert C*-modules in terms of a single operator. 
				Let us fix a finite-length $\alg A$-Hilbert complex $\complex{t}:=\{t_k\in \reg{\md E_k,\md E_{k+1}}\}_{k=0}^N$. 
				Put
				\begin{equation*}
					\md E_{\text{ev}}:=\bigoplus_{k} \md E_{2k}\quad\text{and}\quad\md E_{\text{odd}}:=\bigoplus_{k} \md E_{2k+1}
				\end{equation*}
				as (finite) orthogonal sums of Hilbert C*-modules. 
																
				\begin{definition}\label{def:even-odd-Dirac}
					Given an $\alg A$-Hilbert complex $\complex{t}$, we consider the submodules 
					\begin{align*}
						\dom{\evDirac{t}}&:=\bigoplus_{k}\dom{t_{2k}}\cap \dom{t_{2k-1}^*}\subset \md E_{\mathrm{ev}} , \\ 
						\dom{\Dirac{t}^-}&:=\bigoplus_{k}\dom{t_{2k}^*}\cap \dom{t_{2k+1}}\subset \md E_{\mathrm{odd}} .
				    \end{align*}
				    We define the \emph{even Dirac operator} $\evDirac{t}\colon \dom{\evDirac{t}}\to \md E_{\mathrm{odd}}$ and the \emph{odd Dirac operator} $\oddDirac{t}\colon  \dom{\oddDirac{t}}\to \md E_{\mathrm{ev}}$ given in matrix forms by
					\begin{equation*}
					      \evDirac{t}:=\begin{pmatrix}
					        t_0&t_1^*&0&\cdots\\
					        0&t_2&t_3^*&\cdots\\
						      0&0&t_4&\cdots\\
						      \vdots&\vdots&\vdots&\ddots
					    \end{pmatrix}
					    \quad\text{and}\quad
					    \oddDirac{t}:=\begin{pmatrix}
						      t_0^*&0&0&\cdots\\
							t_1&t_2^*&0&\cdots\\
							0&t_3&t_4^*&\cdots\\
							\vdots&\vdots&\vdots&\ddots
						\end{pmatrix}.
					\end{equation*}
				\end{definition}
																
				The following statement will be proven in \Cref{subsec:dirac-operator}. 
				\begin{theorem}                                                
					Let $\complex{t}$ be a finite-length $\alg A$-Hilbert complex. Then the even Dirac operator $\evDirac{t}$ is a regular operator with adjoint $\oddDirac{t}$. 
				\end{theorem}
			
        \subsection{Weak and strong Hodge decomposition}
												
				Consider an $\alg A$-Hilbert complex $\complex{t}$. Recall that the symbol $\oplus$ denotes the orthogonal direct sum of pre-Hilbert or Hilbert C*-modules. As a consequence of the equality $\ker{t_k^*} = (\ran{t_k})^\perp$ (see \Cref{thm:kert_perp_rant^*}) and the complex property \eqref{eq1:subsec:C*-hilbert-complexes}	of $\complex{t}$, we observe that, for each $k$, the three subspaces
				\[                                                
					\ker{t_k}\cap \ker{t_{k-1}^*} , \quad 
					\ran{t_k^*} , \quad\text{and}\quad 
					\ran{t_{k-1}}
				\]
				of $\md E_k$ are mutually orthogonal. Hence, we have the inclusions 
				\begin{multline*}
					\big(\ker{t_k}\cap \ker{t_{k-1}^*}\big) \oplus \ran{t_k^*} \oplus \ran{t_{k-1}} \\
					\subset \big(\ker{t_k}\cap \ker{t_{k-1}^*}\big) \oplus \overline{\ran{t_k^*}} \oplus \overline{\ran{t_{k-1}}} 
					\subset \md E_k .
				\end{multline*}
				We now consider the special cases where these inclusions might be equalities. 
				\begin{definition}\label{def:hodge-decomposition}
					Let $\complex{t}$ be an $\alg A$-Hilbert complex. 
					\begin{enumerate}[(i)]
						\item The complex $\complex{t}$ has \emph{weak Hodge decomposition} if 
						      \begin{equation*}
                              \md E_k=(\ker{t_k}\cap \ker{t_{k-1}^*})\oplus \overline{\ran{t_k^*}}\oplus \overline{\ran{t_{k-1}}}, \quad \forall k\in\mathbb Z.
						      \end{equation*}
						\item The complex $\complex{t}$ has \emph{strong Hodge decomposition} if
							\begin{equation*}
								\md E_k=(\ker {t_k}\cap \ker {t_{k-1}^*})\oplus {\ran{t_k^*}}\oplus {\ran{t_{k-1}}}, \quad \forall k\in\mathbb Z.
							\end{equation*}
					\end{enumerate}
				\end{definition}
																
				We observe that the weak Hodge decomposition for a complex of length one, $0 \to \md E_0 \xrightarrow{t} \md E_1 \to 0$, simply means that $\overline{\ran t}$ and $\overline{\ran{t^*}}$ are orthogonally complemented, which is equivalent to $t$ being polar decomposable (see \Cref{thm:polar-decomposition}). Similarly, the strong Hodge decomposition simply requires $\ran t$ and $\ran{t^*}$ to be closed. In general, we can characterize weak or strong Hodge decompositions as follows (the proofs are given in \Cref{thm:weak-hodge-decomposition-equivalences} and \Cref{thm:hodge-decomposition-equivalences}). 
																
				\begin{theorem}[Weak Hodge decomposition]\label{thm:weak-hodge-decomposition-equivalences-main}
					Let $\complex{t}$ be an $\alg A$-Hilbert complex. 
					\begin{enumerate}[(i)]
						\item The complex $\complex{t}$ has weak Hodge decomposition if and only if, for every $k$, the differential map $t_k$ is polar decomposable. 
						\item If $\complex{t}$ has finite length, then $\complex{t}$ has weak Hodge decomposition if and only if the even Dirac operator $\evDirac{t}$ is polar decomposable. 
					\end{enumerate}
				\end{theorem}
																
				\begin{theorem}[Strong Hodge decomposition]\label{thm:hodge-decomposition-equivalences-main}
					Let $\complex{t}$ be an $\alg A$-Hilbert complex. 
					\begin{enumerate}[(i)]
						\item The complex $\complex{t}$ has strong Hodge decomposition if and only if, for every $k$, the differential map $t_k$ has closed range. 
						\item If $\complex{t}$ has finite length, then $\complex{t}$ has strong Hodge decomposition if and only if the even Dirac operator $\evDirac{t}$ has closed range. 
					\end{enumerate}
				\end{theorem}
								
				In the case of Hilbert spaces, every densely defined, closed operator is polar decomposable, and, therefore, every Hilbert complex has weak Hodge decomposition (see also \cite[Lemma 2.1]{BruningLesch1992}). In general, it is known that not every (adjointable or regular) operator between  Hilbert C*-modules has a polar decomposition (see \Cref{example:Fredholm-no-Hodge}). Therefore, we can not expect every $\alg A$-Hilbert complex to have weak Hodge decomposition.  However, we will see in \Cref{thm:weak-hodge-decomposition-subalgebra-compact-operators} that, if the C*-algebra $\alg A$ is $*$-isomorphic to a subalgebra of compact operators on a Hilbert space, then every $\alg A$-Hilbert complex has weak Hodge decomposition.
																
				As before, not every $\alg A$-Hilbert complex has strong Hodge decomposition. However, if a complex $\complex{t}$ has strong Hodge decomposition, then its cohomology groups
				\begin{equation*}
				    \coH{k}{\complex{t}}:=\ker{t_k}/\ran{t_{k-1}}, \quad k\in \mathbb Z,
				\end{equation*}
				are well-defined Hilbert $\alg A$-modules. The following statement is proven in \Cref{subsec:topological-complexes}. 
				\begin{theorem}                                                
					If $\complex{t}$ is an $\alg A$-Hilbert complex with strong Hodge decomposition, then we have the isomorphism 
					\[                                                
						\coH{k}{\complex{t}} \cong \ker{t_k} \cap \ker{t_{k-1}^*} . 
				    \]
				\end{theorem}
		
		\subsection{The Fredholm property}\label{sec:main:Fredholm-property}
													
		      	A Hilbert complex $\complex[\md H]{d}:=\{d_k\in\reg[\mathbb C]{\md H_{k},\md H_{k+1}}\}_{k=0}^N$ (with $\alg A = \mathbb{C}$) is called a \emph{Fredholm complex} if its cohomology groups
				\begin{equation}\label{eq2:subsec:C*-hilbert-complexes}
					\coH{k}{\complex[\md H]{d}}:=\ker{d_k}/\ran{d_{k-1}},\quad 0\leq k\leq N+1,
				\end{equation}
				are finite-dimensional Hilbert spaces. In this case, every Fredholm complex has strong Hodge decomposition \cite[Corollary 2.5]{BruningLesch1992}, and its \emph{Fredholm index} can be defined as follows:
				\begin{equation}\label{eq3:subsec:C*-hilbert-complexes}
					\ind{\complex[\md H]{d}}:=\sum_{k=0}^{N+1} (-1)^k \dim \coH{k}{\complex[\md H]{d}}\in \mathbb Z\cong \K[\mathbb C].
				\end{equation}
				For a general $\alg A$-Hilbert complex, the cohomology groups might not be well-defined, and we consider instead an alternative definition of the Fredholm property which is independent of the cohomology groups and which does not require a (strong or weak) Hodge decomposition. 
																
				The Fredholm property of an adjointable Hilbert complex can be equivalently formulated in terms of \emph{parametrices}, see, e.g, \cite[Section 3.2.3.1]{RempelSchulze1982}. Similarly, adjointable operators between Hilbert C*-modules are called C*-Fredholm if they have a \emph{parametrix} (\Cref{def:adjointable-fredholm-operators} below). Moreover, Joachim defined in \cite[Definition 2.1]{Joachim2003} that a regular operator $t\in \reg{\md E,\md F}$ is Fredholm if it has a \emph{pseudo-left} and a \emph{pseudo-right inverse}. Here we will consider the following general definition for the Fredholm property of a finite-length $\alg A$-Hilbert complex. 
				\begin{definition}\label{def:joint-parametrix}
					Let $\complex{t}:=\{t_k\in \reg{\md E_k,\md E_{k+1}}\}_{k=0}^N$ be a finite-length $\alg A$-Hilbert complex. 
					A collection of adjointable operators $\{P_{l,k},P_{r,k}\in \adj{\md E_{k+1},\md E_k}\}_{k=0}^N$ is called a \emph{joint parametrix} of $\complex{t}$ if there exist compact operators $\{C_k\in \comp{\md E_k}\}_{k=0}^{N+1}$ such that
					\begin{align}
				        P_{r,k-1}(\dom{t_k})&\subset\dom{t_{k-1}},\label{eq1:def:joint-parametrix-complex}\\
						P_{l,k}t_k+t_{k-1}P_{r,k-1}&=1_{\md E_k}-C_k \text{ in }\dom{t_k},\label{eq2:def:joint-parametrix-complex}
					\end{align}		
					for every $0\leq k\leq N+1$. We assume that $P_{l,k}$, $P_{r,k}$, and $t_k$ are zero if $k<0$ or $N<k$, and we follow this assumption in similar situations. The joint parametrix will be denoted by $\complex{P_l,P_r}:=\{P_{l,k},P_{r,k}\}_{k=0}^N$, where the sequence $\{C_k\}_{k=0}^N$ is considered implicitly. We refer to conditions \eqref{eq1:def:joint-parametrix-complex} and \eqref{eq2:def:joint-parametrix-complex} as \emph{the Fredholm property} of $\complex{t}$.
				\end{definition}
																
				\begin{definition}\label{def:fredholm-complex}
					A  Fredholm complex of Hilbert $\alg A$-modules, or an \emph{$\alg A$-Fredholm complex} for short, is a finite-length $\alg A$-Hilbert complex that has a joint parametrix.
				\end{definition}
					We will prove in \Cref{thm:fredholm_complexes_of_hilbert_spaces} that our definition of an $\alg A$-Fredholm complex is compatible (in the special case $\alg A = \mathbb C$) with the definition of a Fredholm complex of Hilbert spaces. 
													
				\begin{remark}
					If we consider a regular operator $t\in \reg{\md E,\md F}$ as a complex of length one, then $t$ is an $\alg A$-Fredholm complex if and only if there exist adjointable operators $P_l,P_r\in \adj{\md F,\md E}$, and compact operators $C_l\in \comp{\md E}$ and $C_r\in \comp{F}$ such that
					\begin{equation*}
				    	P_lt\subset 1_{\md E}-C_l\quad\text{and}\quad tP_r=1_{\md F}-C_r,
					\end{equation*}
					where $\ran{P_r}\subset \dom{t}$. In this case, we simply say that $t$ is $\alg A$-Fredholm or an $\alg A$-Fredholm operator, and we denote by $(P_l,P_r)$ its joint parametrix. 
					We will see in \Cref{thm:consistency-joachim} that this definition of the Fredholm property for regular operators is consistent with \cite[Definition 2.1]{Joachim2003}. 
				\end{remark}
																
				It is a very useful result that the Fredholm property of a complex can be characterized by the Fredholm property of a single regular operator, namely its even Dirac operator, which we will prove in \Cref{thm:fredholm-goal-I}. 
																
				\begin{theorem}
					A finite-length $\alg A$-Hilbert complex $\complex{t}$ is $\alg A$-Fredholm if and only if its even Dirac operator $\evDirac{t}$ is $\alg A$-Fredholm.
				\end{theorem}
																
		\subsection{The Fredholm index}
				We consider the Fredholm index map defined by \cite[Section 3]{Exel1993}, which assigns to any adjointable $\alg A$-Fredholm operator $T\in \adj{\md E,\md F}$ an element $\ind{T} \in \K$. Since a regular operator $t\in \reg{\md E,\md F}$ is $\alg A$-Fredholm if and only if its bounded transform $F_t\in \adj{\md E,\md F}$ is $\alg A$-Fredholm (see \Cref{thm:fredholm-equivalence-complex-and-bounded-transform}), we may define the index of a regular $\alg A$-Fredholm operator $t$ as 
				\begin{equation*}
					\ind{t}:=\ind{F_t}\in \K.
				\end{equation*}
				Now, if $\complex{t}$ is an $\alg A$-Fredholm complex, we may then apply this index map to the even Dirac operator $\evDirac{t}$. This leads us to the following definition. 
				\begin{definition}\label{def:Fredholm-index-complex}
					The \emph{Fredholm index} of an $\alg A$-Fredholm complex $\complex{t}$ is defined as 
					\begin{equation*}
						\ind{\complex{t}}:=\ind{\evDirac{t}}\in \K.
				    \end{equation*}
				\end{definition}
				We observe that for a complex of length one, $0 \to \md E_0 \xrightarrow{t} \md E_1 \to 0$, we have $\evDirac{t} = t$, so that the index of the complex agrees with the index of the regular operator $t\in \reg{\md E_0,\md E_1}$. 
																
				The following statements list some of the main properties of the Fredholm index, all of which will be proven in \Cref{sec:fredholm_index_map,sec:stability}.
																
				\begin{theorem}[{see \Cref{thm:index-adjoint}}]
				    If $\complex{t}$ is an $\alg A$-Fredholm complex, then the adjoint complex $\complex{t^\sharp}$, which is the (chain) complex whose differential maps are the adjoints of the differential maps of $\complex{t}$, is also $\alg A$-Fredholm and
					\begin{equation*}
						\ind{\complex{t^\sharp}}=(-1)^{N+1}\ind{\complex{t}}\in \K.
					\end{equation*}
				\end{theorem}
																
				\begin{theorem}[{see \Cref{thm:index_weak_hodge_decomposition} and \Cref{thm:index_topological_complex}}]
					Let $\complex{t}$ be an $\alg A$-Fredholm complex. 
					\begin{enumerate}[(i)]
						\item If $\complex{t}$ has weak Hodge decomposition, then
				    		\begin{equation*}
								\ind{\complex{t}}=\Kclass{\ker{\evDirac{t}}}-\Kclass{\ker{\oddDirac{t}}}\in \K.
							\end{equation*}
						\item If $\complex{t}$ has strong Hodge decomposition, then
					        \begin{equation*}
								\ind{\complex{t}}=\sum_{k=0}^{N+1} (-1)^k \Kclass{\coH{k}{\complex{t}}}\in \K.
							\end{equation*}
					\end{enumerate}
					Here, the symbol $\Kclass{\md F}$ denotes the K-theory class in $\K$ of a finite-rank Hilbert $\alg A$-module $\md F$. 
				\end{theorem}
																
				\begin{theorem}[{see \Cref{thm:small_stability_index_complex} and \Cref{thm:homotopy_fredholm_complexes}}]
					The index of an $\alg A$-Fredholm complex $\complex{t}$ is stable under small perturbations with respect to the gap topology: there exists an $\epsilon>0$ such that
					\begin{equation*}
						\ind{\complex{t}}=\ind{\complex{s}}\in \K
					\end{equation*}
					for any $\alg A$-Hilbert complex $\complex{s}$ that satisfies the inequality
					\begin{equation*}
						\sup_{0\leq k\leq N} \gapmet{t_k}{s_k}<\epsilon.
					\end{equation*}
					In particular, homotopic $\alg A$-Fredholm complexes (with respect to the gap topology) have the same index.
				\end{theorem}
																
				\begin{theorem}[{see \Cref{thm:relative-compact-stability-index-complex}}]
					The index of an $\alg A$-Fredholm complex $\complex{t}$ is stable under relatively compact perturbations: given a sequence of semi-regular operators $\{r_k:\dom{r_k}\subset \md E_{k}\to \md E_{k+1}\}_{k}$ such that each $r_k$ is relatively $t_k$-compact and each $r_k^*$ is relatively $t_k^*$-compact, we have the equality 
					\begin{equation*}
					   \ind{\complex{t}}=\ind{\complex{t+r}}\in \K
					\end{equation*} 
					whenever $\complex{t+r}:=\{t_k+r_k\}_{k\in\mathbb Z}$ satisfies the complex property \eqref{eq1:subsec:C*-hilbert-complexes}.
				\end{theorem}
	
    \section{Preliminaries}\label{sec:notation}
									
		\subsection{Hilbert C*-modules}\label{subsec:hilbert-C*-modules}
				A (right) \emph{pre-Hilbert C*-module} over $\alg A$ is a complex vector space $\md E$ that is a right $\alg A$-module (with compatible scalar multiplication) together with an $\alg A$-valued inner product $\ev[\md E]{\cdot}{\cdot}\colon \md E\times \md E\to \alg A$. This map is $\alg A$-linear in the second variable, and, for every $x,y\in \md E$, it satisfies that
				\begin{equation*}
					\ev[\md E]{x}{y}=\ev[\md E]{y}{x}^*\quad \text{and}\quad \ev[\md E]{x}{x}\geq 0\text{ with equality only if }x=0. 
				\end{equation*}
				Any pre-Hilbert $\alg A$-module $\md E$ induces a norm $x\mapsto \norm[\md E]{x}:=\norm{\ev[\md E]{x}{x}}^{1/2}$. We say that $\md E$ is a \emph{Hilbert C*-module} over $\alg A$, or a Hilbert $\alg A$-module for short, if $\md E$ is complete with respect to this norm. The symbols $\md E$ and $\md F$ will stand for Hilbert $\alg A$-modules, unless stated otherwise. We rarely use subscripts to distinguish norms or inner products; for example, we write $\ev{\cdot}{\cdot}$ for the $\alg A$-valued inner product of $\md E$ that induces the norm $\norm{x}:=\norm{\ev{x}{x}}^{1/2}$ for all $x\in \md E$. 
				
				For a submodule $\md S \subset \md E$, we denote by $\md S^\perp$ its orthogonal submodule:
				\[                        
				    \md S^\perp := \{ e\in\md E : \ev[\md E]{e}{s}=0 \;\forall s\in\md S \} .
				\]
				The \emph{orthogonal sum} $\md E\oplus \md F$ is the Hilbert $\alg A$-module of ordered pairs $x\oplus y$ with inner product
				\begin{equation*}
					\ev{x_1\oplus y_1}{x_2\oplus y_2}:=\ev{x_1}{x_2}+\ev{y_1}{y_2}, \quad \forall x\in \md E,\; y\in \md F.
				\end{equation*}
				For an introduction to Hilbert C*-modules, we refer to \cite{Lance1995}.
												
		\subsection{Operators on Hilbert C*-modules}\label{subsec:operators-hilbert-C*-modules}
				
				The space of adjointable operators from $\md E$ into $\md F$ is denoted by $\adj{\md E,\md F}$. We write $T^*$ for the adjoint of $T$. We say that an operator $T\in \adj{\md E,\md F}$ is unitary if it is bijective and $T^*=T^{-1}$. Given elements $x\in \md F$ and $y\in \md E$, the map $\theta_{x,y}(z):=x\ev{y}{z}$, for $z\in \md E$, is an adjointable operator in $\adj{\md E,\md F}$ with adjoint $\theta_{y,x}\in \adj{\md F,\md E}$, and their linear span defines the collection $\frank{\md E,\md F}$ of all \emph{finite-rank operators} from $\md E$ into $\md F$. The closure of $\frank{\md E,\md F}$ in $\adj{\md E,\md F}$ with respect to the operator norm is denoted by $\comp{\md E,\md F}$, and we refer to them as \emph{compact operators}. 
				
				Following \cite[Definitions 2.2 and 2.4]{Exel1993}, we say that $\md E$ is a \emph{finite-rank Hilbert $\alg A$-module} if the identity operator $1_{\md E}$ is in $\comp{\md E,\md E}$. In this case, we denote by $\Kclass{\md E}$ the $K$-class of $\md E$ in $\K$. Equivalently, we note that $\md E$ is finite-rank if and only if there exists an idempotent matrix $p \in M_n(\alg A)$ such that $\md E \cong p \alg A^{\oplus n}$ \cite[Proposition 2.3]{Exel1993}. 
				
				A densely defined ($\alg A$-linear) operator $t\colon\dom{t}\subset \md E\to \md F$ is called \emph{semi-regular} if the adjoint $t^*$ is also densely defined. 
                A semi-regular operator $t$ is called \emph{regular} if $t$ is closed and $1_{\md E}+t^*t$ has dense range. We write $\reg{\md E,\md F}$ for the set of all regular operators from $\md E$ to $\md F$. 
				Given a regular operator $t\in \reg{\md E,\md F}$, we write $\dom{t}$ for the domain of $t$ in $\md E$, $\ran{t}$ for its range in $\md F$, and $\ker t$ for its kernel in $\md E$. Given two densely defined operators $t$ and $s$, we write $s\subset t$ whenever $\dom{s} \subset \dom{t}$ and $sx=tx$ for any $x\in \dom{s}$. In particular, $t\in \reg{\md E,\md E}$ is self-adjoint whenever $t=t^*$ (i.e., $t\subset t^*$ and $t^*\subset t$).
				
				\begin{remark}
					We shall write $\adj{\md E}$ instead of $\adj{\md E,\md E}$, and we will follow this convention for the set of finite-rank, compact, and regular operators.
				\end{remark}

		\subsubsection{The bounded transform}
				
				For any regular operator $t\in \reg{\md E,\md F}$, its adjoint, $t^*$, is regular in $\reg{\md F,\md E}$ (with $t^{**}=t$) and $t^*t$ is self-adjoint in $\reg{\md E}$. Moreover, $1_{\md E}+t^*t$ is invertible, and we can define the operators						
				\begin{equation}\label{eq1:sec:preliminaries}
					\begin{aligned}
					F_t&:=t(1_{\md E}+t^*t)^{-1/2} \in \adj{\md E,\md F},\\
					Q_t&:=(1_{\md E}+t^*t)^{-1/2} \in \adj{\md E}.
				\end{aligned}
				\end{equation}
				We refer to $F_t$ as the \emph{bounded transform} of $t$ and to $Q_t^2$ as the \emph{resolvent} of $t$. The operator $Q_t$ is positive (hence self-adjoint) with $0\leq Q_t\leq 1$ such that $Q_t\colon\md E\to \dom t$ is bijective. Moreover, $\norm{F_t}\leq 1$ and $t=F_t(1_{\md E}-F_t^*F_t)^{-1/2}$. The bounded transform is adjoint preserving: the bounded transform of $t^*$ is given by $F_t^*$. 
				For future reference, we list here the following properties of the bounded transform and the resolvent, and we refer to \cite[Chapters 9 and 10]{Lance1995} for further details.
				
				\begin{proposition}\label{thm:properties-bounded-transform}
					Consider $t\in \reg{\md E,\md F}$. The following relations hold.
					\begin{enumerate}[(i)]
						\item $F_t=tQ_t$. \label{item1:thm:properties-bounded-transform}
						\item $t=F_tQ_t^{-1}$. \label{item2:thm:properties-bounded-transform}
						\item $Q_t^2=1_{\md E}-F_t^*F_t$. \label{item3:thm:properties-bounded-transform}
						\item $Q_tQ_t^{-1}=1_{\dom t}$ and $Q_t^{-1}Q_t=1_{\md E}$. \label{item4:thm:properties-bounded-transform}
						\item $Q_{t}^2Q_{t}^{-2}=1_{\dom{t^*t}}$ and $Q_t^{-2}Q_t^2=1_{\md E}$. \label{item5:thm:properties-bounded-transform}
						\item $Q_t t^*\subset F_t^*$. \label{item6:thm:properties-bounded-transform}
						\item $Q_{t^*}F_t=F_tQ_t$. \label{item7:thm:properties-bounded-transform}
						\item $F_tQ_t^{-1}=Q_{t^*}^{-1}F_t|_{\dom t}$. \label{item8:thm:properties-bounded-transform}
					\end{enumerate}
					Moreover, analogous relations hold for $t^*\in \reg{\md F,\md E}$.
				\end{proposition}
				
				In the following results, we consider basic but important properties of regular operators that we shall need in the sequel.
				
				\begin{proposition}\label{thm:kert_perp_rant^*}\label{thm:invariance-kernel-rank-bounded-transform}
					For any $t\in \reg{\md E,\md F}$, we have the following equalities: 
					\begin{gather*}
						\ker{t^*}=(\ran t)^\perp , \\
						\ker{F_t^*F_t}=\ker{F_t}=\ker{t}=\ker{t^*t} , \\ 
						\ran{F_t}=\ran{t} \quad\text{and}\quad \ran{F_t^*F_t}=\ran{t^*t} , \\
						\overline{\ran{t^*t}}=\overline{\ran{t^*}} .
					\end{gather*}
					\begin{proof}
						Since $\dom{t}$ is dense in $\md E$, we have $\dom{t}^\perp = \md E^\perp=0$. Then $y\in \ker{t^*}$ if and only if $0=\ev{t^*y}{x}=\ev{y}{tx}$ for every $x\in \dom t$, which proves $\ker{t^*}=(\ran t)^\perp$. Since $\ran{Q_t}=\dom{t}$, \Cref{thm:properties-bounded-transform}\ref{item1:thm:properties-bounded-transform} implies that $\ran{t}=\ran{F_t}$ (and similarly for the adjoint). Therefore, using the first part of the proof, $\ker{F_t}=(\ran{F_t^*})^\perp=(\ran{t^*})^\perp=\ker{t}$. 
						The equality $\ker{t} = \ker{t^*t}$ (and similarly the equality $\ker{F_t} = \ker{F_t^*F_t}$) is straightforward to verify. 
						By \Cref{thm:properties-bounded-transform}, $F_{t}^*F_{t}=t^*tQ_{t}^2$ and $\dom{t^*t}=\ran{Q_{t}^2}$. Thus, $\ran{F_t^*F_{t}}=\ran{t^*t}$. 
						Finally, by \cite[Proposition 3.7]{Lance1995}, we have $\overline{\ran{F_t^*F_t}}=\overline{\ran{F_t^*}}$, and the final equality then follows from the previous ones. 
					\end{proof}
				\end{proposition}
				
				\begin{theorem}\label{thm:closed_range}
					Let $t\in \reg{\md E,\md F}$. If any of the operators $t$, $t^*$, $t^*t$, and $tt^*$ has closed range, then all these operators have closed range. 
					In this case, 
					\[                        
						\md E = \ker{t} \oplus \ran{t^*} 
						\quad\text{and}\quad 
						\md F = \ker{t^*} \oplus \ran{t} .
					\]
					\begin{proof}                           
						By \Cref{thm:invariance-kernel-rank-bounded-transform}, it suffices to consider the bounded transforms. 
						We consider the following lemma, cf.\ \cite[Ch.\ IV, Theorem 5.2]{Kato1995}. An adjointable operator $T\in \adj{\md E,\md F}$ has closed range if and only if
						\begin{equation*}
							\gamma(T):=\inf_{x\not\in \ker T} \frac{\norm{Tx}}{\mathrm{dist}(x,\ker{T})}>0.
						\end{equation*} 
						Since $\norm{F_t^*F_tx}\leq\norm{F_t^*}\norm{F_tx}$ for any $x\in \md E$, and $\ker{F_t}=\ker{F_t^*F_t}$, we deduce that $F_t$ has closed range whenever $F_t^*F_t$ has closed range. Moreover, by \cite[Theorem 3.2]{Lance1995}, if $F_t$ has closed range, then $F_t^*$ and $F_t^*F_t$ have closed range. By also considering $F_t^{**}=F_t$, we have shown that the operators $F_t$, $F_t^*$, $F_tF_t^*$, and $F_t^*F_t$ all have closed range if and only if any of them has closed range. Finally, if $F_t$ has closed range, then we know by \cite[Theorem 3.2]{Lance1995} that $\ker{F_t}$ is complemented with orthogonal complement ${\ran{F_t^*}}$, and a similar statement holds for $F_t^*$.  
					\end{proof}
				\end{theorem}

		\subsubsection{The graph inner-product}\label{sec:graph-inner-product}
				
				For any $t\in \reg{\md E,\md F}$, we denote by $\md E_{\graph{t}}$ the Hilbert $\alg A$-module whose underlying vector space is $\dom t$ endowed with the \emph{graph inner-product}:
				\begin{equation*}
					\ev[\md E_{\graph{t}}]{x}{y}:=\ev[\md E]{x}{y}+\ev[\md F]{tx}{ty},\quad \forall x,y\in \dom t.
				\end{equation*}
				The following result is well known.
				
				\begin{proposition}\label{thm:adjoint-with-graph-norm}
					For any $t\in \reg{\md E,\md F}$, the operator $Q_t\colon  \md E\to \md E_{\graph{t}}$ is unitary.  In particular, $t\colon \md E_{\graph{t}}\to\md F$ is an adjointable operator with adjoint $Q_tF_t^*\colon \md F\to \md E_{\graph{t}}$.
					\begin{proof}
						We know $Q_t\colon \md E\to \md E_{\graph{t}}$ is surjective, so it suffices to check it is an isometry:
						\[                                                      
							\ev[\md E_\graph{t}]{Q_tx}{Q_ty} = \ev[\md E]{Q_tx}{Q_ty}+\ev[\md F]{F_tx}{F_ty} = \ev[\md E]{x}{(Q_t^2+F_t^*F_t)y} = \ev[\md E]{x}{y} ,
						\]
						where in the last step we used \Cref{thm:properties-bounded-transform}.\ref{item3:thm:properties-bounded-transform}. 															
						Finally, observe that $t=F_tQ^{-1}_t$ (\Cref{thm:properties-bounded-transform}\ref{item2:thm:properties-bounded-transform}) whence $t\colon \md E_{\graph{t}}\to\md F$ is adjointable with the desired adjoint.
					\end{proof}
				\end{proposition}
				
				\begin{remark}
					We emphasize that $Q_t$ is self-adjoint as an operator in $\adj{\md E}$ but unitary as an operator in $\adj{\md E,\md E_{\graph{t}}}$. To avoid confusion, we do not write $Q_t^*$ for $Q_{t}^{-1}$ when $Q_t$ is in $\adj{\md E,\md E_{\graph{t}}}$, and $Q_t$ will be in $\adj{\md E}$, unless stated otherwise. 
				\end{remark}

		\subsubsection{Polar decomposition}
				
				Recall that an adjointable operator $T\in \adj{\md E,\md F}$ is called \emph{polar decomposable} if there exists a \emph{polar decomposition} $T = V |T|$, where $V \in \adj{\md E,\md F}$ is a partial isometry for which $\ker{V} = \ker{T}$, $\ker{V^*} = \ker{T^*}$, $\ran{V} = \overline{\ran{T}}$, and $\ran{V^*} = \overline{\ran{|T|}}$. Instead of this definition, we will mostly make use of the following characterization of polar decomposability, cf.\ \cite[Propositions 3.7 and 3.8]{Lance1995} and \cite[Proposition 15.3.7]{WeggeOlsen1993}. 
				\begin{lemma}\label{thm:polar-decomp}
					An adjointable operator $T\in \adj{\md E,\md F}$ is polar decomposable if and only if the closures of $\ran{T^*}$ and $\ran{T}$ are orthogonally complementable in $\md E$ and $\md F$, respectively. 
				\end{lemma}
				
				\begin{definition}\label{def:polar-decomposition}
					A regular operator $t\in \reg{\md E,\md F}$ is called \emph{polar decomposable} if its bounded transform $F_t := t(1+t^*t)^{-\frac12} \in \adj{\md E,\md F}$ is polar decomposable.
				\end{definition}
				
				As in the adjointable case, a regular operator with polar decomposition induces a partial isometry, and we refer the reader to \cite{FrankSharifi2010} for this approach. Instead, we will use the following characterization, which follows from \Cref{thm:invariance-kernel-rank-bounded-transform} and \cite[Propositions 3.7 and 3.8]{Lance1995} (or \cite[Proposition 15.3.7]{WeggeOlsen1993}).
				
				\begin{proposition}\label{thm:polar-decomposition}
					The operator $t\in \reg{\md E,\md F}$ is polar decomposable if and only if
				\begin{align*}
					\md E&=\ker {F_t}\oplus \overline{\ran{F_t^*}}=\ker t\oplus \overline{\ran {t^*}},\\
					\md F&=\ker {F_t^*}\oplus \overline{\ran{F_t}}=\ker {t^*}\oplus \overline{\ran t}.
				\end{align*}
				\end{proposition}

	\section{Complexes of Hilbert C*-modules}\label{sec:complexes_hilbert_modules}
				In this section, we consider complexes of Hilbert C*-modules as defined in \Cref{subsec:C*-hilbert-complexes}. Given such a C*-Hilbert complex, we construct in \Cref{subsec:bounded-transform,subsec:adjoint-complex,subsec:graph-norm-complex} its bounded transform, its adjoint complex, and its graph-norm complex. Subsequently, we define the Dirac and Laplace operators of a finite-length complex in \Cref{subsec:dirac-operator,subsec:laplace-operator}. Finally, we introduce \emph{quasi}complexes (with adjointable differential maps) in \Cref{subsec:quasicomplex-hilbert-modules}. 
			
		\subsection{The bounded transform  complex}\label{subsec:bounded-transform}
		
				\begin{definition}\label{def:adjointable-C*-hilbert-complexes}
					An \emph{adjointable} $\alg A$-Hilbert complex is an $\alg A$-Hilbert complex whose differential maps are all adjointable operators.
				\end{definition}
				
				Let $\complex{t}$ be an $\alg A$-Hilbert complex. We claim that $\complex{t}$ induces an adjointable $\alg A$-Hilbert complex via the bounded transform. Indeed, consider for each $k$ the bounded transform of $t_k$ given by 
				\begin{equation}\label{eq1:subsec:bounded-transform}
					F_{t_k}=t_kQ_{t_k}, 
				\end{equation}
				where $Q_{t_k}=(1_{\md E_k}+t_k^*t_k)^{-1/2}$. To ease the notation, we often write $Q_k \equiv Q_{t_k}$ and $Q_{k^*} \equiv Q_{t_k^*}$. As a consequence of \Cref{thm:invariance-kernel-rank-bounded-transform},  we immediately obtain the following result.
				\begin{corollary}\label{thm:equivalence-complex-property-bounded-transform}
					A sequence  $\{t_k\in \reg{\md E_k,\md E_{k+1}}\}_{k\in \mathbb Z}$ is an $\alg A$-Hilbert complex if and only if $\{F_{t_k}\in \adj{\md E_{k},\md E_{k+1}}\}_{k\in\mathbb Z}$ is an adjointable $\alg A$-Hilbert complex.
				\end{corollary}				
				
				\begin{definition}\label{def:bounded-transform-complex}
					The bounded transform of an $\alg A$-Hilbert complex $\complex t$ is the adjointable $\alg A$-Hilbert complex $\complex{F_t}$ whose differential maps are given by \eqref{eq1:subsec:bounded-transform}.
				\end{definition}																											
				For future convenience of reference, we collect several properties of the operators $F_{t_k}$ and $Q_k$ in the following technical lemma. 
				
				\begin{lemma}\label{thm:properties-resolvent-sequence}
					Let $\complex{t}$ be an $\alg A$-Hilbert complex. 
					For every $k\in \mathbb Z$, the following statements hold.
					\begin{enumerate}[(i)]
						\item $Q_k^2Q_{(k-1)^*}^2=1_{\md E_k}-F_{t_k}^*F_{t_k}-F_{t_{k-1}}F_{t_{k-1}}^*=Q_{(k-1)^*}^2Q_{k}^2$. \label{item1:thm:properties-resolvent-sequence}
						\item The operators $Q_{(k-1)^*}$ and $Q_k$ commute. \label{item2:thm:properties-resolvent-sequence}
						\item $F_{t_k}Q_{(k-1)^*}=F_{t_k}=Q_{{k+1}}F_{t_k}$. \label{item3:thm:properties-resolvent-sequence}
						\item $\norm{1_{\md E_k}-Q_k^2Q_{(k-1)^*}^2}\leq 1$. \label{item4:thm:properties-resolvent-sequence}
						\item $\ran{Q_kQ_{(k-1)^*}}=\dom{t_k}\cap \dom{t_{k-1}^*}$. \label{item5:thm:properties-resolvent-sequence}
						\item $\ran{Q_k^2Q_{(k-1)^*}^2}=\dom{t_k^*t_k}\cap \dom{t_{k-1}t_{k-1}^*}$. \label{item6:thm:properties-resolvent-sequence}
						\item The operator $Q_{k}^2Q_{(k-1)^*}^2$ has dense range. \label{item7:thm:properties-resolvent-sequence}
				    \end{enumerate}
				    \begin{proof}
						Fix $k\in \mathbb Z$. By \Cref{thm:properties-bounded-transform}\ref{item3:thm:properties-bounded-transform} and the complex property of $\complex{F_t}$, a straightforward computation shows (i). Moreover, since $\adj{\md E_k}$ is a C*-algebra, the operator $Q_{k}$ commutes with any operator that commutes with $Q_{k}^2$, see, e.g., \cite[Proposition 4.4.8]{Pedersen1989} and then use Gelfand–Naimark theorem. A similar statement holds for $Q_{(k-1)^*}$. This proves (ii).
						
						Again by \Cref{thm:properties-bounded-transform}\ref{item3:thm:properties-bounded-transform} and the complex property of $\complex{F_t}$, one can verify that $F_{t_k}Q_{(k-1)^*}^2=F_{t_k}=Q_{k+1}^2F_{t_k}$. Thus, $F_{t_k}p(Q_{(k-1)^*}^2)=p(1)F_{t_k}=p(Q_{k+1}^2)F_{t_k}$ for any polynomial $p\in \mathbb C[z]$. By continuity, it also holds for any function $p\in C([0,1])$, so in particular for the square root function. This proves (iii).
						
						We recall the following result. Let $a,b,$ and $c$ be elements of a C*-algebra $\alg B$. If $0\leq a\leq b$, then $\norm{a}\leq \norm{b}$ and $0\leq c^*ac\leq c^*bc$, see, e.g., \cite[Theorem 2.2.5]{Murphy1990}. In our case, by definition, $0\leq Q_k^2 \leq 1_{\md E_k}$ and $0\leq Q_{(k-1)^*}^2\leq 1_{\md E_k}$. Thus,
						\begin{equation*}
							0\leq Q_{k}Q_{(k-1)^*}^2Q_k\leq Q_{k}^2\leq 1_{\md E_k}.
						\end{equation*}
						Hence, $0\leq 1_{\md E_k}-Q_k^2Q_{(k-1)^*}^2\leq 1_{\md E_k}$ because $Q_k$ and $Q_{(k-1)^*}$ commute. This proves (iv).
																							
						Recall that $Q_k\colon  \md E_k\to \dom{t_k}$ and $Q_{(k-1)^*}\colon \md E\to \dom{t_{k-1}^*}$ are bijective (\Cref{thm:properties-bounded-transform}\ref{item4:thm:properties-bounded-transform}). Moreover, since $Q_k$ and $Q_{(k-1)^*}$ commute,
						\begin{align*}
							\ran{Q_kQ_{(k-1)^*}}\subset \dom{t_k}\cap \dom{t_{k-1}^*}\quad\text{and}\quad Q_k(\dom{t_{k-1}^*})\subset \dom{t_{k-1}^*}.
						\end{align*}
                        To prove the reverse inclusion, take $x\in \dom{t_k}\cap \dom{t_{k-1}^*}$. Thus, there exists an element $y\in \md E_k$ such that $Q_ky=x$. From the complex property of $\complex{F_t}$ and \Cref{thm:invariance-kernel-rank-bounded-transform}, we deduce that $\ran{F_{t_k}^*}\subset \dom{t_{k-1}^*}$. Computing $Q_k^2y$ with \Cref{thm:properties-bounded-transform}\ref{item3:thm:properties-bounded-transform}, we obtain the following relation:
                        \begin{equation*}
						 	y=Q_kx+F_{t_k}^*F_{t_k}y\in \dom{t_{k-1}^*}.
                        \end{equation*}
                        Hence, there exists an element $z\in \md E_k$ such that $Q_{(k-1)^*}z=y$. All together, $x=Q_kQ_{(k-1)^*}z$, as desired. This proves (v). Moreover, since $\ker{t_{k-1}^*} = \ker{t_{k-1}t_{k-1}^*}$, a similar argument works for (vi).
						
						Finally, if $f\colon X\to Y$ is a continuous function between topological spaces and $A\subset X$, then $\overline{f(A)}=\overline{f(\overline A)}$. Therefore, (vii) holds because $Q_{k}^2$ and $Q_{(k-1)^*}^2$ have dense range.
					\end{proof}
				\end{lemma}

        \subsection{The adjoint complex} \label{subsec:adjoint-complex}
				Given a regular operator $t\in \reg{\md E,\md F}$, the adjoint operator $t^*\in \reg{\md F,\md E}$ is also regular. Analogously, we want to derive an \emph{adjoint complex} from any $\alg A$-Hilbert complex, and one should expect to obtain a chain (instead of a cochain) complex; by duality, we can extend the previous theory to chain complexes. In the remainder, we will usually not distinguish between chain and cochain complexes.
						
				\begin{proposition}\label{thm:adjoint-complex}
					Let $\complex{t}$ be an $\alg A$-Hilbert complex. The sequence of adjoint operators $\{t_k^*\in \reg{\md E_{k+1},\md E_k}\}_{k\in \mathbb Z}$ is a (chain) $\alg A$-Hilbert complex. 
					\begin{proof}
						We only need to prove the complex property. Fix $k\in \mathbb Z$. By \Cref{thm:equivalence-complex-property-bounded-transform}, $\complex{F_t}$ is a complex. Since $F_{t_{k}}^*F_{t_{k+1}}^*=(F_{t_{k+1}}F_{t_k})^*=0$, we get that $\{F_{t_k}^*\}_{k\in \mathbb Z}$ is a (chain) $\alg A$-Hilbert complex. Therefore, the result follows by  \Cref{thm:equivalence-complex-property-bounded-transform} because $F_{t_k}^*=F_{t_k^*}$ for every $k\in \mathbb Z$.
					\end{proof}
				\end{proposition}
																		
				\begin{definition}\label{def:adjoint-complex}
					The adjoint of $\complex{t}$ is the (chain) $\alg A$-Hilbert complex $\complex{t^\sharp}$ whose differential maps are $\{t_k^*\in \reg{\md E_{k+1},\md E_{k}}\}_{k\in \mathbb Z}$.
				\end{definition}

                \begin{remark}
                    Observe that the adjoint of complexes also reverses the sequence of Hilbert C*-modules. According to \cite[p.90]{BruningLesch1992}, a more descriptive notation for the adjoint of $\complex{t}$ could be $\complex[\md E^*]{t^*}$. However, we do not follow this notation.
                \end{remark}
				
				It follows immediately that the adjoint of complexes defines an involution (i.e., if $\complex{t}$ is an $\alg A$-Hilbert complex, then the adjoint of $\complex{t^\sharp}$ is $\complex{t}$), and that the complex $\complex{F_t^\sharp}$ is the bounded transform of $\complex{t^\sharp}$. Moreover, a sequence  of regular operators is an $\alg A$-Hilbert complex if and only if its sequence of adjoint operators is an $\alg A$-Hilbert complex.
				
        \subsection{The graph-norm complex}\label{subsec:graph-norm-complex}
				As an alternative to the bounded transform complex, we can use the \emph{graph inner-product} to transfer the situation to the adjointable case. 
				Consider an $\alg A$-Hilbert complex $\complex{t}$. 
				Fix $k\in \mathbb Z$ and recall from \Cref{sec:graph-inner-product} that $\md E_{\graph{t_k}} := \dom{t_k}$ is a Hilbert $\alg A$-module endowed with the graph inner-product. 
				Take $x\in \md E_{\graph{t_k}}$ and $y\in \md E_{\graph{t_{k+1}}}$. 
				By the complex property of $\complex{t}$ and \Cref{thm:adjoint-with-graph-norm} applied to the operator $t_{k}$, we obtain 
				\begin{equation}\label{eq1:subsec:graph-norm-complex}
					\ev[\md E_{\graph{t_{k+1}}}]{t_kx}{y}=\ev[\md E_{k+1}]{t_kx}{y}=\ev[\md E_{\graph{t_k}}]{x}{Q_kF_{t_k}^*y}.
				\end{equation}
				Hence, the operator $t_k\colon  \md E_{\graph{t_k}}\to \md E_{\graph{t_{k+1}}}$ is adjointable for every $k\in \mathbb Z$. Furthermore, this sequence of operators clearly still satisfies the complex property.
				
				\begin{definition}\label{def:graph-norm-complex}
					The \emph{graph-norm complex} of $\complex{t}$ is the adjointable $\alg A$-Hilbert complex $\complex[\md E_{\graph{t}}]{t}$ whose differential maps are $\{t_k\in \adj{\md E_{\graph{t_k}},\md E_{\graph{t_{k+1}}}}\}_{k\in\mathbb Z}$.
				\end{definition}
				
				The following result establishes the relation between the bounded transform complex and the graph-norm complex. 
				
				\begin{proposition}\label{thm:isomorphism_graph_and_boundedtransform}
					Let $\complex{t}$ be an $\alg A$-Hilbert complex with bounded transform $\complex{F_t}$ and graph-norm complex $\complex[\md E_{\graph{t}}]{t}$. Then the diagram
					\[\begin{tikzcd}
						\cdots & {\md E_k} & {\md E_{k+1}} & \cdots \\
						\cdots & {\md E_{\graph{t_k}}} & {\md E_{\graph{t_{k+1}}}} & \cdots
						\arrow[from=1-1, to=1-2]
						\arrow["{F_{t_k}}", from=1-2, to=1-3]
						\arrow["{Q_{t_k}}"', from=1-2, to=2-2]
						\arrow[from=1-3, to=1-4]
						\arrow["{Q_{t_{k+1}}}", from=1-3, to=2-3]
						\arrow[from=2-1, to=2-2]
						\arrow["{t_k}", from=2-2, to=2-3]
						\arrow[from=2-3, to=2-4]
					\end{tikzcd}\]
					commutes for every $k\in \mathbb Z$. 
					\begin{proof}
						Fix $k\in \mathbb Z$. 
						Since $t_k=F_{t_k}Q_k^{-1}$ (see \Cref{thm:properties-bounded-transform}\ref{item2:thm:properties-bounded-transform}) and $Q_{k+1}F_{t_k}=F_{t_k}$ (see \Cref{thm:properties-resolvent-sequence}), we observe that \(t_k = Q_{k+1} F_{t_k} Q_k^{-1}\), as desired. 
					\end{proof}
				\end{proposition}					

        \subsection{The Dirac operator}\label{subsec:dirac-operator}
				Our main aim in this subsection is to prove the regularity of the (even) Dirac operator described in \Cref{sec:main:Dirac}. 
				Our methods are different from the Hilbert space case \cite[Section 2]{BruningLesch1992}, and we start by considering the Dirac operator of the bounded transform complex. 
				
				Let us fix a finite-length $\alg A$-Hilbert complex $\complex{t}:=\{t_k\in \reg{\md E_k,\md E_{k+1}}\}_{k=0}^N$ with bounded transform $\complex{F_t}$. 
				Recall the notation $\md E_{\text{ev}}:=\bigoplus_{k} \md E_{2k}$ and $\md E_{\text{odd}}:=\bigoplus_{k} \md E_{2k+1}$. 
				We now use \Cref{thm:properties-resolvent-sequence} to obtain an adjointable operator by \airquotes{rolling up} the complex $\complex{F_t}$.
				
				\begin{theorem}\label{thm:bounded-even-dirac-operator}
					The operator
					\begin{align*}
						\evDirac{F_t}\colon \md E_{\mathrm{ev}}&\to \md E_{\mathrm{odd}} , & 
						(x_{2k})_{k}&\mapsto (F_{t_{2k}}x_{2k}+F_{t_{2k+1}}^*x_{2k+2})_{k}
					\end{align*}
					is adjointable in $\adj{\md E_{\mathrm{ev}},\md E_{\mathrm{odd}}}$, and we denote its adjoint by $\oddDirac{F_t}$. Moreover, there exists a regular operator $\evDirac{t}\in \reg{\md E_{\mathrm{ev}},\md E_{\mathrm{odd}}}$ such that $\evDirac{F_t}$ is its bounded transform.
					\begin{proof}
						The operator $\evDirac{F_t}$ is clearly adjointable with adjoint defined as follows:
						\begin{equation*}
							\oddDirac{F_t}(y_{2k+1})_{k}:=(F_{t_{2k}}^*y_{2k+1}+F_{t_{2k-1}}y_{2k-1})_{k},\quad \forall (y_{2k+1})_{k}\in \md E_{\mathrm{odd}}.
						\end{equation*}
						To prove the existence of a regular operator $\evDirac{t}$ whose bounded transform is $\evDirac{F_t}$, we need to show that $\norm{\evDirac{F_t}}\leq 1$ and the operator $1_{\md E_{\mathrm{ev}}}-\oddDirac{F_t}\evDirac{F_t}$ has dense range, see, e.g.,  \cite[Theorem 10.4]{Lance1995}. Using \Cref{thm:properties-resolvent-sequence}\ref{item1:thm:properties-resolvent-sequence}, we derive that
						\begin{equation}\label{eq2:subsec:dirac-operator}
							1_{\md E_{\mathrm{ev}}}-\oddDirac{F_t}\evDirac{F_t}=\bigoplus_{k} Q_{2k}^2Q_{(2k-1)^*}^2.
						\end{equation}
						Thus, by \Cref{thm:properties-resolvent-sequence}\ref{item7:thm:properties-resolvent-sequence}, $1_{\md E_{\mathrm{ev}}}-\oddDirac{F_t}\evDirac{F_t}$ has dense range, and \Cref{thm:properties-resolvent-sequence}\ref{item4:thm:properties-resolvent-sequence} gives
						\begin{equation*}
							\norm{\evDirac{F_t}}^2=\norm{\oddDirac{F_t}\evDirac{F_t}}\leq \sup_{k}\norm{1_{\md E_{2k}}-Q_{2k}^2Q_{(2k-1)^*}^2}\leq 1.\qedhere
						\end{equation*}
					\end{proof}
				\end{theorem}
																				
				\begin{definition}\label{def:dirac-operator}
					The \emph{Dirac operator} of $\complex{t}$ is the regular self-adjoint operator
					\begin{equation*}
						\Dirac{t}:=\begin{pmatrix}
							0&\oddDirac{t}\\
							\evDirac{t}&0
						\end{pmatrix}\in \reg{\md E_{\mathrm{ev}}\oplus \md E_{\mathrm{odd}}},
					\end{equation*}
					where $\evDirac{t}$ is the regular operator obtained from \Cref{thm:bounded-even-dirac-operator}, and $\oddDirac{t}$ denotes the adjoint of $\evDirac{t}$. We call $\evDirac{t}$ the even part of the Dirac operator of $\complex{t}$, or the \emph{even Dirac operator} for short, and $\oddDirac{t}$ the odd part.
				\end{definition}
				
				The following theorem shows that the operator $\evDirac{t}$ obtained as the `unbounded transform' of $\evDirac{F_t}$ agrees with the operator defined in \Cref{def:even-odd-Dirac}. Here the direct sum of the intersected domains below should be interpreted as (finite) orthogonal sums of pre-Hilbert C*-modules.
				
				\begin{theorem}\label{thm:dirac-operator}
					We have that
					\begin{equation*}
						\dom{\evDirac{t}}=\bigoplus_{k}\dom{t_{2k}}\cap \dom{t_{2k-1}^*}\subset \md E_{\mathrm{ev}}
					\end{equation*}
					and
					\begin{equation*}
						\dom{\Dirac{t}^-}=\bigoplus_{k}\dom{t_{2k}^*}\cap \dom{t_{2k+1}}\subset \md E_{\mathrm{odd}}.
					\end{equation*}
					Moreover, the operators $\evDirac{t}\colon \dom{\evDirac{t}}\to \md E_{\mathrm{odd}}$ and $\oddDirac{t}\colon  \dom{\oddDirac{t}}\to \md E_{\mathrm{ev}}$ have the following matrix forms:
					\begin{equation}\label{eq3:subsec:dirac-operator}
						\evDirac{t}=\begin{pmatrix}
							t_0&t_1^*&0&\cdots\\
							0&t_2&t_3^*&\cdots\\
							0&0&t_4&\cdots\\
							\vdots&\vdots&\vdots&\ddots
						\end{pmatrix}
						\quad\text{and}\quad
						\oddDirac{t}=\begin{pmatrix}
							t_0^*&0&0&\cdots\\
							t_1&t_2^*&0&\cdots\\
							0&t_3&t_4^*&\cdots\\
							\vdots&\vdots&\vdots&\ddots
						\end{pmatrix}.
					\end{equation}
					\begin{proof}
						By \Cref{thm:properties-resolvent-sequence}, $\ran{Q_{(2k-1)^*}Q_{2k}}=\dom{t_{2k}}\cap\dom{t_{2k-1}^*}$ for every $k$. Moreover, since $\ran{(1_{\md E_{\mathrm{ev}}}-\oddDirac{F_t}\evDirac{F_t})^{1/2}}=\dom{\evDirac{t}}$, then \eqref{eq2:subsec:dirac-operator} yields the identity:
						\begin{equation*}
							\dom{\evDirac{t}}=\bigoplus_{k}\dom{t_{2k}}\cap\dom{t_{2k-1}^*}.
						\end{equation*}
						From \Cref{thm:properties-resolvent-sequence}\ref{item2:thm:properties-resolvent-sequence}-\ref{item3:thm:properties-resolvent-sequence}, we have for each $k$ that 
						\begin{align*}                                                      
							Q_{(2k-1)^*}^{-1}Q_{2k}^{-1} &= Q_{2k}^{-1}Q_{(2k-1)^*}^{-1} , & 
							F_{2k}Q_{(2k-1)^*}^{-1} &\subset F_{2k} , & 
							F_{2k-1}^*Q_{2k}^{-1} &\subset F_{2k-1}^* .
						\end{align*}
						Furthermore, with \Cref{thm:properties-bounded-transform}\ref{item2:thm:properties-bounded-transform} and \eqref{eq2:subsec:dirac-operator}, a straightforward computation shows that
						\begin{equation}\label{eq:Dirac-bounded-transform}
							\evDirac{t}=\evDirac{F_t}\left(\bigoplus_{k} Q_{2k}^{-1}Q_{(2k-1)^*}^{-1}\right)=\begin{pmatrix}
								t_0&t_1^*&0&\cdots\\
								0&t_2&t_3^*&\cdots\\
								0&0&t_4&\cdots\\
								\vdots&\vdots&\vdots&\ddots
							\end{pmatrix}.
						\end{equation}
						A similar argument works for $\Dirac{t}^-$. The result follows.	
					\end{proof}
				\end{theorem}
				
				Since $\evDirac{F_t}$ is the bounded transform of $\evDirac{t}$, we may write $\evDirac{F_t} = \evDirac{t} Q_{\evDirac{t}}$. We then see from \eqref{eq:Dirac-bounded-transform} (and a similar computation for $\oddDirac{t}$) that 
				\begin{align}\label{eq:Q-Dirac}
					Q_{\evDirac{t}} &= \bigoplus_{k} Q_{2k} Q_{(2k-1)^*} , & 
					Q_{\oddDirac{t}} &= \bigoplus_{k} Q_{2k^*} Q_{2k-1} .
				\end{align}
				
				The explicit description of the Dirac operator given in \eqref{eq3:subsec:dirac-operator} allows us to obtain the Dirac operator of the bounded transform complex, the adjoint complex, and the graph-norm complex. This is the content of the following corollaries.
				
				\begin{corollary}\label{thm:dirac-operator-bounded-complex}
					The Dirac operator of $\complex{F_t}$ is
					\begin{equation*}
						\Dirac{F_t}:=\begin{pmatrix}
							0&\oddDirac{F_t}\\
							\evDirac{F_t}&0
						\end{pmatrix}\in \adj{\md E_{\mathrm{ev}}\oplus \md E_{\mathrm{odd}}}.
					\end{equation*}
					Moreover, $\Dirac{F_t}$ is the bounded transform of $\Dirac{t}$.
				\end{corollary}
				
				The previous corollary makes evident that no information is lost by going from a complex to its Dirac operator. Indeed, assume that $\evDirac{t}$ is the even Dirac operator of a finite-length $\alg A$-Hilbert complex $\complex{t}$. By construction, we know that $\evDirac{t}$ is a regular operator and its bounded transform, $\evDirac{F_t}$, is the even Dirac operator of the bounded transform complex $\complex{F_t}$. Using projections and inclusions, we can extract every differential map $F_{t_k}$ from $\evDirac{F_t}$. This allows us to reconstruct each differential map of $\complex{t}$ from the operators $\{F_{t_{k}}\}_{k}$ by \Cref{thm:properties-bounded-transform}.
				
				\begin{corollary}\label{thm:dirac-operator-adjoint-complex}
					The even Dirac operator of the adjoint complex $\complex{t^\sharp}$ (viewed as a \emph{cochain} complex) is
					\begin{equation*}
						\evDirac{t^\sharp}:=\begin{pmatrix}
							t_{N}^*&t_{N-1}&0&\cdots\\
							0&t_{N-2}^*&t_{N-3}&\cdots\\
							0&0&t_{N-4}^*&\cdots\\
							\vdots&\vdots&\vdots&\ddots
						\end{pmatrix}.
					\end{equation*}
					Moreover, up to unitary transformations, we obtain that
					\begin{equation}\label{eq4:subsec:dirac-operator}
						\evDirac{t^\sharp} \simeq \begin{cases}
							\oddDirac{t},& N\text{ even},\\
							\evDirac{t},& N \text{ odd}.
						\end{cases}
					\end{equation}
					\begin{proof}
						It follows from direct computations by using \eqref{eq3:subsec:dirac-operator}.
					\end{proof}
				\end{corollary}
				
				Finally, we investigate the relation between $\Dirac{t_{\Gamma}}$ and $\Dirac{F_t}$, where $\Dirac{t_\Gamma}$ is the Dirac operator of $\complex[\md E_{\graph{t}}]{t}$. Recall from \Cref{thm:adjoint-with-graph-norm} that each $Q_{k}\colon \md E_k\to \md E_{\graph{t_k}}$ is unitary. Let us denote by $Q_{\mathrm{ev}}$ the operator $\bigoplus_{k} Q_{2k}\colon \md E_{\mathrm{ev}}\to \bigoplus_{k} \md E_{\graph{t_{2k}}}$, and we define $Q_{\mathrm{odd}}$ similarly.
				
				\begin{corollary}\label{thm:dirac-operator-graph-complex}
					The Dirac operator of $\complex[\md E_{\graph{t}}]{t}$ satisfies that
					\begin{equation*}
						\Dirac{t_\Gamma}\begin{pmatrix}
							Q_{\mathrm{ev}}&0\\
							0&Q_{\mathrm{odd}}
						\end{pmatrix}=\begin{pmatrix}
							Q_{\mathrm{ev}}&0\\
							0&Q_{\mathrm{odd}}
						\end{pmatrix}\Dirac{F_t}.
					\end{equation*}
					\begin{proof}
						Recall that $Q_{t_k}F_{t_k}^*\colon \md E_{\graph{t_{k+1}}}\to\md E_{\graph{t_k}}$ is the adjoint of $t_k\colon \md E_{\graph{t_k}}\to \md E_{\graph{t_{k+1}}}$. We then obtain the Dirac operator of $\complex[\md E_{\graph{t}}]{t}$ using \eqref{eq3:subsec:dirac-operator}. By \Cref{thm:properties-bounded-transform} and \Cref{thm:properties-resolvent-sequence}, $Q_{\mathrm{odd}}\evDirac{F_t}=\evDirac{t_\Gamma}Q_{\mathrm{ev}}$. Since $Q_k$ is unitary in $\adj{\md E_k,\md E_{\graph{t_k}}}$, taking adjoints yields the identity $\oddDirac{F_t}Q_{\mathrm{odd}}^{-1}=Q_{\mathrm{ev}}^{-1}\oddDirac{t_\Gamma}$, and the result follows.
					\end{proof}
				\end{corollary}

        \subsection{The Laplace operator}\label{subsec:laplace-operator}
				We now define the \emph{Laplace operator} of a C*-Hilbert complex and investigate some of its properties. While it is closely related to the Dirac operator, the Laplace operator is also an important piece of information on its own, and it will play a crucial role in \Cref{sec:fredholm-equivalences}. Let us fix a finite-length $\alg A$-Hilbert complex $\complex{t}:=\{t_k\in \reg{\md E_k,\md E_{k+1}}\}_{k=0}^N$ with bounded transform $\complex{F_t}$. 
				
				\begin{definition}\label{def:laplace-operator}
					The \emph{Laplace operator} of $\complex{t}$ is given by 
					\[                                                   
						\Laplace{t} := \bigoplus_{k} \Laplace{t}_k \colon \bigoplus_{k} \md E_k \to \bigoplus_{k} \md E_k , 
					\]
					where the $k$-th Laplace operator of $\complex{t}$ is defined as 
					\begin{equation*}
						\Laplace{t}_k:=t_k^*t_k+t_{k-1}t_{k-1}^* \colon (\dom{t_k^*t_k}\cap \dom{t_{k-1}t_{k-1}^*})\subset \md E_{k}\to \md E_{k} , \quad 0\leq k\leq N+1.
					\end{equation*}
					We also define
					\begin{equation*}
						\Laplace{t}_{\mathrm{ev}}:=\bigoplus_{k} \Laplace{t}_{2k} \quad \text{and}\quad \Laplace{t}_{\mathrm{odd}}:=\bigoplus_{k} \Laplace{t}_{2k+1}.
					\end{equation*}
				\end{definition}							
									
				While we focus our attention on a finite-length complex, we note that the ($k$-th) Laplace operator is also well defined without the finite-length assumption.
				
				\begin{theorem}[{\cite[Theorem 1.1]{LeschMesland2019}}]\label{thm:LeschMesland}
					Let $t,s\in \reg{\md E}$ be self-adjoint and regular operators. Assume that $t$ and $s$ satisfy the following statements.
					\begin{enumerate}[(i)]
						\item There exist positive constants $C_0,C_1$, and $C_2$ such that the inequality
						\begin{equation*}
							\ev{(st+ts)x}{(st+ts)x}\leq C_0\ev{x}{x}+C_1\ev{sx}{sx}+C_2\ev{tx}{tx}
						\end{equation*}
						holds for all $x\in \mathcal F(s,t):=\{x\in \dom t\cap \dom s\mid sx\in \dom{t}, tx\in\dom s\}$.
						\item There exists a core $\mathscr D$ of $t$ such that $(s+i\lambda)^{-1}(\mathscr{D})\subset \mathcal F(s,t)$, for $\lambda\in \mathbb R\setminus [-\lambda_0,\lambda_0]$.
					\end{enumerate}
					Then $(s+t)$ is regular and self-adjoint in $\reg{\md E}$ with domain $\dom{t}\cap \dom{s}$.
				\end{theorem}
				
				\begin{theorem}\label{thm:regularity-kth-laplace}
					The Laplace operator $\Laplace{t}$ of a finite-length $\alg A$-Hilbert complex $\complex{t}$ is self-adjoint and regular. 
					\begin{proof}
						It suffices to show that, for any $k$, the sum $\Laplace{t}_k$ is regular and self-adjoint. By regularity of $t_k$ and $t_{k-1}^*$, the operators $t:=t_k^*t_k$ and $s:=t_{k-1}t_{k-1}^*$ are self-adjoint in $\reg{\md E_k}$. We now use \Cref{thm:LeschMesland} and follow its notation. Observe that $\ran{t}\subset \ker s$ and $\ran s\subset \ker t$. Thus, $\mathcal F(s,t)=\dom s\cap \dom t$ and $st=ts=0$ on $\mathcal F(s,t)$. Hence, condition (i) is trivially satisfied. We now show (ii). Since the spectrum of $s$ is real, the operator $(s+i\lambda)\colon  \dom{s}\to \md E_k$ is bijective for $\lambda\in \mathbb R\setminus \{0\}$ \cite[Proposition 5.20]{Kustermans1997}.  Put $\mathscr{D}:=\dom{t}$ and take $x\in \mathscr{D}$. By construction, $y:=(s+i\lambda )^{-1}x$ belongs to $\dom{s}$. Since $\ran{s}\subset \ker t$, we get $y=-i\lambda^{-1}(x-sy)\in \dom t$. In conclusion, $(s+i\lambda)^{-1}x$ belongs to $\mathcal F(s,t)$. The result follows.
					\end{proof}
				\end{theorem}
				
				\begin{corollary}\label{thm:dirac-operator-and-laplace-operator}
					For a finite-length $\alg A$-Hilbert complex $\complex{t}$, we have 
					\begin{equation*}
						\Laplace{t}\simeq\Laplace{t}_{\mathrm{ev}}\oplus \Laplace{t}_{\mathrm{odd}}=\Dirac{t}^-\evDirac{t}\oplus \evDirac{t}\Dirac{t}^-.
					\end{equation*}
					\begin{proof}
                        Using \eqref{eq3:subsec:dirac-operator}, one computes that $\oddDirac{t}\evDirac{t} \subset \Laplace{t}_{\mathrm{ev}}$ and $\evDirac{t}\oddDirac{t} \subset \Laplace{t}_{\mathrm{odd}}$. Since these operators are all self-adjoint, we must have equalities. 
					\end{proof}
				\end{corollary}
				
				As before, we can use this explicit description to obtain the Laplace operator of the bounded transform complex, the adjoint complex, and the graph-norm complex. First, we investigate the relation between the Laplace operator of $\complex{F_t}$ and of $\complex{t}$. Unfortunately, the bounded transform of $\Laplace{t}_k$ is not $\Laplace{F_t}_k$. In fact, one can prove that
				\begin{equation*}
					F_{\Laplace{t}_k}=F_{t_k^*t_k}+F_{t_{k-1}t_{k-1}^*}, \quad 0\leq k\leq N+1.
				\end{equation*}
				However, we do not need these relations and the following corollary will be enough.
				
				\begin{corollary}\label{thm:laplace-operator-bounded-complex}
					Let $\Delta^{\complex{F_t}}$ be the Laplace operator of $\complex{F_t}$. Then
					\begin{equation*}
						\Laplace{F_t}_{\mathrm{ev}}=\Laplace{t}_{\mathrm{ev}}\left(\bigoplus_{k} Q_{2k}^2Q_{(2k-1)^*}^2\right)\quad\text{and}\quad \Laplace{F_t}_{\mathrm{odd}}=\Laplace{t}_{\mathrm{odd}}\left(\bigoplus_{k} Q_{(2k)^*}^2Q_{2k+1}^2\right).
					\end{equation*}
					\begin{proof}
						For any $t\in \reg{\md E,\md F}$, we have $F_t^*F_t=t^*Q_{t^*}tQ_t=t^*tQ_t^2$ by \Cref{thm:properties-bounded-transform} items \ref{item1:thm:properties-bounded-transform} and \ref{item7:thm:properties-bounded-transform}. 
                        Applying this to $\evDirac{t}$, we obtain 
                        \[
                        \oddDirac{F_t}\evDirac{F_t} 
                        = F^*_{\evDirac{t}} F_{\evDirac{t}} 
                        = \oddDirac{t} \evDirac{t} Q_{\evDirac{t}}^2 
                        = \oddDirac{t} \evDirac{t} \left( \bigoplus_k Q_{2k}^2Q_{(2k-1)^*}^2 \right) ,
                        \]
                        where in the last step we used \eqref{eq:Q-Dirac} and \Cref{thm:properties-resolvent-sequence}\ref{item2:thm:properties-resolvent-sequence}. 
                        The first equality then follows from \Cref{thm:dirac-operator-and-laplace-operator}. 
                        The second equality is similar. 
					\end{proof}
				\end{corollary}
				
				\begin{corollary}\label{thm:laplace-operator-adjoint-complex}
					The $k$-th Laplace operator of $\complex{t^\sharp}$ satisfies the following:
						\begin{equation*}
							\Laplace{t^\sharp}_{k}=t_{N+1-k}^*t_{N+1-k}+t_{N-k}t_{N-k}^*=\Laplace{t}_{N+1-k}, \quad 0\leq k\leq N+1.
						\end{equation*}
				\end{corollary}
				
				\begin{corollary}\label{thm:laplace_operator_graph_complex}
					The Laplace operator of $\complex[\md E_{\graph{t}}]{t}$ satisfies that
					\begin{equation*}
						\Laplace[\md E_{\graph{t}}]{t}=\begin{pmatrix}
							Q_{\mathrm{ev}}&0\\
							0&Q_{\mathrm{odd}}
						\end{pmatrix}\Laplace{F_t} \begin{pmatrix}
						Q_{\mathrm{ev}}^{-1}&0\\
						0&Q_{\mathrm{odd}}^{-1}
						\end{pmatrix},
					\end{equation*}
					where $Q_{\mathrm{ev}}$ and $Q_{\mathrm{odd}}$ are the unitary operators defined for \Cref{thm:dirac-operator-graph-complex}.
				\end{corollary}

        \subsection{C*-Hilbert quasicomplexes}\label{subsec:quasicomplex-hilbert-modules}
				As described in \Cref{sec:main:Fredholm-property}, we define the Fredholm property of an $\alg A$-Hilbert complex without referring to its cohomology groups. 
				Therefore, the complex property is not essential (at least, in the adjointable case), and it can be replaced with the following notion (attributed to Putinar \cite{Putinar1982} in the context of Banach spaces). 
				
				\begin{definition}\label{def:C*-hilbert-quasicomplexes}
					A quasicomplex of Hilbert $\alg A$-modules, or an \emph{$\alg A$-Hilbert quasicomplex} for short, is a sequence $\complex{S}:=\{S_k\in \adj{\md E_{k},\md E_{k+1}}\}_{k\in\mathbb Z}$ of adjointable operators, defined on Hilbert $\alg A$-modules, such that
					\begin{equation}\label{eq1:subsec:quasicomplex-hilbert-modules}
						S_{k+1}S_{k}\in \comp{\md E_k,\md E_{k+2}}, \quad \forall k\in\mathbb Z.
					\end{equation}
				\end{definition}
				
				Observe that we are using the same notation for complexes and quasicomplexes but no confusion should arise in the future. Let us set some notation. For adjointable operators $T,S\in \adj{\md E,\md F}$, we define the following equivalence relation:
				\begin{equation*}
					T\sim S \quad \text{ if and only if }\quad T-S\in\comp{\md E,\md F}.
				\end{equation*}
				Since the compact operators form an ideal, we note that this equivalence is preserved when taking sums or products. We write $S_{k+1}S_k\sim 0$ instead of \eqref{eq1:subsec:quasicomplex-hilbert-modules}, and we refer to it as the \emph{quasicomplex property} of $\complex{S}$. 
				
				Analogous to C*-Hilbert complexes, one can easily define \emph{finite-length quasicomplexes}, \emph{the adjoint quasicomplex}, and \emph{chain quasicomplexes}. On this note, observe that an adjointable operator can be regarded as a finite-length complex and as a finite-length quasicomplex. More generally, an adjointable $\alg A$-Hilbert complex is a quasicomplex. For this reason, we would like to extend the results from \Cref{subsec:dirac-operator,subsec:laplace-operator} to quasicomplexes. This is the content of the following theorem, which follows straightforwardly from the definitions and the quasicomplex property of $\complex{S}$.
				
				\begin{theorem}\label{thm:laplace-dirac-operators-for-quasicomplex}
					Let $\complex{S}:=\{S_k\in \adj{\md E_{k},\md E_{k+1}}\}_{k=0}^N$ be a finite-length $\alg A$-Hilbert quasicomplex. For every $0\leq k\leq N+1$, the $k$-th Laplace operator of $\complex{S}$, denoted by $\Laplace{S}_k$, is well defined by \Cref{def:laplace-operator}. Similarly,  the even and odd parts of the Dirac operator of $\complex{S}$, denoted by $\evDirac{S}$ and $\oddDirac{S}$, are well defined by \eqref{eq3:subsec:dirac-operator}. 
					
					Moreover, $\oddDirac{S}$ is the adjoint of $\evDirac{S}$ and 
					\begin{equation*}
						\Laplace{S}:=\Laplace{S}_{\mathrm{ev}}\oplus\Laplace{S}_{\mathrm{odd}}\sim \oddDirac{S}\evDirac{S}\oplus \evDirac{S}\oddDirac{S}.
					\end{equation*}
				\end{theorem}			

	\section{Hodge decompositions}\label{sec:weak_hodge_decomposition}		
							
				In this section, we completely characterize C*-Hilbert complexes with weak or strong Hodge decomposition. Such complexes behave very similarly to the Hilbert-space case. In general, however, depending on the choice of C*-algebra $\alg A$, we note that not every $\alg A$-Hilbert complex has (weak or strong) Hodge decomposition. 
								
        \subsection{Weak Hodge decomposition} \label{subsec:weak-hodge-decomposition}

				In this subsection, we first consider the notion of weak Hodge decomposition as defined in \Cref{def:hodge-decomposition}. For a fixed $k\in \mathbb Z$, we say that an $\alg A$-Hilbert complex $\complex{t}$ has \emph{weak Hodge decomposition at $k$} whenever
				\begin{equation}
					\md E_k=(\ker{t_k}\cap \ker{t_{k-1}^*})\oplus \overline{\ran{t_k^*}}\oplus \overline{\ran{t_{k-1}}}.
				\end{equation}
				The following statement gives an alternative description of the first summand in this weak Hodge decomposition. 
				\begin{lemma}\label{thm:kernel-kth-laplace-operator}
					Let $\complex{t}$ be an $\alg A$-Hilbert complex. Then
					\begin{equation*}
						\ker{\Laplace{t}_k}=\ker {t_{k}}\cap \ker {t_{k-1}^*},\quad \forall k\in \mathbb Z.
					\end{equation*}
					\begin{proof}
						Fix $k\in \mathbb Z$. The inclusion $ \ker{t_k}\cap \ker{t_{k-1}^*}\subset \ker{\Laplace{t}_k}$ follows immediately from the definition of $\Laplace{t}_k$. Conversely, consider $x\in \ker{\Laplace{t}_k}$ to get that
						\begin{equation*}
							0=\ev{x}{\Laplace{t}_kx}=\ev{x}{t_{k}^*t_kx+t_{k-1}t_{k-1}^*x}=\ev{t_kx}{t_kx}+\ev{t_{k-1}^*x}{t_{k-1}^*x},
						\end{equation*}
						which implies both $t_kx=0$ and $t_{k-1}^*x=0$. (Observe that this argument does not use completeness of the Hilbert C*-modules or the complex property of the sequence.)
					\end{proof}
				\end{lemma}				
				
				The even Dirac operator and the Laplace operator will play a crucial role to determine if a complex has weak Hodge decomposition. For this reason, we investigate some of their properties involving their range and kernel. Notice that the following theorem is a general version of \cite[Theorem 3.2 and Proposition 3.7]{Lance1995}.
				
				\begin{theorem}\label{thm:kernel-range-laplace-operator}
					Let $\complex{t}$ be a finite-length $\alg A$-Hilbert complex. Then
					\begin{enumerate}[(i)]
					\item \label{item1:thm:kernel-range-laplace-operator}
						$\ker{\Laplace{t}_{\mathrm{ev}}}=\ker{\Laplace{F_t}_\mathrm{ev}}$ and $\ran{\Laplace{t}_{\mathrm{ev}}}=\ran{\Laplace{F_t}_{\mathrm{ev}}}$; 
					\item \label{item2:thm:kernel-range-laplace-operator}
						$\ker{\evDirac{t}}=\ker{\Laplace{t}_{\mathrm{ev}}}$ and $\ker{\oddDirac{t}}=\ker{\Laplace{t}_{\mathrm{odd}}}$;
					\item \label{item3:thm:kernel-range-laplace-operator}
						$\overline{\ran{\evDirac{t}}}=\overline{\ran{\Laplace{t}_{\mathrm{odd}}}}$ and $\overline{\ran{\oddDirac{t}}}=\overline{\ran{\Laplace{t}_{\mathrm{ev}}}}$. 
					\end{enumerate}
					Furthermore, the operators $\evDirac{t}$, $\oddDirac{t}$, $\Laplace{t}_{\mathrm{ev}}$, and $\Laplace{t}_\mathrm{odd}$ all have closed range if and only if any of them has closed range, and in this case 
					\begin{equation}\label{eq3:subsec:weak-hodge-decomposition}
						\begin{aligned}
						\md{E}_{\mathrm{ev}}&=\ker{\evDirac{t}}\oplus \ran{\oddDirac{t}}=\ker{\Laplace{t}_{\mathrm{ev}}}\oplus \ran{\Laplace{t}_{\mathrm{ev}}},\\
						\md{E}_{\mathrm{odd}}&=\ker{\oddDirac{t}}\oplus \ran{\evDirac{t}}=\ker{\Laplace{t}_{\mathrm{odd}}}\oplus \ran{\Laplace{t}_{\mathrm{odd}}}.\
					\end{aligned}
					\end{equation}
					\begin{proof}
						By \Cref{thm:dirac-operator-and-laplace-operator}, we know that $\Laplace{t}_{\mathrm{ev}}=\oddDirac{t}\evDirac{t}$ and $\Laplace{F_t}_{\mathrm{ev}}=\oddDirac{F_t}\evDirac{F_t}$, where $\evDirac{F_t}$ is the bounded transform of $\evDirac{t}$.
						Properties \ref{item1:thm:kernel-range-laplace-operator}-\ref{item3:thm:kernel-range-laplace-operator} then follow from \Cref{thm:invariance-kernel-rank-bounded-transform}. 
						The final statement follows from \Cref{thm:closed_range}. 
					\end{proof}
				\end{theorem}		
				
				The following statement follows straightforwardly from the observation that the Laplace operator is given by a finite orthogonal direct sum of $k$-th Laplace operators. 
				\begin{lemma}\label{thm:range-direc-sum-operators}
					Let $\complex{t}$ be a finite-length $\alg A$-Hilbert complex. The Laplace operator $\Laplace{t}$ has polar decomposition (or closed range) if and only if the $k$-th Laplace operator $\Laplace{t}_k$ has polar decomposition (resp.\ closed range) for every $k$.
				\end{lemma}
				
				We now relate the polar decomposition of the $k$-th Laplace operator of a complex with its weak Hodge decomposition at the step $k$.

				\begin{proposition}\label{thm:weak-hodge-decomposition}
					Let $\complex{t}$ be an $\alg A$-Hilbert complex, such that the $k$-th Laplace operator $\Laplace{t}_k$ of $\complex{t}$ is polar decomposable for some $k\in \mathbb Z$. Then the complex $\complex{t}$ has weak Hodge decomposition at $k$ and 
					\begin{align*}
						\ker{t_k}&=\ker{\Laplace{t}_k}\oplus \overline{\ran{t_{k-1}}},\\
						\ker{t_{k-1}^*}&=\ker{\Laplace{t}_k}\oplus \overline{\ran{t_{k}^*}},\\
						\overline{\ran{\Laplace{t}_k}}&=\overline{\ran{t_k ^*}}\oplus\overline{\ran{t_{k-1}}}.
					\end{align*}
					\begin{proof}
						Fix $k\in\mathbb Z$. From \Cref{thm:kernel-kth-laplace-operator}, we have that
						\begin{equation}\label{eq4:subsec:weak-hodge-decomposition}
							\ker{\Laplace{t}_k}=\ker{t_k}\cap\ker{t_{k-1}^*}.
						\end{equation}																										
						By definition, $\Laplace{t}_k=t_k^*t_k+t_{k-1}t_{k-1}^*$. Hence, $\ran{\Laplace{t}_k}\subset\ran{t_k^*}\oplus_{\mathrm{alg}}\ran{t_{k-1}}$. In particular, 
						\begin{equation*}
							\overline{\ran{\Laplace{t}_k}}\subset\overline{\ran{t_k^*}}\oplus\overline{\ran{t_{k-1}}}.
						\end{equation*}
						We now show the converse inclusion. By \Cref{thm:kert_perp_rant^*} and \eqref{eq4:subsec:weak-hodge-decomposition},
						\begin{align}
							\overline{\ran {t_{k-1}}}&\subset(\ker{t_{k-1}^*})^\perp\subset(\ker{\Laplace{t}_k})^\perp, \label{eq5:subsec:weak-hodge-decomposition}\\
							\overline{\ran {t^*_k}}&\subset(\ker{t_k})^\perp\subset(\ker{\Laplace{t}_k})^\perp.\label{eq6:subsec:weak-hodge-decomposition}
						\end{align}
						If we assume that $\Laplace{t}_k$ has a polar decomposition, so that $\md E_k=\ker {\Laplace{t}_k}\oplus \overline{\ran {\Laplace{t}_k}}$ by \Cref{thm:polar-decomposition}, then $(\ker{\Laplace{t}_k})^\perp=\overline{\ran{\Laplace{t}_k}}$. Hence, \eqref{eq6:subsec:weak-hodge-decomposition} and \eqref{eq5:subsec:weak-hodge-decomposition} yield the inclusion
						\begin{equation*}
							\overline{\ran{t^*_k}}\oplus	\overline{\ran{t_{k-1}}}\subset \overline{\ran{\Laplace{t}_k}}.
						\end{equation*} 
						Therefore, $\complex{t}$ has weak Hodge decomposition at $k$.
						
						Finally, notice that the inclusion $\ker{\Laplace{t}_k}\oplus\overline{\ran{t_{k-1}}}\subset \ker {t_k}$ follows by the complex property of $\complex{t}$ and \eqref{eq4:subsec:weak-hodge-decomposition}. The reverse inclusion follows from the weak Hodge decomposition and the inclusion $\overline{\ran {t_{k}^*}}\subset (\ker {t_k})^\perp$.
                        The proof of the remaining relation is similar.
					\end{proof}
				\end{proposition}

				\begin{remark}
					The proof in the previous theorem is a modified version of \cite[Theorem 11]{Krysl2015}. However, \cite{Krysl2015} establishes the strong Hodge decomposition of what we would call an \emph{adjointable complex of pre-Hilbert C*-modules} under the hypothesis the author calls \emph{self-adjoint parametrix possessing}.
				\end{remark}

				We are in position to characterize complexes with weak Hodge decomposition. 

				\begin{theorem}[Weak Hodge decomposition]\label{thm:weak-hodge-decomposition-equivalences}
					Let $\complex{t}$ be an $\alg A$-Hilbert complex. The following statements are equivalent:
					\begin{enumerate}[(i)]
						\item The complex $\complex{t}$ has weak Hodge decomposition.
						\item  For every $k \in \mathbb Z$, the differential map $t_k$ is polar decomposable.
						\item For every $k \in \mathbb Z$, the $k$-th Laplace operator $\Laplace{t}_k$ is polar decomposable.
					\end{enumerate}
					If $\complex{t}$ has finite length, the above statements are furthermore equivalent to: 
					\begin{enumerate}[(i)]
						\addtocounter{enumi}{3}
						\item The even Dirac operator $\evDirac{t}$ is polar decomposable.
					\end{enumerate}
					\begin{proof}
						In view of \Cref{thm:polar-decomposition}, we know that $\evDirac{t}$ is polar decomposable if and only if $\oddDirac{t}$ is polar decomposable. In consequence, the equivalence of (iii) with (iv) is an immediate application of \Cref{thm:kernel-range-laplace-operator}\ref{item3:thm:kernel-range-laplace-operator} and \Cref{thm:range-direc-sum-operators}.

						(i)$\Rightarrow$(ii) Using the weak Hodge decomposition, we have that each $t_k$ is polar decomposable because the closure of the ranges of $t_k$ and $t_k^*$ are orthogonally complementable. 
						
						(iii)$\Rightarrow$(i) It follows immediately by \Cref{thm:weak-hodge-decomposition}.
						
						(ii)$\Rightarrow$(iii) Fix $k$. By \eqref{eq4:subsec:weak-hodge-decomposition} and since $\Laplace{t}_k$ is self-adjoint, it suffices to show that 
						\begin{equation*}
							\overline{\ran{\Laplace{t}_k}}=\overline{\ran{t_k^*}}\oplus\overline{\ran{t_{k-1}}}
						\end{equation*}
						From definition of $\Laplace{t}_k$, $\overline{\ran{\Laplace{t}_k}}\subset \overline{\ran{t_k^*}}\oplus\overline{\ran{t_{k-1}}}$. We now show the reverse inclusion. By \Cref{thm:kernel-range-laplace-operator}, we have that
						\begin{equation*}															\overline{\ran{t_k^*t_k}}=\overline{\ran{t_k^*}}\quad\text{and}\quad\overline{\ran{t_{k-1}t_{k-1}^*}}=\overline{\ran{t_{k-1}}}.
						\end{equation*}
						Consider an element $z\oplus w\in \overline{\ran{t_k^*}}\oplus\overline{\ran{t_{k-1}}}$. There exist sequences $\{x_j\}_{j\geq 0}$ and $\{y_j\}_{j\geq 0}$ such that $\lim_j t_k^*t_kx_j=z$ and $\lim_jt_{k-1}t_{k-1}^*y_j=w$. For each $j$, put
						\begin{align*}
							x_j&:=x_j^{(0)}\oplus x_j^{(1)}\in \ker{t_k}\oplus\overline{\ran{t_k^*}},\\
							y_j&:=y_j^{(0)}\oplus y_j^{(1)}\in \ker{t_{k-1}^*}\oplus\overline{\ran{t_{k-1}}}.
						\end{align*}
						Notice that $t_k^*t_kx_j^{(1)}=t_k^*t_kx_j$ and $t_{k-1}t_{k-1}^*y_j^{(1)}=t_{k-1}t_{k-1}^*y_j$ whence $\lim_j t_k^*t_kx_j^{(1)}=z$ and $\lim_j t_{k-1}t_{k-1}^*y_j^{(1)}=w$. Moreover, by the complex property, $\overline{\ran{t_{k-1}}}\subset\ker {t_k}$ and $\overline{\ran{t_k^*}}\subset \ker {t_{k-1}^*}$. Thus $t_ky_j^{(1)}=0=t_{k-1}^*x_j^{(1)}$. Therefore,
						\begin{equation*}
							\lim_j\Laplace{t}_k(x_j^{(1)}+y_j^{(1)})= z\oplus w \in \overline{\ran{\Laplace{t}_k}}.
						\end{equation*}
						Hence, $\overline{\ran{\Laplace{t}_k}}=\overline{\ran{t_k^*}}\oplus\overline{\ran{t_{k-1}}}$, as desired.
					\end{proof}
				\end{theorem}

				As a consequence of the Hilbert projection theorem, Hilbert complexes always have weak Hodge decomposition, see, e.g., \cite[Lemma 2.1]{BruningLesch1992}. The next corollary shows examples of C*-algebras, other than the complex numbers, for which any C*-Hilbert complex has weak Hodge decomposition.
					
				\begin{corollary}\label{thm:weak-hodge-decomposition-subalgebra-compact-operators}
					Suppose the C*-algebra $\alg B$ admits a faithful representation $\alg B \hookrightarrow \alg K(\md H)$ as compact operators on some Hilbert space $\md H$. 
					Then every $\alg B$-Hilbert complex $\complex{t}$ has weak Hodge decomposition.
					\begin{proof}
						In \cite{Magajna1997}, the author shows that every closed Hilbert $\alg B$-module is orthogonally complementable. Hence, each differential map is polar decomposable by \Cref{thm:polar-decomposition}.
					\end{proof}
				\end{corollary}	

				In the next corollary, we investigate the weak Hodge decomposition of the adjoint complex, the bounded transform complex, and the graph-norm complex. Observe that we do not include the finite-length assumption because the proof does not use the Dirac operator.

				\begin{corollary}\label{thm:weak-hodge-decomposition-derived-complexes}
					Let $\complex{t}$ be an $\alg A$-Hilbert complex. All of the following objects have weak Hodge decomposition if and only if any of them has weak Hodge decomposition.
					\begin{enumerate}[(i)]
						\item The complex $\complex{t}$.
						\item The adjoint complex $\complex{t^\sharp}$.
						\item The bounded transform complex $\complex{F_t}$.
						\item The graph-norm complex $\complex[\md E_{\graph{t}}]{t}$.														
					\end{enumerate}
					\begin{proof}
						The equivalence of (i) with (ii) and (iii) is an immediate consequence of the definition and \Cref{thm:invariance-kernel-rank-bounded-transform}, respectively. To prove the equivalence  of (iii) with (iv), observe that the equation
						\begin{equation*}
							\Laplace[\md E_{\graph{t}}]{t}_k=Q_{k}\Laplace{F_t}_kQ_{k}^{-1},\quad k\in \mathbb Z,
						\end{equation*}
						follows by similar computations to those in \Cref{thm:laplace_operator_graph_complex}, where $Q_k\colon \md E_k\to \md E_{\graph{t_k}}$ is unitary; although the complex does not have finite-length, we can carry out the computations at each step using that $Q_kF_{t_{k-1}}=F_{t_{k-1}}$ and $F_{t_{k-1}}^*Q_{k}^{-1}\subset F_{t_{k-1}}^*$ by \Cref{thm:properties-resolvent-sequence}. Fix $k\in \mathbb Z$. Since $\Laplace{F_t}_k$ and $\Laplace[\md E_{\graph{t}}]{t}_k$ are unitarily equivalent, we deduce that
						\begin{equation*}
							\ker{\Laplace[\md E_{\graph{t}}]{t}_k} = Q_k(\ker{\Laplace{F_t}_k}) \quad\text{and}\quad \ran{\Laplace[\md E_{\graph{t}}]{t}_k}=Q_k(\ran{\Laplace{F_t}_k}).
						\end{equation*}
						Moreover,
						\begin{equation*}
							\overline{\ran{\Laplace[\md E_{\graph{t}}]{t}_k}} = Q_k(\overline{\ran{\Laplace{F_t}_k}})\quad\text{and}\quad Q_k^{-1}(\overline{\ran{\Laplace[\md E_{\graph{t}}]{t}_k}})=\overline{\ran{\Laplace{F_t}_k}}
						\end{equation*}
						because $Q_k$ is a homeomorphism. Now, by \Cref{thm:properties-bounded-transform}\ref{item3:thm:properties-bounded-transform}, we conclude that $Q_k|_{\ker{F_{t_k}}}=1_{\ker{F_{t_k}}}$ and $Q_k^{-1}|_{\ker{F_{t_k}}}=1_{\ker{F_{t_k}}}$. Thus, if $\Laplace{F_t}_k$ is polar decomposable and we take $x\in \md E_{\graph{t_k}}$, then there exist an element $y_1\oplus y_2\in \ker{\Laplace{F_t}_k}\oplus\overline{\ran{\Laplace{F_t}_k}}$ such that
						\begin{equation*}
							x=Q_k(y_1\oplus y_2)=y_1\oplus Q_ky_2\in \ker{\Laplace[\md E_{\graph{t}}]{t}}\oplus \overline{\ran{\Laplace[\md E_{\graph{t}}]{t}}}.
						\end{equation*}
						The converse implication is similar. The result readily follows by \Cref{thm:weak-hodge-decomposition-equivalences}.
					\end{proof}
				\end{corollary}

        \subsection{Strong Hodge decomposition}\label{subsec:strong_hodge_decomposition}

				Recall from \Cref{def:hodge-decomposition} that an $\alg A$-Hilbert complex $\complex{t}$ has strong Hodge decomposition if
					\begin{equation*}
						\md E_k=(\ker {t_k}\cap \ker {t_{k-1}^*})\oplus {\ran{t_k^*}}\oplus {\ran{t_{k-1}}}, \quad \forall k\in\mathbb Z.
					\end{equation*}
				
				\begin{theorem}[Strong Hodge decomposition]\label{thm:hodge-decomposition-equivalences}
					Let $\complex{t}$ be an $\alg A$-Hilbert complex. The following statements are equivalent:
					\begin{enumerate}[(i)]
						\item The complex $\complex{t}$ has strong Hodge decomposition.
						\item  For every $k\in\mathbb Z$, the differential map $t_k$ has closed range.
						\item For every $k\in\mathbb Z$, the $k$-th Laplace operator $\Laplace{t}_k$ has closed range.
					\end{enumerate}
					If $\complex{t}$ has finite length, the above statements are furthermore equivalent to: 
					\begin{enumerate}[(i)]
						\addtocounter{enumi}{3}
						\item The even Dirac operator $\evDirac{t}$ has closed range.
					\end{enumerate}
					If any of these conditions is satisfied, then, for every $k$, we have
					\begin{align*}
						\ker{t_k}&=\ker{\Laplace{t}_k}\oplus \ran{t_{k-1}},\\
						\ker{t_{k-1}^*}&=\ker{\Laplace{t}_k}\oplus \ran{t_{k}^*},\\
						\ran{\Laplace{t}_k}&=\ran{t_k ^*}\oplus\ran{t_{k-1}}.
					\end{align*}
					\begin{proof}
						It follows by \Cref{thm:kernel-range-laplace-operator} that $\evDirac{t}$ has closed range if and only if $\oddDirac{t}$ has closed range. Thus, the equivalence of (iii) with (iv) is an immediate consequence of \Cref{thm:kernel-range-laplace-operator} and \Cref{thm:range-direc-sum-operators}.
						
						(i)$\Leftrightarrow$(ii) Using the strong Hodge decomposition, we have that each $t_k$ has closed range because their ranges are orthogonally complementable, see, e.g.,  \cite[Lemma 15.3.4]{WeggeOlsen1993}. 
						
						 (iii)$\Rightarrow$(i) This is a direct modification of \Cref{thm:weak-hodge-decomposition}. Using \eqref{eq3:subsec:weak-hodge-decomposition}, the reader should simply remove the \airquotes{overline} and replace \airquotes{polar decomposition} by \airquotes{closed range}.
						
						(ii)$\Rightarrow$(iii) This is also a direct modification of the analogous implication in \Cref{thm:weak-hodge-decomposition-equivalences}. In fact, again by \eqref{eq3:subsec:weak-hodge-decomposition}, the reader should simply remove the \airquotes{overline} and consider constant sequences $\{x:=x_j\}_{j\geq 0}$ and $\{y:=y_j\}_{j\geq 0}$.
						
						The last statement follows by its counterpart in \Cref{thm:weak-hodge-decomposition}.
					\end{proof}
				\end{theorem}							

				\begin{corollary}\label{thm:hodge-decomposition-derived_complexes}
					Let $\complex{t}$ be an $\alg A$-Hilbert complex. All of the following objects have strong Hodge decomposition if and only if any of them has strong Hodge decomposition.
					\begin{enumerate}[(i)]\setlength{\itemsep}{0pt}
						\item The complex $\complex{t}$.
						\item The adjoint complex $\complex{t^\sharp}$.
						\item The bounded transform $\complex{F_t}$.
						\item The graph-norm complex $\complex[\md E_{\graph{t}}]{t}$.														
					\end{enumerate}
					\begin{proof}
						It is a simple modification of \Cref{thm:weak-hodge-decomposition-derived-complexes}.
					\end{proof}
				\end{corollary}

        \subsection{Cohomology groups}\label{subsec:topological-complexes}
				If $\complex{t}$ is an $\alg A$-Hilbert complex, then $\ran{t_{k-1}}$ neither need to be orthogonally complementable nor a closed subspace of $\ker{t_k}$. Thus, the \emph{cohomology groups} $\{\ker{t_k}/\ran{t_{k-1}}\}_{k\in \mathbb Z}$ are not necessarily Hilbert $\alg A$-modules. The following theorem shows that for complexes with strong Hodge decomposition the cohomology groups are well defined.

				\begin{theorem}[{\cite[Corollary 14]{Krysl2015}}]\label{thm:self-adj-poss-thm-Kry}
					Let $\complex{T}$ be an adjointable $\alg A$-Hilbert complex with strong Hodge decomposition. For each $k\in\mathbb Z$, the \emph{$k$-th cohomology group} 
					\begin{equation*}
						\coH{k}{\complex{T}}:=\ker {T_k}/\ran{T_{k-1}}
					\end{equation*}
					is well defined and isomorphic, as Hilbert $\alg A$-modules, to $\ker{\Laplace{T}_k}$.
				\end{theorem}

				\begin{remark}
						In \cite{Krysl2015}, Kr\'ysl works with the hypothesis \airquotes{self-adjoint parametrix possessing.} In our context, an adjointable complex is self-adjoint parametrix possessing if and only if it has strong Hodge decomposition, see \cite[Remark 5]{Krysl2016}.
				\end{remark}
				
				As a consequence of \Cref{thm:invariance-kernel-rank-bounded-transform}, \Cref{thm:kernel-kth-laplace-operator}, and \Cref{thm:kernel-range-laplace-operator,thm:self-adj-poss-thm-Kry}, we obtain the following result.
				\begin{corollary}\label{coro:cohomology-kernel}
					Let $\complex{t}$ be an $\alg A$-Hilbert complex with strong Hodge decomposition. Then we have the following isomorphisms of Hilbert $\alg A$-modules:
					\begin{equation*}
						\coH{k}{\complex{t}} := \ker{t_k}/\ran{t_{k-1}} \cong \ker{\Laplace{t}_k} \cong \ker{t_k} \cap \ker{t_{k-1}^*} , \quad\forall k\in\mathbb Z.
					\end{equation*}
				\end{corollary}
				
				\begin{remark}                                                
					If each differential map in an $\alg A$-Hilbert complex has closed range, then the complex is sometimes called a \emph{topological} complex, referring to the fact that $\coH{k}{\complex{t}}\cong \ker{\Laplace{t}_k}$ is a \emph{topological} isomorphism of Hilbert $\alg A$-modules (cf.\ \cite[p.138]{SchulzeTarkhanov1998}). By \Cref{thm:hodge-decomposition-equivalences}, an $\alg A$-Hilbert complex is topological if and only if it has strong Hodge decomposition, and we will only be using the latter terminology. 
				\end{remark}

	\section{Maps of complexes}\label{sec:maps_complexes}														
				In this section, we briefly discuss results that we shall need in \Cref{sec:fredholm_property} regarding maps of complexes. Particularly, we focus on complexes with weak Hodge decomposition; in this case, the arguments are very similar to the Hilbert-space case, and we present the statements without proofs (further details can however be found in \cite[Section 2.3]{Villegas-Villalpando2024}).
		
				\begin{definition}
					Let $\complex{t}$ and $\complex[\md F]{s}$ be $\alg A$-Hilbert complexes. A collection of $\alg A$-module maps $\{g_k\colon \md E_k\to \md F_k\}_{k\in\mathbb Z}$ will be called a \emph{map of complexes} if, for each $k\in \mathbb Z$, the inclusion $g_k(\dom{t_k})\subset \dom{s_k}$ holds and the diagram
					\[\begin{tikzcd}
						\cdots & {\dom {t_k}} & {\dom {t_{k+1}}} & \cdots \\
						\cdots & {\dom {s_k}} & {\dom {s_{k+1}}} & \cdots
						\arrow[from=1-1, to=1-2]
						\arrow["{t_k}", from=1-2, to=1-3]
						\arrow[from=2-1, to=2-2]
						\arrow["{s_k}", from=2-2, to=2-3]
						\arrow[from=1-3, to=1-4]
						\arrow[from=2-3, to=2-4]
						\arrow["{g_k}"', from=1-2, to=2-2]
						\arrow["{g_{k+1}}", from=1-3, to=2-3]
					\end{tikzcd}\]
					commutes. We might assume that each $g_k$ is bounded, unless stated otherwise, and we use the symbol $g\colon \complex{t}\to \complex[\md F]{s}$ to denote the collection of maps $\{g_k\}_{k\in \mathbb Z}$.
				\end{definition}
				
				If $g\colon \complex{t}\to \complex[\md F]{s}$ is a map of complexes, then it is straightforward to show that $g_k(\ran{t_{k-1}})\subset\ran{s_{k-1}}$ and $g_k(\ker{t_k})\subset\ker{s_k}$ for every $k\in \mathbb Z$.
				
				\begin{definition}
					An \emph{adjointable map of complexes} $G\colon \complex{t}\to \complex[\md F]{s}$ is a map of complexes  whose collection of maps $\{G_k\}_{k\in \mathbb Z}$ consists of adjointable operators.
				\end{definition}
							
				\begin{lemma}
					Let $G\colon  \complex{t}\to \complex[\md F]{s}$ be an adjointable map of complexes. Then $G$ induces a map of complexes $\{G_k^*\}_{k\in \mathbb Z}$ denoted by $G^*\colon \complex[\md F]{s^\sharp}\to \complex{t^\sharp}$.
					\begin{proof}
						Let $k\in\mathbb Z$. One can show that $G_k^*(\dom{s_{k-1}^*})$ is contained in $\dom{t_{k-1}^*}$. Indeed, for any $x\in \dom{s_{k-1}^*}$, since $G_kt_{k-1}\subset s_{k-1}G_{k-1}$, we get that
						\begin{equation*}
							\ev{t_{k-1}y}{G_k^*x}=\ev{y}{G_{k-1}^*s_{k-1}^*x},\quad \forall y\in \dom{t_{k-1}}.
						\end{equation*}
						Then $G$ induces an adjointable map of complexes $ G^*\colon \complex[\md F]{s^\sharp}\to\complex{t^\sharp}$.
					\end{proof}
				\end{lemma}
				
				\begin{example}\label{example:isomorphism_graph_and_boundedtransform}
					Let $\complex{t}$ be an $\alg A$-Hilbert complex with bounded transform $\complex{F_t}$ and graph-norm complex $\complex[\md E_{\graph{t}}]{t}$. Then 
					the maps $\{Q_{t_k}\in \adj{\md E_k,\md E_{\graph{t_k}}}\}_{k\in \mathbb Z}$ yield a \emph{unitary} map of complexes $Q\colon \complex{F_t}\to \complex[\md E_{\graph{t}}]{t}$ (see \Cref{thm:isomorphism_graph_and_boundedtransform}).
				\end{example}
								
				The following lemma and its corollary will be briefly used to consider the \emph{weak index} of \emph{weakly Fredholm complexes} in \Cref{subsec:fredholm_complex_with_WHD}.

				\begin{lemma}[{\cite[Lemma 2.7]{BruningLesch1992}}]\label{thm:functiorial_homomorphism_map_complexes}
					Let $g\colon \complex{t}\to\complex[\md F]{s}$ be a map of complexes. Denote by $\Laplace{t}_k$ and $\Laplace[\md F]{s}_k$ the $k$-th Laplace operators of $\complex{t}$ and $\complex[\md F]{s}$, respectively. Suppose that both complexes have weak Hodge decomposition. Consider the orthogonal projections $\hat{\mathcal P}_k$, $\mathcal P_k^{(1)}$, and $\mathcal P_k^{(2)}$ from $\md F_k$ onto $\ker{\Laplace[\md F]{s}_k}$, $\overline{\ran{s_{k-1}}}$, and $\overline{\ran{s_{k}^*}}$, respectively. Then $g$ induces functorial homomorphisms for all $k\in \mathbb Z$, namely
					\begin{align*}
						\hat g_k:=\hat{\mathcal P}_k g_k|_{\ker{\Laplace{t}_k}}&\colon \ker{\Laplace{t}_k}\to\ker{\Laplace[\md F]{s}_k},\\
						g_k^{(1)}:=\mathcal P_k^{(1)}g_k|_{\overline{\ran{t_{k-1}}}}&\colon \overline{\ran{t_{k-1}}}\to\overline{\ran{s_{k-1}}},\\
						g_k^{(2)}:=\mathcal P_k^{(2)}g_k|_{\overline{\ran{t_{k}^*}}}&\colon \overline{\ran{t_{k}^*}}\to\overline{\ran{s_{k}^*}}.
					\end{align*}
					If we further assume that $g$ is unitary, then so are $\hat g$, $g_k^{(1)}$ and $g_k^{(2)}$.
				\end{lemma}
				
				\begin{corollary}\label{thm:isomorphic_regular_complexes}
					Let $\complex{t}$ and $\complex[\md F]{s}$ be $\alg A$-Hilbert complexes having weak Hodge decomposition. If there exist $h\colon \complex{t}\to\complex[\md F]{s}$ and $g\colon \complex[\md F]{s}\to\complex{t}$ maps of complexes such that $gh=1_{\complex{t}}$ and $hg=1_{\complex[\md F]{s}}$, then
					\begin{equation*}
						\ker{\Laplace{t}_k}\cong\ker{\Laplace[\md F]{s}_k},\quad\forall k\in\mathbb Z,
					\end{equation*} 
					as $\alg A$-modules. If $G:=g$ and $H:=h$ are adjointable, then we have an isomorphism of Hilbert C*-modules.
				\end{corollary}
	
				We now consider the cohomology maps of a complex when its cohomology groups are well defined as Hilbert C*-modules (\Cref{subsec:topological-complexes}). 
				
				\begin{theorem}
					Let $G\colon \complex{t}\to\complex[\md F]{s}$ be an adjointable map of complexes with strong Hodge decomposition. For any $k\in\mathbb Z$, the map
					\begin{equation*}
						\mathcal H_k(G)(x+\ran {t_{k-1}}):=G_kx+\ran{s_{k-1}} ,\quad \forall x\in \ker{t_k},
					\end{equation*}
					induces a collection $\mathcal H(G):=\{\mathcal H_k(G)\in \adj{\coH{k}{\complex{t}},\coH{k}{\complex{s}}}\}_{k\in\mathbb Z}$ of adjointable operators between the cohomology groups of $\complex{t}$ and $\complex[\md F]{s}$.
				\end{theorem}
				
				Let $g,h\colon  \complex{t}\to \complex[\md F]{s}$ be maps of complexes. A \emph{chain-homotopy map} for $g$ and $h$ is a collection $p:=\{p_k\colon \md E_{k+1}\to \md F_k\}_{k\in\mathbb Z}$ of $\alg A$-module maps such that, for every $k\in \mathbb Z$, the inclusion $p_{k-1}(\dom{t_k})\subset \dom{s_{k-1}}$ holds and
				\begin{equation*}
					g_k-h_k=p_{k}t_{k}+s_{k-1}p_{k-1}\text{ in }\dom{t_k}.
				\end{equation*}
				We might assume that each $p_k$ is bounded, unless stated otherwise. In this case, $g$ and $h$ are said to be \emph{chain-homotopic}.
				
				\begin{definition}
					Two complexes $\complex{t}$ and $\complex[\md F]{s}$ are called \emph{chain-homotopic equivalent} if there exist maps of complexes $g\colon \complex{t}\to \complex[\md F]{s}$ and $h\colon \complex[\md F]{s}\to\complex{t}$ such that $hg$ and $gh$ are both chain homotopic to the identity.
				\end{definition}
				
				\begin{lemma}[{\cite[Lemma 2.9]{BruningLesch1992}}]\label{thm:homology_chain_homotopic_morphisms}
					Let $\complex{t}$ and $\complex[\md F]{s}$ be $\alg A$-Hilbert complexes with strong Hodge decomposition. Assume that $G,H\colon \complex{t}\to \complex[\md F]{s}$ are adjointable maps of complexes, and that $G$ and $H$ are chain-homotopic. Then the induced maps coincide on homology, that is, $\mathcal H(G)$ and $\mathcal H(H)$ are equal. We also have $\hat G=\hat H \colon \ker{\Laplace{t}}\to\ker{\Laplace[\md F]{s}}$.
				\end{lemma}
				
				\begin{remark}
					The maps coincide on homology as $\alg{A}$-modules without any extra conditions on the complexes or assuming that $g:=G$ and $h:=H$ are adjointable. Similarly, we have assumed that the complexes have strong Hodge decomposition, but to prove that $\hat h=\hat g$ we only need weak Hodge decomposition. With the assumptions of the lemma we ensure that all modules are Hilbert and the maps are adjointable.
				\end{remark}

	\section{The Fredholm property}\label{sec:fredholm_property}
				In this section, we consider the \emph{Fredholm property} for C*-Hilbert complexes and quasicomplexes. However, we shall focus on developing tools to investigate the Fredholm property of C*-Hilbert \emph{quasi}complexes, which will later be used in \Cref{sec:fredholm-equivalences} to investigate the general case.
									
	    \subsection{C*-Fredholm complexes}\label{subsec:fredholm-complexes}
				
				Recall from \Cref{def:joint-parametrix} the notion of a joint parametrix. We recall that, for adjointable operators, the Fredholm property is characterized by the existence of a \emph{single} parametrix. 
				
				\begin{definition}[cf. {\cite[Definition 3.1]{Exel1993}}]\label{def:adjointable-fredholm-operators}
					An adjointable operator $T\in \adj{\md E,\md F}$ is called $\alg A$-Fredholm if it has a joint parametrix $(P_l,P_r)$ that satisfies the identity $P_l=P_r$. In this case, we say that $P_l\in \adj{\md F,\md E}$ is a parametrix of $T$.
				\end{definition}														
				
				However, by \cite[Definition 2.1]{Joachim2003}, we know that the unbounded case is more involved, namely a regular operator $t\in \reg{\md E,\md F}$ is called \emph{Fredholm} if it has a pseudo-left and a pseudo-right inverse in the following sense:
				\begin{quote}
					An adjointable operator $P_l\in\adj{\md F,\md E}$ is a \emph{pseudo-left inverse} of $t$ if $P_lt$ is closable and its closure is adjointable such that $1_{\md E}-\overline{P_lt}\in \comp{\md E}$.   Similarly, an adjointable operator $P_r\in \adj{\md F,\md E}$ is a \emph{pseudo-right inverse} of $t$ if $tP_r$ is closable and its closure is adjointable such that $1_{\md F}-\overline{tP_r}\in \comp{\md F}$. 
				\end{quote}
				We will prove in \Cref{thm:consistency-joachim} that a regular operator is Fredholm in the sense of \cite{Joachim2003} if and only if there exists a joint parametrix for $t$. 
				
				For the moment, we tackle the adjointable case. Due to the following theorem and its corollary, the Fredholm property of an adjointable operator, viewed as a finite-length complex, is consistent with \Cref{def:adjointable-fredholm-operators}.
				
				\begin{lemma}\label{thm:parmetrices-fredholm-operators}
					For any $T\in \adj{\md E,\md F}$, the following statements hold.
					\begin{enumerate}[(i)]
						\item If $P_1$ and $P_2$ are parametrices of $T$, then $P_1-P_2$ is a compact operator.
						\item If $(P_l,P_r)$ is a joint parametrix of $T$, then $P_l-P_r$ is a compact operator. In particular, $P_l$ and $P_r$ are parametrices of $T$. 
					\end{enumerate}
					\begin{proof}
						It suffices to show (ii). Assume that $(P_l,P_r)$ is a joint parametrix of $T$. Then there exist compact operators $C_l$ and $C_r$ such that $P_l(1_{\md F}+C_r)=P_lTP_r=(1_{\md E}+C_l)P_r$. This shows that $P_l-P_r$ is compact. The result readily follows.
					\end{proof}
				\end{lemma} 
				
				\begin{corollary}\label{thm:parametrix-joint-parametrix-operators}
					An adjointable operator has a parametrix if and only if it has a joint parametrix.
				\end{corollary}
				
				The following lemma shows that, at least for adjointable Fredholm operators, the compact operators appearing in \eqref{eq2:def:joint-parametrix-complex} can be taken as finite-rank operators. We later show in \Cref{thm:finite-rank-complex-parametrix-II} that the same holds in general.
				
				\begin{lemma}[{\cite[Lemma 4.4]{GraciaVarillyFigueroa2001}}]\label{thm: parametrix-finite-rank}
					If $T\in \adj{\md E,\md F}$ is an $\alg A$-Fredholm operator, then there exists an adjointable operator $S\in \adj{\md F,\md E}$ such that $1_{\md E}-ST$ and $1_{\md F}-TS$ are finite-rank operators.
				\end{lemma}
	
	    \subsection{C*-Fredholm quasicomplexes}\label{subsec:fredholm-quasicomplexes}
				We now define the Fredholm property for C*-Hilbert quasicomplexes by following closely the presentation given in \cite[Section 5.1]{SchulzeTarkhanov1998} for quasicomplexes of Fréchet spaces.
				
				Let $\complex{S}=\{S_k\in \adj{\md E_k,\md E_{k+1}}\}_{k=0}^N$ be an arbitrary finite-length $\alg A$-Hilbert quasicomplex. A \emph{parametrix} of $\complex{S}$ is a sequence of adjointable operators $\{P_k\in\adj{\md E_{k+1},\md E_{k}}\}_{k=0}^N$, and, of compact operators $\{C_k\in \comp{\md E_k}\}_{k=0}^{N+1}$ such that
				\begin{equation}\label{eq1:def:Fredholm_complexes}
					P_{k}S_k+S_{k-1}P_{k-1}=1_{\md E_k}-C_k,\quad 0\leq k\leq N+1.
				\end{equation}		
				Such sequence will be denoted by $\complex P:=\{P_k\}_{k=0}^N$, where  the collection $\{C_k\}_{k=0}^N$ is considered implicitly. Let us mention that a parametrix $\complex P$ of $\complex{S}$ is not necessarily a quasicomplex, but if $\complex{P}$ is a parametrix of $\complex{S}$ that is also a quasicomplex, then we refer to it as a \emph{quasicomplex parametrix}.
				
				\begin{definition}\label{def:fredholm-quasicomplex}
					A Fredholm quasicomplex of Hilbert $\alg A$-modules, or an $\alg A$-Fredholm quasicomplex for short, is a finite-length $\alg A$-Hilbert quasicomplex with a parametrix.
				\end{definition}
				
				\begin{remark}
					By \Cref{def:adjointable-fredholm-operators}, an adjointable $\alg A$-Fredholm operator, viewed as a finite-length quasicomplex, is an $\alg A$-Fredholm quasicomplex.
				\end{remark}
				
				From here on, we extensively use the equivalence of compact perturbations defined in \Cref{subsec:quasicomplex-hilbert-modules}. In particular, the relations
				\begin{equation*}
					P_{k}S_k+S_{k-1}P_{k-1}\sim 1_{\md E_k}, \quad 0\leq k\leq N+1,
				\end{equation*}
				have an equivalent meaning to \eqref{eq1:def:Fredholm_complexes}. Moreover, we say that an $\alg A$-Hilbert quasicomplex $\complex{\tilde S}$ is a compact perturbation of another $\alg A$-Hilbert quasicomplex $\complex{S}$ if $\tilde S_k\sim S_k$ for every $k\in \mathbb Z$. 
				The following result follows straightforwardly from the fact that the compact operators form a $*$-ideal, so that the equivalence $\sim$ is preserved when taking sums, products, and adjoints. 
				
				\begin{proposition}\label{thm:basic-properties-fredholm-quasicomplexes}
					Let $\complex{S}$ be an $\alg A$-Fredholm quasicomplex, and let $\complex{P}$ be a parametrix of $\complex{S}$.
					\begin{enumerate}[(i)]
						\item The adjoint quasicomplex $\complex{S^\sharp}$  is an $\alg A$-Fredholm quasicomplex with para\-metrix $\complex{P^\sharp}:=\{P_k^*\}_{k=0}^N$.\label{item1:thm:basic-properties-fredholm-quasicomplexes}
						\item If $\complex{P}$ is a \emph{quasicomplex}, then $\complex{P}$ is an $\alg A$-Fredholm quasicomplex with quasicomplex parametrix $\complex{S}$.\label{item2:thm:basic-properties-fredholm-quasicomplexes}
						\item Any compact perturbation of $\complex{S}$ is an $\alg A$-Fredholm quasicomplex, and any compact perturbation of a parametrix of $\complex{S}$ is again a parametrix of $\complex{S}$.\label{item3:thm:basic-properties-fredholm-quasicomplexes}
					\end{enumerate}
				\end{proposition}															
	
	    \subsection{Associated quasicomplexes}\label{subsec:associated_quasicomplexes}
				In the next lemma, we show how to construct a parametrix from \emph{associated maps}. This method can be extracted from the proof of \cite[Proposition 6.1]{AtiyahBott1967}, and it is of particular interest when considering the Laplace operator of a quasicomplex, cf. \cite[Theorem 5.1.4]{SchulzeTarkhanov1998}.
	
				\begin{lemma}\label{thm:associated-quasicomplex}
					Let $\complex{S}=\{S_k\in \adj{\md E_k,\md E_{k+1}}\}_{k=0}^N$ be a finite-length $\alg A$-Hilbert quasicomplex. Assume there exists a sequence $\{P_k\in \adj{\md E_{k+1},\md E_k}\}_{k=0}^N$ such that the operators
					\begin{equation*}
						D_k:=P_kS_k+S_{k-1}P_{k-1}, \quad 0\leq  k\leq N+1,
					\end{equation*} 
					are $\alg A$-Fredholm. Then $\complex{S}$ is an $\alg A$-Fredholm quasicomplex.
					\begin{proof}
						For each $k$, let $R_k$ be a parametrix of $D_k$ and define $\hat P_{k}:=R_kP_k$. We claim that the collection $\{\hat P_k\}_{k}$ is a parametrix of $\complex S$. Indeed, from the quasicomplex property, that is, $S_kS_{k-1}\sim 0$, we derive that
						\begin{equation}\label{eq1:thm:associated-quasicomplex}
							S_{k-1}D_{k-1}\sim S_{k-1}P_{k-1}S_{k-1} \sim D_{k}S_{k-1},\quad \forall k.
						\end{equation}
						Since $R_k$ is a parametrix of $D_k$, then
						\begin{equation}\label{eq2:thm:associated-quasicomplex}
							R_kD_k\sim1_{\md E_k}\quad \text{and}\quad D_kR_k\sim 1_{\md E_k},\quad \forall k.
						\end{equation}
						We now multiply the first equivalence in \eqref{eq2:thm:associated-quasicomplex} by $S_{k-1}R_{k-1}$ from the right. By \eqref{eq1:thm:associated-quasicomplex} and the second equivalence in \eqref{eq2:thm:associated-quasicomplex}, we get the relation $R_kS_{k-1}\sim S_{k-1}R_{k-1}$. As a result, for any $k$, 
						\begin{align*}
							\hat P_k S_k+S_{k-1}\hat P_{k-1}&=R_kP_kS_k+(S_{k-1}R_{k-1})P_{k-1}\\
							&\sim R_kP_kS_k+R_kS_{k-1}P_{k-1}=R_kD_k\sim 1_{\md E_k}.
						\end{align*}
						Therefore, $\complex{S}$ has a parametrix.
					\end{proof}
				\end{lemma}
				
				A quasicomplex $\complex{P}$ that satisfies \Cref{thm:associated-quasicomplex} is called an \emph{associated quasicomplex} of $\complex{S}$, cf. \cite[p.114]{SchulzeTarkhanov1998}. If $\complex{P}$ is not a quasicomplex, then we simply say that the operators $\{P_k\}_{k=0}^N$ are \emph{associated} with $\complex{S}$. With this terminology, a finite-length quasicomplex is C*-Fredholm if and only if it has associated maps. We now show, as suggested by Segal \cite[preamble of Prop. 5.2]{Segal1970} and Putinar \cite[Eq.2.2]{Putinar1982}, that any parametrix induces a quasicomplex parametrix.
				
				\begin{theorem}\label{thm:quasicomplex-parametrix}
					Any $\alg A$-Fredholm quasicomplex has a quasicomplex parametrix.
					\begin{proof}
						Let $\complex{S}$ be an $\alg A$-Fredholm quasicomplex with parametrix $\complex{\tilde{P}}$. By the quasicomplex property, observe that $S_k\sim S_k\tilde P_kS_k$ and $S_{k-1}\sim S_{k-1}\tilde P_{k-1}S_{k-1}$. The sequence $\{P_k:=\tilde P_{k}S_{k}\tilde P_{k}\}_{k=0}^N$ is a parametrix of $\complex{S}$ with the desired property.
					\end{proof}					
				\end{theorem}
				
				\begin{corollary}\label{thm:fredholm-property-associated-quasicomplex}
					Any finite-length quasicomplex is an $\alg A$-Fredholm quasicomplex if and only if it has an associated quasicomplex.
				\end{corollary}
				
				Let us prove two corollaries of \Cref{thm:associated-quasicomplex} and \Cref{thm:quasicomplex-parametrix}.  First, in \Cref{thm:finite-rank-complex-parametrix-I} below, we generalize \Cref{thm: parametrix-finite-rank} for a certain type of quasicomplexes.
			
				\begin{corollary}\label{thm:finite-rank-complex-parametrix-I}
					If $\complex{S}$ is an $\alg A$-Fredholm quasicomplex such that $S_{k+1}S_k$ is a finite-rank operator for each $k$, then there exists a quasicomplex parametrix $\complex{P}$ such that $P_kP_{k+1}$ and the operator $C_k$ in \eqref{eq1:def:Fredholm_complexes} are finite-rank for every $k$. In particular, this holds whenever $\complex{S}$ is an adjointable, finite-length  $\alg A$-Hilbert complex with a parametrix.
					\begin{proof}
						The main idea is to maintain all the equivalences with respect to finite-rank operators, and we know from \Cref{thm: parametrix-finite-rank} that one can take the operators $\{R_k\}_{k}$ in \Cref{thm:associated-quasicomplex} such that the equivalences \eqref{eq2:thm:associated-quasicomplex} are modulo finite-rank operators. The result readily follows.
					\end{proof}
				\end{corollary}
				
				An adjointable C*-Hilbert complex can satisfy the Fredholm property in the sense of complexes or in the sense of quasicomplexes. We show that these two notions are equivalent, which is the general version of \Cref{thm:parametrix-joint-parametrix-operators}. Contrary to \Cref{thm:finite-rank-complex-parametrix-I}, the following proof does not depend on its particular case.
				
				\begin{corollary}\label{thm:fredholm-equivalence-quasicomplex-complex}
					Let $\complex{S}$  be an adjointable, finite-length $\alg A$-Hilbert complex. Then $\complex{S}$ has a joint parametrix if and only if it has a parametrix.
					\begin{proof}
						If $\complex{S}$ has a parametrix, then it clearly has a joint parametrix. Conversely, let us denote by $\complex{R_l,R_k}$ a joint parametrix of $\complex{S}$. We claim that
						\begin{equation}\label{eq1:thm:fredholm-equivalence-quasicomplex-complex}
							R_{r,k}S_k+S_{k-1}R_{r,k-1}
						\end{equation}
						is an  $\alg A$-Fredholm operator for each $k$. In this case, we can apply \Cref{thm:associated-quasicomplex} to obtain a parametrix of $\complex{S}$. Fix $k$. Since $\complex{R_l,R_r}$ is a joint parametrix, then
						\begin{align}
							R_{l,k}S_{k}+S_{k-1}R_{r,k-1}&\sim 1_{\md E_{k}},\label{eq2:thm:fredholm-equivalence-quasicomplex-complex}\\
							R_{l,k+1}S_{k+1}+S_{k}R_{r,k}&\sim 1_{\md   E_{k+1}}.\label{eq3:thm:fredholm-equivalence-quasicomplex-complex}
						\end{align}
						We want to \airquotes{replace} $R_{l,k}$  with $R_{r,k}$ in \eqref{eq2:thm:fredholm-equivalence-quasicomplex-complex}. From \eqref{eq3:thm:fredholm-equivalence-quasicomplex-complex},
						\begin{equation*}
							S_{k}R_{r,k}\sim 1_{\md E_{k+1}}-R_{l,k+1}S_{k+1}.
						\end{equation*}
						Multiplying \eqref{eq2:thm:fredholm-equivalence-quasicomplex-complex} by $R_{r,k}$ form the right, we get that
						\begin{equation}\label{eq3bis:thm:fredholm-equivalence-quasicomplex-complex}
							\begin{aligned}
								R_{r,k}&\sim R_{l,k}(S_{k}R_{r,k})+S_{k-1}R_{r,k-1}R_{r,k}\\
							&\sim R_{l,k}-R_{l,k}R_{l,k+1}S_{k+1}+S_{k-1}R_{r,k-1}R_{r,k}.
							\end{aligned}
						\end{equation}
						Thus, 
						\begin{equation}\label{eq4:thm:fredholm-equivalence-quasicomplex-complex}
							R_{r,k}S_{k}\sim R_{l,k}S_k+S_{k-1}R_{r,k-1}R_{r,k}S_k.
						\end{equation}
						Using \eqref{eq2:thm:fredholm-equivalence-quasicomplex-complex}, we derive that
						\begin{align*}
							R_{r,k}S_k+S_{k-1}R_{r,k-1}&\sim 1_{\md  E_k}+(S_{k-1}R_{r,k-1}R_{r,k}S_{k})=:\Gamma_k.
						\end{align*}
						Since $(S_{k-1}R_{r,k-1}R_{r,k}S_{k})^2=0$ (because $\complex{S}$ is a complex), then $\Gamma_k$ is invertible with inverse $1_{\md E_k}-(S_{k-1}R_{r,k-1}R_{r,k}S_{k})$. Therefore, $R_{r,k}S_k+S_{k-1}R_{r,k-1}$ is a compact perturbation of an invertible operator, so it is an $\alg A$-Fredholm operator. We can perform this process a finite number of times to obtain our claim.						
					\end{proof}
				\end{corollary}
				
				With the arguments of \Cref{thm:fredholm-equivalence-quasicomplex-complex}, one can also show that a quasicomplex $\complex{S}$ has a parametrix if and only if it has a \emph{joint parametrix}, which can be defined in the obvious way. Indeed, the map $\Gamma_k$ is $\alg A$-Fredholm with parametrix $1_{\md E_k}-(S_{k-1}R_{r,k-1}R_{r,k}S_{k})$. Hence, via associated maps, we conclude that $\complex{S}$ is $\alg A$-Fredholm. However, we do not need this notion.
					
				One might be tempted to apply again the arguments of \Cref{thm:fredholm-equivalence-quasicomplex-complex} to show that a similar result can be obtain for arbitrary C*-Fredholm complexes. However, recall that the symbol $\sim$ indicates that there is a compact operator, say $C$, coming from \eqref{eq3bis:thm:fredholm-equivalence-quasicomplex-complex} such that \eqref{eq4:thm:fredholm-equivalence-quasicomplex-complex} can be written equivalently as
				\begin{equation*}
						R_{r,k}S_{k}= R_{l,k}S_k+S_{k-1}R_{r,k-1}R_{r,k}S_k+CS_{k}.
				\end{equation*}  
				For adjointable operators, the map $CS_{k}$ is compact whereas for regular operators it might fail to be compact or even adjointable.

	\section{Fredholm equivalences}\label{sec:fredholm-equivalences}\noindent
								
				Throughout this section, we will consider a finite-length $\alg A$-Hilbert complex $\complex{t}:=\{t_k\in \reg{\md E_k,\md E_{k+1}}\}_{k=0}^N$ with bounded transform $\complex{F_t}$. 
				We recall that 
				\begin{equation*}
					F_{t_k}=t_kQ_{k}\quad \text{and}\quad F_{t_k}^*=t_k^*Q_{k^*}, \quad  0\leq k\leq N,
				\end{equation*}
				where $Q_k:=(1_{\md E_k}+t_k^*t_k)^{-1/2}$ and $Q_{k^*}:=(1_{\md E_{k+1}}+t_kt_k^*)^{-1/2}$. Both $\complex{t}$ and $\complex{F_t}$ have even Dirac operators, $\evDirac{t}$ and $\evDirac{F_t}$, and Laplace operators, $\Laplace{t}$ and $\Laplace{F_t}$, respectively. 
                The following result provides a Hilbert $\alg A$-module analogue of \cite[Theorem 2.4]{BruningLesch1992}. 
				\begin{theorem}\label{thm:fredholm-goal-I}
					Let $\complex{t}$ be a finite-length $\alg A$-Hilbert complex. All of the following objects are $\alg A$-Fredholm if and only if any of them is.
					\begin{enumerate}[(i)]
						\item The complex $\complex{t}$.
						\item The even Dirac operator $\evDirac{t}$.
						\item The Laplace operator $\Laplace{t}$.
						\item The bounded transform complex $\complex{F_t}$.
						\item The adjoint complex $\complex{t^\sharp}$.
						\item The graph-norm complex $\complex[\md E_{\graph{t}}]{t}$.
					\end{enumerate}
					\begin{proof}
						The proofs can be found in this section: the equivalence of (i) with (iv) is in \Cref{thm:fredholm-equivalence-complex-and-bounded-transform}, of (i) with (v) is in \Cref{thm:fredholm-equivalence-adjoint-complex}, of (ii) with (iii) is in \Cref{thm:fredholm-equivalence-laplace-dirac-operators}, of (iv) with (vi) is in \Cref{thm:fredholm-equivalence-graph-bounded-transform-complex}, and of (i) with (ii) is in \Cref{thm:fredholm-equivalence-complex-and-dirac-operator}.
					\end{proof}
				\end{theorem}

		\subsection{The complex and its bounded transform}\label{subsec:equivalence-bounded-transform-complex}
				The next statement proves the equivalence of (i) and (iv) in \Cref{thm:fredholm-goal-I}, and it will allow us to reduce several arguments to the case of adjointable complexes.

				\begin{proposition}\label{thm:fredholm-equivalence-complex-and-bounded-transform}
					The complex $\complex{t}$ is an $\alg A$-Fredholm complex if and only if its bounded transform $\complex{F_t}$ is an $\alg A$-Fredholm complex.
					\begin{proof}
						Assume the bounded transform $\complex{F_t}$ is an $\alg A$-Fredholm complex. By \Cref{thm:fredholm-equivalence-quasicomplex-complex}, there exists a parametrix $\complex{P}$ of $\complex{F_t}$, that is, there exists a collection of adjointable operators $\complex{P}=\{P_k\}_{k=0}^N$ such that
						\begin{equation}\label{eq1:thm:fredholm-equivalence-complex-and-bounded-transform}
							P_{k}F_{t_k}+F_{t_{k-1}}P_{k-1}=1_{\md E_k}-C_k,\quad 0\leq k\leq N+1,
						\end{equation}
						for some collection of compact operators $\{C_k\}_{k=0}^N$. For each $k$, let us define $P_{l,k}:=P_kQ_{k^*}$ and $P_{r,k}:=Q_kP_k$. We claim that $\complex{P_l,P_r}$ induces a joint parametrix for $\complex{t}$. Indeed, recall that $F_{t_k}=t_kQ_k$, and $Q_{k^*}t_k\subset F_{t_k}$; moreover, $\ran{Q_{k}}=\dom{t_k}$. Thus, $\ran{P_{r,k}}\subset \dom{t_k}$ and, from \eqref{eq1:thm:fredholm-equivalence-complex-and-bounded-transform}, we get that
						\begin{equation*}
							P_kQ_{t^*}t_k+t_{k-1}Q_{k-1}P_{k-1}=1_{\md E_k}-C_k\text{ in }\dom{t_k}, \quad  0\leq k\leq N+1.
						\end{equation*}
						Hence, $\complex{P_l,P_r}$ is a joint parametrix of $\complex{t}$. (Observe that \Cref{thm:fredholm-equivalence-quasicomplex-complex} is not necessary, since a joint parametrix of $\complex{F_t}$ does the work with the same arguments. However, it will be useful in \Cref{thm:finite-rank-complex-parametrix-II,thm:fredholm-equivalence-adjoint-complex} below.)
						
						Conversely, assume that $\complex{P_l,P_r}$ is a joint parametrix of $\complex{t}$. We want to induce a joint parametrix for $\complex{F_t}$. Fix $k$. By \Cref{thm:properties-bounded-transform}, we get the relations $F_{t_k}Q_k=Q_{k^*}F_{t_k}$, $F_{t_k}Q_{k}^{-1}\subset Q_{k^*}^{-1}F_{t_k}$, and $F_{t_k}=t_kQ_k$. Thus,
						\begin{align}
							P_{l,k}t_k+t_{k-1}P_{r,k}&=P_{l,k}F_{t_k}Q_k^{-1}+F_{t_{k-1}}Q_{k-1}^{-1}P_{r,k}\\ 	&=P_{l,k}F_{t_k}Q_k^{-1}+Q_{(k-1)^*}^{-1}F_{t_{k-1}}P_{r,k}=1_{\dom{t_k}}-C_k|_{\dom{t_k}}.\label{eq3:thm:fredholm-equivalence-complex-and-bounded-transform}
						\end{align}
						By \Cref{thm:properties-resolvent-sequence}, $Q_{(k-1)^*}$ and $Q_k$ commute, and $Q_k^2Q_{(k-1)^*}^2=1_{\md E_k}-F_{t_k}^*F_{t_k}-F_{t_{k-1}}F_{t_{k-1}}^*$. Multiplying \eqref{eq3:thm:fredholm-equivalence-complex-and-bounded-transform} by $Q_{(k-1)^*}^2$ from the left, and, by $Q_k^2$ from the right, we obtain that
						\begin{align*}
							Q_{(k-1)^*}^2P_{l,k}F_{t_k}Q_k+Q_{(k-1)^*}F_{t_{k-1}}P_{r,k}Q_k^2
							&=Q_{(k-1)^*}^2P_{l,k}Q_{k^*}F_{t_k}+F_{t_{k-1}}Q_{(k-1)}P_{r,k}Q_k^2\\
							&\sim 1_{\md E_k}-F_{t_k}^*F_{t_k}-F_{t_{k-1}}F_{t_{k-1}}^*
						\end{align*}
						Thus $(Q_{(k-1)^*}^2P_{l,k}Q_{k^*}+F_{t_k}^*)F_{t_k}+F_{t_{k-1}}(Q_{(k-1)}P_{r,k}Q_k^2+F_{t_{k-1}}^*)\sim 1_{\md E_k}$ for every $k$, so that $\complex{F_t}$ has a joint parametrix, as desired.
					\end{proof}
				\end{proposition}
				
				The next corollary is the unbounded version of \Cref{thm:finite-rank-complex-parametrix-I}.
				
				\begin{corollary}\label{thm:finite-rank-complex-parametrix-II}
					If $\complex{t}$ is an $\alg A$-Fredholm complex, then there exists a joint parame\-trix of $\complex{t}$ such that the operators $\{C_k\}_{k=0}^N$ in \eqref{eq2:def:joint-parametrix-complex} are finite-rank operators.
					\begin{proof}
						By the proof of \Cref{thm:fredholm-equivalence-complex-and-bounded-transform}, a parametrix of $\complex{F_t}$ induces a joint parame\-trix for $\complex{t}$ with the same collection of compact operators. Thus, it suffices to prove the statement for adjointable complexes, that is, \Cref{thm:finite-rank-complex-parametrix-I}.
					\end{proof}
				\end{corollary}	

				We now show that our definition of the Fredholm property for regular operators is equivalent to \cite[Definition 2.1]{Joachim2003}.																
				\begin{corollary}\label{thm:consistency-joachim}
					Consider a regular operator $t\in \reg{\md E,\md F}$. Then $t$ is $\alg A$-Fredholm (as a finite-length complex) if and only if it has a pseudo-left and a pseudo-right inverse.
					\begin{proof}
						By \Cref{thm:fredholm-equivalence-complex-and-bounded-transform}, the operator $t$ is $\alg A$-Fredholm if and only if its bounded transform $F_t$ is $\alg A$-Fredholm. With the same arguments from \cite[Theorem 2.2]{Joachim2003}, one can show that $t$ has a pseudo-left and a pseudo-right inverse if and only if $F_t$ is $\alg A$-Fredholm.
					\end{proof}
				\end{corollary}
				
        \subsection{The complex and its adjoint}\label{subsec:equivalence-adjoint-complex}
				\begin{corollary}\label{thm:fredholm-equivalence-adjoint-complex}
					The complex $\complex{t}$ is an $\alg A$-Fredholm complex if and only if its adjoint complex $\complex{t^\sharp}$ is an $\alg A$-Fredholm complex. Moreover, there exists a joint parametrix $\complex{P_l,P_r}$ of $\complex{t}$ such that $\complex{P_r^\sharp,P_l^\sharp}:=\{(P_{r,k}^*,P_{l,k}^*)\}_{k=0}^N$ is a joint parametrix of $\complex{t^\sharp}$.
					\begin{proof}
						By \Cref{thm:fredholm-equivalence-quasicomplex-complex} and \Cref{thm:fredholm-equivalence-complex-and-bounded-transform}, any joint parametrix of $\complex{t}$ induces a parametrix for its bounded transform. In turn, this parametrix of the bounded transform induces the desired joint parametrix for $\complex{t}$ (as constructed in \Cref{thm:fredholm-equivalence-complex-and-bounded-transform}) and for $\complex{t^\sharp}$; see \Cref{thm:basic-properties-fredholm-quasicomplexes}\ref{item1:thm:basic-properties-fredholm-quasicomplexes}.
					\end{proof}
				\end{corollary}

        \subsection{The Dirac operator and the Laplace operator}\label{equivalence:dirac-laplace-operators}
				For adjointable operators, one can show that $T$ is $\alg A$-Fredholm if and only if both $T^*T$ and $TT^*$ are $\alg A$-Fredholm. We want to obtain a similar statement for the even Dirac operator. The following statements breakdown the proof.

				\begin{lemma}\label{thm:fredholm-equivalence-dirac-operators}
					The operator $\evDirac{t}$ is $\alg A$-Fredholm if and only if $\evDirac{F_t}$ is $\alg A$-Fredholm.
					\begin{proof}
						By \Cref{thm:bounded-even-dirac-operator}, $\evDirac{F_t}$ is the bounded transform of $\evDirac{t}$. The result follows by \Cref{thm:fredholm-equivalence-complex-and-bounded-transform}.
					\end{proof}
				\end{lemma}			
				
				As we mentioned in the preamble of \Cref{thm:laplace-operator-bounded-complex}, the Laplace operator of $\complex{F_t}$ is not the bounded transform of the Laplace operator of $\complex{t}$, so their Fredholm equivalence is not an immediate consequence of \Cref{thm:fredholm-equivalence-complex-and-bounded-transform}. Nevertheless, we can prove the following result. 
				\begin{lemma}\label{thm:fredholm-equivalence-laplace-operators}
					The operator $\Laplace{F_t}$ is $\alg A$-Fredholm if and only if $\Laplace{t}$ is $\alg A$-Fredholm.
					\begin{proof}
						It suffices to prove, for each $k\in\{0,1,\ldots,N\}$, that $\Laplace{t}_k$ is $\alg A$-Fredholm if and only if $\Laplace{F_t}_k$ is $\alg A$-Fredholm. By \Cref{thm:laplace-operator-bounded-complex}, we know that $\Laplace{F_t}_k=\Laplace{t}_kQ_k^2Q_{(k-1)^*}^2$. Moreover,
						\begin{align*}
							Q_k^2Q_{(k-1)^*}^2\Laplace{t}_k&\subset \Laplace{F_t}_k
						\end{align*}
						by taking adjoints. Hence, if $P$ is a parametrix of $\Laplace{F_t}_k$, then $P_l:=PQ_k^2Q_{(k-1)^*}^2$ and $P_r:=Q_k^2Q_{(k-1)^*}^2P$ define a joint parametrix for $\Laplace{t}_k$. Conversely, if $\Laplace{t}_k$ has a joint parametrix $(P_l,P_r)$, then 
						\begin{align}
							P_l\Laplace{F_t}_kQ_{(k-1)^*}^{-2}Q_{k}^{-2}&=P_l\Laplace{t}_k\subset 1_{\md E_k}+C_l,\label{eq1:thm:fredholm-equivalence-laplace-operators}\\
							Q_{(k-1)^*}^{-2}Q_{k}^{-2}\Laplace{F_t}_kP_r&=\Laplace{t}_kP_r=1_{\md E_k}+C_r.\label{eq2:thm:fredholm-equivalence-laplace-operators}
						\end{align}
						Since $Q_k^2Q_{(k-1)^*}^2=1_{\md E_k}-\Laplace{F_t}_k$ (\Cref{thm:properties-resolvent-sequence}\ref{item1:thm:properties-resolvent-sequence}), then $(P_l+1_{\md E_k})\Laplace{F_t}_k\sim 1_{\md E_k}$ by \eqref{eq1:thm:fredholm-equivalence-laplace-operators}, and $\Laplace{F_t}_k(1_{\md E_k}+P_r)\sim 1_{\md E_k}$ by \eqref{eq2:thm:fredholm-equivalence-laplace-operators}. Therefore, $\Laplace{F_t}_k$ is $\alg A$-Fredholm.
					\end{proof}
				\end{lemma}

				\begin{proposition}\label{thm:fredholm-equivalence-laplace-dirac-operators}
					The Laplace operator $\Laplace{t}$ is $\alg A$-Fredholm if and only if the even Dirac operator $\evDirac{t}$ is $\alg A$-Fredholm.
					\begin{proof}
						Using \Cref{thm:parmetrices-fredholm-operators}, one can show that, for any $T$ adjointable, $T$ is $\alg A$-Fredholm if and only if $TT^*$ and $T^*T$ are $\alg A$-Fredholm. In particular, $\evDirac{F_t}$ is $\alg A$-Fredholm if and only if both $\Laplace{F_t}_{\mathrm{ev}}$ and $\Laplace{F_t}_{\mathrm{odd}}$ are $\alg A$-Fredholm (\Cref{thm:dirac-operator-and-laplace-operator}). 
                        The statement then follows from \Cref{thm:fredholm-equivalence-dirac-operators,thm:fredholm-equivalence-laplace-operators}.
					\end{proof}
				\end{proposition}			
				
				\begin{remark}\label{rmk:fredholm-equivalence-laplace-dirac-operators-quasicomplexes}
					Let $\complex{S}$ be a finite-length $\alg A$-Hilbert quasicomplex. By \Cref{thm:laplace-dirac-operators-for-quasicomplex}, the Laplace operator is a compact perturbation of $\oddDirac{S}\evDirac{S}\oplus \evDirac{S}\oddDirac{S}$. Applying the first part of the proof in \Cref{thm:fredholm-equivalence-laplace-dirac-operators} and using stability under compact perturbations (\Cref{thm:basic-properties-fredholm-quasicomplexes}\ref{item3:thm:basic-properties-fredholm-quasicomplexes}), we see that the statement of \Cref{thm:fredholm-equivalence-laplace-dirac-operators} also holds for quasicomplexes. 
				\end{remark}

        \subsection{The bounded transform complex and the graph-norm complex}\label{eqivalence:the-graph-complex}
				The graph-norm complex was advertised in \Cref{subsec:graph-norm-complex} as an alternative to the bounded transform complex. However, observe that its Fredholm property is proved below with the aid of the bounded transform, and, as far as our discussion shows, there is no direct connection with the Fredholm property of the complex.												
				
				\begin{lemma}\label{thm:isomorphic-fredholm-complexes}
					Let $\complex{r}$ and $\complex[\md F]{s}$ be  isomorphic finite-length $\alg A$-Hilbert complexes. Then $\complex{r}$ is an $\alg A$-Fredholm complex if and only if $\complex[\md F]{s}$ is an $\alg A$-Fredholm complex.
					\begin{proof}
						Assume that $G$ and $H$ are adjointable maps of complexes such that $HG=1_{\complex{r}}$ and $GH=1_{\complex[\md F]{s}}$, and that $\complex{r}$ is an $\alg A$-Fredholm complex with joint parametrix $\complex{P_l,P_r}$. Put $\tilde P_{l,k}:=G_kP_{l,k}H_{k+1}$ and $\tilde P_{r,k}:=G_kP_{r,k}H_{k+1}$ for every $0\leq k\leq N$. Then $\complex[\md F]{\tilde P_l,\tilde P_r}$ is a joint parametrix of $\complex[\md F]{s}$. By symmetry, the result readily follows.
					\end{proof}
				\end{lemma}
															
				\begin{corollary}\label{thm:fredholm-equivalence-graph-bounded-transform-complex}
					The graph-norm complex $\complex[\md E_{\graph{t}}]{t}$ is an $\alg A$-Fredholm complex if and only if the bounded transform complex $\complex{F_t}$ is an $\alg A$-Fredholm complex.
					\begin{proof}
						By \Cref{thm:isomorphism_graph_and_boundedtransform} and \Cref{example:isomorphism_graph_and_boundedtransform}, $\complex{F_t}$ and $\complex[\md E_{\graph{t}}]{t}$ are isomorphic.
					\end{proof}
				\end{corollary}

        \subsection{The complex and its Dirac operator} \label{subsec:fredholm-equivalence-complex-dirac-operator}
				We now want to investigate the Fredholm property of a complex by means of its even Dirac operator. Every known proof of this equivalence for Hilbert complexes and quasicomplexes
				relies on the strong Hodge decomposition and Atkinson’s theorem (or Riesz’s lemma), see, e.g., \cite[Theorem 3, Sec.3.2.3.1]{RempelSchulze1982} and \cite[Proposition 2.6]{Putinar1982}. In this subsection, we provide a non-standard proof that applies to the general case of C*-Hilbert quasicomplexes, and we obtain, as a corollary, the result for arbitrary C*-Hilbert complexes.

				\subsubsection{Putinar's functor}
				We must take a detour to define a map $\phi_\Sigma$ acting as a \emph{covariant functor}. This map has been studied in \cite{Putinar1982} and \cite{SchulzeTarkhanov1998} for Banach and Fréchet spaces, respectively. Our goal is to extend it to the case of Hilbert C*-modules. By \Cref{thm:fredholm-equivalence-complex-and-bounded-transform} and \Cref{thm:fredholm-equivalence-dirac-operators}, it suffices to consider adjointable complexes. Let $\md E$, $\md F$ and $\Sigma$ be Hilbert $\alg A$-modules. Put
				\begin{equation*}
					\phi_{\Sigma }(\md E):=\adj{\Sigma,\md E}/\comp{\Sigma,\md E},
				\end{equation*}
				and, for $T\in\adj{\md E,\md F}$, we define a map $\phi_{\Sigma}(T)\colon  \phi_{\Sigma}(\md E)\to \phi_{\Sigma}(\md F)$ given by
				\begin{equation*}
					\phi_{\Sigma}(T)(e+\comp{\Sigma, \md E}):=Te+\comp{\Sigma,\md F}, \quad\forall  e\in\adj{\Sigma,\md E}.
				\end{equation*}
				
				Observe that the Calkin algebra $\mathcal C(\Sigma):=\phi_\Sigma(\Sigma)$ is a C*-algebra with the quotient structure. Thus, we would like to view $\phi_\Sigma(\md E)$ as a Hilbert $\mathcal C(\Sigma)$-module. Let us first consider the following lemma.

				\begin{lemma}\label{thm:C*-like-property-quotient}
					Let $\md E$ and $\Sigma$ be Hilbert $\alg A$-modules. If $e\in \adj{\Sigma,\md E}$, then
					\begin{equation*}
						\norm[\mathcal C(\Sigma)]{e^*e+\comp{\Sigma}}\leq \norm[\adj{\Sigma,\md E}/\comp{\Sigma,\md E}]{e+\comp{\Sigma,\md E}}^2\leq 2\norm[\mathcal C(\Sigma)]{e^*e+\comp{\Sigma}}. 
					\end{equation*}
					\begin{proof}
						Let $e\in \adj{\Sigma,\md E}$. By the identity $\norm{e}^2=\norm{e^*e}$ and definition of the quotient norm, one can show that $\norm{e^*e+\comp{\Sigma}} \leq \norm{e+\comp{\Sigma,\md E}}^2$. For the second inequality, let $\{u_\lambda\}_{\lambda\in \Lambda}$ be a (positive, contractive) approximate unit of $\comp{\Sigma}$. By \cite[Theorem 3.1.3]{Murphy1990},
						\begin{equation*}
							\norm{e^*e+\comp{\Sigma}}=\lim_{\lambda}\norm{e^*e-e^*eu_{\lambda}}.
						\end{equation*}
						On the other hand, for every $\lambda\in \Lambda$,
						\begin{align*}
							\norm{e+\comp{\Sigma,\md E}}^2\leq  \norm{e-eu_\lambda}^2&=\norm{e^*e-e^*eu_\lambda-u_\lambda e^*e+u_\lambda e^*eu_\lambda}\\
							&\leq 2\norm{e^*e-e^*eu_\lambda}.
						\end{align*}
						The result follows after taking the limit.									
					\end{proof}
				\end{lemma}
				
				We now prove the desired properties of \emph{Putinar's functor}; for this, we first need to endow $\phi_{\Sigma}(\md E)$ with a Hilbert $\mathcal C(\Sigma)$-module structure. This is the content of the following theorem.

				\begin{theorem}\label{thm:properties-putinar's-functor}
					Let $\md E, \md F$, and $\Sigma$ be arbitrary Hilbert $\alg A$-modules, and let us define on $\phi_{\Sigma}(\md E)$ the following structure.
					\begin{itemize}
						\item A right $\mathcal C(\Sigma)$-action: For $e\in \adj{\Sigma,\md E}$ and $\alpha\in \adj{\Sigma}$,
						\begin{equation*}
							(e+\comp{\Sigma,\md E}) \cdot (\alpha+\comp{\Sigma}):= e\alpha+\comp{\Sigma,\md E}.
						\end{equation*}
						\item A $\mathcal C(\Sigma)$-valued map: For $e_1,e_2\in \adj{\Sigma,\md E}$,
						\begin{equation*}
							\ev{e_1+\comp{\Sigma,\md E}}{e_2+\comp{\Sigma,\md E}}:=e_1^*e_2+\comp{\Sigma}.
						\end{equation*}
					\end{itemize}
					Then $\phi_{\Sigma}(\md E)$ is a Hilbert $\mathcal{C}(\Sigma)$-module. Moreover, for every $T\in \adj{\md E,\md F}$, the map $\phi_{\Sigma}(T)$ is adjointable in $\adj[\mathcal{C}(\Sigma)]{\phi_{\Sigma}(\md E),\phi_{\Sigma}(\md F)}$ with adjoint $\phi_{\Sigma}(T^*)$. In addition, the map $\phi_{\Sigma}$ is a Banach functor from the category of Hilbert $\alg A$-modules into the category of Hilbert $\mathcal C(\Sigma)$-modules:
					\begin{enumerate}[(i)]
						\item If $\md E$ is a Hilbert $\alg A$-module, then $\phi_{\Sigma}(1_{\md E})=1_{\phi_{\Sigma}(\md E)}$.
						\item For $T\in\adj{\md E,\md F}$ and $S\in \adj{\md F,\md G}$, we have $\phi_{\Sigma}(ST)=\phi_{\Sigma}(S)\phi_{\Sigma}(T)$.														
						\item The map 
						\[\phi_\Sigma(\cdot)\colon \adj{\md E,\md F}\to \adj[\mathcal{C}(\Sigma)]{\phi_\Sigma(\md E),\phi_\Sigma(\md F)}\]
						is linear and continuous.
					\end{enumerate}
					Furthermore,
					\begin{enumerate}
						\item[(iv)] An operator $T\in \adj{\md E,\md F}$ is compact if and only if $\phi_{\Sigma}(T)=0$ for every $\Sigma$.
						\item[(v)] If  $\complex{T}$ is an $\alg A$-Hilbert quasicomplex, then $\phi_{\Sigma}$ induces an adjointable $\mathcal C(\Sigma)$-Hilbert complex $\complex[\phi_{\Sigma}(\md E)]{\phi_{\Sigma}(T)}:=\{\phi_{\Sigma}(T_k)\}_{k\in \mathbb Z}$.
					\end{enumerate}
					\begin{proof}
						Let $\Sigma$ and $\md E$ be Hilbert $\alg A$-modules. Recall that $\comp{\Sigma,\md E}$ is a subspace of $\adj{\Sigma,\md E}$, so $\phi_{\Sigma}(\md E)$ is a vector space. Moreover, the right $\mathcal C(\Sigma)$-action is well defined and makes $\phi_{\Sigma}(\md E)$ a $\mathcal{C}(\Sigma)$-module. We now show that the $\mathcal{C}(\Sigma)$-valued map is a $\mathcal{C}(\Sigma)$-valued inner product. It is straightforward to check that it is a well-defined sesquilinear map that respects the module action. Take $e\in \adj{\Sigma,\md E}$. Since $e^*e$ is positive in $\adj{\Sigma}$, see, e.g., \cite[Proposition 15.2.5]{WeggeOlsen1993}, we have that
						\begin{equation*}
							\ev{e+\comp{\Sigma,\md E}}{e+\comp{\Sigma,\md E}}=e^*e+\comp{\Sigma}
						\end{equation*}
						is positive in $\mathcal C(\Sigma)$. Moreover, if $e^*e$ is compact, then, by \Cref{thm:C*-like-property-quotient}, $e$ is compact. Therefore, the $\mathcal{C}(\Sigma)$-valued map is positive definite. In conclusion, $\phi_{\Sigma}(\md E)$ is a pre-Hilbert $\mathcal C(\Sigma)$-module.  Finally, observe that the norm induced by the inner product is
						\begin{equation*}
							\norm{e+\comp{\Sigma,\md E}}:=\norm[\mathcal C(\Sigma)]{e^*e+\comp{\Sigma}}^{1/2},\quad \forall e\in \adj{\Sigma,\md E},
						\end{equation*}
						which by \Cref{thm:C*-like-property-quotient} is equivalent to the norm of a Banach space. Therefore, $\phi_{\Sigma}(\md E)$ is complete, and we conclude that it is a Hilbert $\mathcal C(\Sigma)$-module. One can easily check that $\phi_\Sigma(T)$ is adjointable for any $T\in \adj{\md E,\md F}$, with adjoint $\phi_\Sigma(T^*)$.  								
						
						Properties (i) and (ii) follow immediately from the definition. To show (iii), it suffices to prove that $\phi_\Sigma(\cdot)\colon \adj{\md E}\to \adj[\mathcal{C}(\Sigma)]{\phi_\Sigma(\md E)}$ is linear continuous because $\norm{T}^2=\norm{T^*T}$ for every adjointable operator. By \cite[Theorem 2.1.7]{Murphy1990}, we know that $\norm{\phi_\Sigma(\cdot)}\leq 1$ because $\phi_\Sigma(\cdot)$ is a *-homomorphism between C*-algebras.
						
						(iv) Take $T\in\adj{\md E,\md F}$ such that $\phi_{\Sigma}(T)=0$ for every Hilbert $\alg A$-module $\Sigma$. In particular, for $\Sigma=\md E$, we have that $T+\comp{\md E,\md F}=\phi_{\Sigma}(T)(1_{\md E})=\comp{\md E,\md F}$ whence $T$ is compact. The converse implication is trivial.
						
						(v) Let $\complex{T}=\{T_k\}_{k\in\mathbb Z}$ be an $\alg A$-Hilbert quasicomplex, and let $\Sigma$ be a Hilbert $\alg A$-module. By properties (ii) and (iv) of $\phi_{\Sigma}$, we have that
						\begin{equation*}
							\phi_{\Sigma}(T_{k+1})\phi_{\Sigma}(T_k)=\phi_{\Sigma}(T_{k+1}T_k)=0,\quad  \forall k\in \mathbb Z.
						\end{equation*}
						Therefore, $\complex[\phi_{\Sigma}(\md E)]{\phi_{\Sigma}(T)}:=\{\phi_{\Sigma}(T_k)\}_{k\in \mathbb Z}$ is a $\mathcal C(\Sigma)$-Hilbert complex.
					\end{proof}
				\end{theorem}

				\subsubsection{The Fredholm property via Putinar's functor}
				Putinar provides in \cite[Proposition 2.3]{Putinar1982} an equivalent definition of the Fredholm property for complexes of Banach spaces using the functor $\phi_{\Sigma}$. The same description can be applied to our case without modifying the arguments. This is the content of the following theorem.

				\begin{theorem}[cf. {\cite[Proposition 2.3]{Putinar1982}}]\label{thm:fredholm-equivalence-putinar-criterion}
					Let $\complex{T}$ be a finite-length $\alg A$-Hilbert quasicomplex. Then $\complex{T}$ is an $\alg A$-Fredholm quasicomplex if and only if, for each Hilbert $\alg A$-module $\Sigma$, the induced $\mathcal{C}(\Sigma)$-Hilbert complex
					\begin{equation*}
						\complex[\phi_{\Sigma}(\md E)]{\phi_{\Sigma}(T)}
					\end{equation*}
                    is an exact sequence of Hilbert $\mathcal C(\Sigma)$-modules.
					\begin{proof}
						The proof is similar to \cite[Theorem 5.1.3]{SchulzeTarkhanov1998}, using the results for chain homotopies in \Cref{sec:maps_complexes}.
					\end{proof}
				\end{theorem}

				The reader can verify that the proof of \Cref{thm:fredholm-equivalence-putinar-criterion} does not use the Hilbert C*-module structure of the complex $\complex[\phi_\Sigma(\md E)]{\phi_\Sigma(T)}$ but only its vector space structure. However, due to the following statement, cf. \cite[p.396]{AtiyahBott1967} and \cite[Lemma 9]{Krysl2015}, the Hilbert C*-module structure will be crucial in the proof of \Cref{thm:fredholm-equivalence-quasicomplex-and-dirac-operator} below.

				\begin{lemma}\label{thm:exact-automorphism}
					Let $\md E_0, \md E_1$, and $\md E_2$ be Hilbert $\alg A$-modules, and let $T_0$ and $T_1$ be adjointable operators such that the sequence
					\[\begin{tikzcd}
						{\md E_0} & {\md E_1} & {\md E_2}
						\arrow["T_0", from=1-1, to=1-2]
						\arrow["T_1", from=1-2, to=1-3]
					\end{tikzcd}\]
					is exact at $\md E_1$. If $\ran{T_1}$ is closed, then the operator $\Delta:=T_1^*T_1+T_0T_0^*\in \adj{\md E_1}$ is self-adjoint and invertible.
					\begin{proof}
						Observe that $\ran{T_0}=\ker{T_1}$, so that $T_0$ and $T_1$ both have closed range. By \Cref{thm:kernel-kth-laplace-operator} and \Cref{thm:hodge-decomposition-equivalences}, we deduce that $\Delta$ has closed range and $\ker \Delta=\ker{T_0^*}\cap \ker{T_1}$. Thus, $\ker \Delta=0$, by exactness at $\md E_1$ and the equation $\ker{T_0^*}=(\ran{T_0})^\perp$. From \eqref{eq3:subsec:weak-hodge-decomposition}, we know that $\md E_1=\ker \Delta\oplus \ran{\Delta}$. Therefore, $\Delta$ is a bijective self-adjoint operator.
					\end{proof}
				\end{lemma}

				\begin{corollary}\label{thm:exact-complex-isomorphism-laplacians}
					Let $\complex{S}$ be an adjointable $\alg A$-Hilbert complex. The complex $\complex{S}$ is exact if and only if $\Laplace{S}_k$ is an isomorphism for every $k\in \mathbb Z$.
					\begin{proof}
                        Fix $k\in \mathbb Z$. Assume that $\Laplace{S}_k$ is an isomorphism. Then $S_k$ and $S_{k-1}$ have closed range, and $\ker{S_k}=\ker{\Laplace{S}_k}\oplus \ran{S_{k-1}}=\ran{S_{k-1}}$ (\Cref{thm:hodge-decomposition-equivalences}). The converse implication is a consequence of \Cref{thm:exact-automorphism}.
					\end{proof}
				\end{corollary}

				We are now in position to prove the main result of this subsection. The reader might want to compare our use of the functor $\phi_\Sigma$  in \eqref{eq1:thm:fredholm-equivalence-quasicomplex-and-dirac-operator} with the use of  the \emph{symbol of elliptic operators} in \cite[p.396]{AtiyahBott1967} or with \eqref{eq:symbol_elliptic_operators} below.

				\begin{proposition}\label{thm:fredholm-equivalence-quasicomplex-and-dirac-operator}
					Let $\complex{S}$ be a finite-length $\alg A$-Hilbert quasicomplex. Then $\complex{S}$ is an $\alg A$-Fredholm quasicomplex if and only if its even Dirac operator $\evDirac{S}$ is $\alg A$-Fredholm.
					\begin{proof}
						Assume that $\evDirac{S}$ is $\alg A$-Fredholm. By \Cref{thm:laplace-dirac-operators-for-quasicomplex}, $\oddDirac{S}\evDirac{S}\oplus \evDirac{S}\oddDirac{S}\sim \Laplace{S}_{\mathrm{ev}}\oplus\Laplace{S}_{\mathrm{odd}}$. From \Cref{rmk:fredholm-equivalence-laplace-dirac-operators-quasicomplexes}, we know that each $k$-th Laplace operator $\Laplace{S}_k$ is $\alg A$-Fredholm. Thus, $\complex S$ is $\alg A$-Fredholm by \Cref{thm:associated-quasicomplex}. 
                        
                        Conversely, suppose that $\complex{S}:=\{S_k\}_{k=0}^N$ is $\alg A$-Fredholm. Let $\Sigma$ be a Hilbert $\alg A$-module. By \Cref{thm:fredholm-equivalence-putinar-criterion}, the complex $\complex[\phi_{\Sigma}(\md E)]{\phi_{\Sigma}(S)}$ is exact. Hence,
						\begin{equation}\label{eq1:thm:fredholm-equivalence-quasicomplex-and-dirac-operator}
							\phi_{\Sigma}(S_k)^*\phi_{\Sigma}(S_k)+\phi_{\Sigma}(S_{k-1})\phi_{\Sigma}(S_{k-1})^*=\phi_{\Sigma}(\Laplace{S}_k), \quad 0\leq k\leq N+1,
						\end{equation}
						is an isomorphism on $\phi_{\Sigma}(\md E_k)$ by \Cref{thm:exact-automorphism}. 
                        Since $\Sigma$ was arbitrary, it follows from \Cref{thm:fredholm-equivalence-putinar-criterion} that the operator $\Laplace{S}_k$ (viewed as a complex of length one) is $\alg A$-Fredholm for every $k$, so $\evDirac{S}$ is $\alg A$-Fredholm by \Cref{rmk:fredholm-equivalence-laplace-dirac-operators-quasicomplexes}.
					\end{proof}
				\end{proposition}

				\begin{corollary}\label{thm:fredholm-equivalence-complex-and-dirac-operator}
					The complex $\complex{t}$ is an $\alg A$-Fredholm complex if and only if its even Dirac operator $\evDirac{t}$ is an $\alg A$-Fredholm operator.
					\begin{proof}
						By \Cref{thm:fredholm-equivalence-complex-and-bounded-transform} and \Cref{thm:fredholm-equivalence-dirac-operators}, it suffices to prove the statement for adjointable complexes, that is, \Cref{thm:fredholm-equivalence-quasicomplex-and-dirac-operator}.
					\end{proof}
				\end{corollary}

				Observe that  we have also proven the following result.
				
				\begin{theorem}\label{thm:fredholm-goal-I_quasicomplexes}
					Let $\complex{S}$ be a finite-length $\alg A$-Hilbert quasicomplex. All of the following objects are $\alg A$-Fredholm if and only if any of them is.
					\begin{enumerate}[(i)]
						\item The quasicomplex $\complex{S}$.
						\item The even Dirac operator $\evDirac{S}$.
						\item The Laplace operator $\Laplace{S}$.
						\item The adjoint quasicomplex $\complex{S^\sharp}$.
					\end{enumerate}
					\begin{proof}
						The equivalence of (i) with (ii) is in \Cref{thm:fredholm-equivalence-quasicomplex-and-dirac-operator}, of (ii) with (iii) is in \Cref{rmk:fredholm-equivalence-laplace-dirac-operators-quasicomplexes}, and of (i) with (iv) is in \Cref{thm:basic-properties-fredholm-quasicomplexes}\ref{item1:thm:basic-properties-fredholm-quasicomplexes}.
					\end{proof}
				\end{theorem}
				
	\section{The Fredholm index}\label{sec:fredholm_index_map}\noindent
		\subsection{Properties of the Fredholm index}
				Recall the Fredholm index map, denoted by $\ind{\cdot}$, such that $\ind{T}\in \K$ for any $\alg A$-Fredholm operator $T\in \adj{\md E,\md F}$ \cite[Section 3]{Exel1993}. In the next theorem, we gather some properties about the index for Fredholm operators that we shall need in the sequel.
																
				\begin{theorem}[{\cite[Section 3]{Exel1993}}] \label{thm:index-map-adjointable-operators}
					Let $\md X$ and $\md Y$ be Hilbert $\alg A$-modules. For any $\alg A$-Fredholm operator $T\in \adj{\md E,\md F}$, the following statements hold.
					\begin{enumerate}[(i)]
						\item $\ind{T^*}=-\ind{T}$. \label{item1:thm:index-map-adjointable-operators}
						\item If $T$ has closed range, then $\ind{T}=\Kclass{\ker T}-\Kclass{\ker {T^*}}$. \label{item2:thm:index-map-adjointable-operators}
						\item If $C\in \comp{\md E,\md F}$, then $\ind{T}=\ind{T+C}$. \label{item3:thm:index-map-adjointable-operators}
						\item If $P\in \adj{\md F,\md E}$ is a parametrix of $T$, then $\ind{P}=-\ind{T}$.\label{item4:thm:index-map-adjointable-operators}
						\item If $U\in \adj{\md X,\md E}$ and $V\in \adj{\md F,\md Y}$ are invertible, then $\ind{VTU}=\ind{T}$.\label{item5:thm:index-map-adjointable-operators}
						\item If $S\in\adj{\md X,\md Y}$ is $\alg A$-Fredholm, then $\ind{T\oplus S}=\ind{T}+\ind{S}$.\label{item6:thm:index-map-adjointable-operators}
						\item There exists an $\epsilon>0$ such that $S\in \adj{\md E,\md F}$ is $\alg A$-Fredholm and $\ind{S}=\ind{T}$ whenever $\norm{S-T}<\epsilon$. \label{item7:thm:index-map-adjointable-operators}
						\item If $S\in \adj{\md F,\md X}$ is $\alg A$-Fredholm, then $\ind{ST}=\ind{S}+\ind{T}$.\label{item8:thm:index-map-adjointable-operators}
					\end{enumerate}
				\end{theorem}
				
				\begin{remark}\label{rmk:regular_operators_exel}
					Exel uses in \cite{Exel1993} the term \airquotes{regular} for an adjointable operator which has a generalized inverse (we refrain from using this terminology, since we reserve the word \airquotes{regular} for densely defined operators as described in \Cref{subsec:operators-hilbert-C*-modules}). 
					It can be shown that an operator has closed range if and only if it has a generalized inverse. In particular, if $T$ is $\alg A$-Fredholm with closed range, then $1_{\ker T}$ and $1_{\ker {T^*}}$ are compact \cite[Proposition 3.3]{Exel1993}, and the $K$-classes of $\ker T$ and $\ker {T^*}$ are well defined in $\K$. Moreover, one can prove that if $T$ has closed range, and $1_{\ker T}$ and $1_{\ker {T^*}}$ are compact, then $T$ is $\alg A$-Fredholm. The converse implication (giving the equivalence) is known in the Hilbert-space case as Atkinson's theorem, but it fails to hold in general on Hilbert C*-modules (not every $\alg A$-Fredholm operator has closed range). 
				\end{remark}
				
				Recall that for a regular $\alg A$-Fredholm operator $t\in \reg{\md E,\md F}$, the Fredholm index is defined via its bounded transform: 
				\begin{equation*}
					\ind{t}:=\ind{F_t}\in \K.
				\end{equation*}
				If $T\in \adj{\md E,\md F}$ is $\alg A$-Fredholm, then its bounded transform $F_T=TQ_T$ satisfies the equation $\ind{F_T}=\ind{T}$ because $Q_T$ is invertible (see \Cref{thm:index-map-adjointable-operators}\ref{item5:thm:index-map-adjointable-operators}). In other words, the index map for regular operators is well defined.
				
				Recall the definition of the Fredholm index for an $\alg A$-Fredholm complex $\complex{t}$ from \Cref{def:Fredholm-index-complex}:
				\begin{equation*}
					\ind{\complex{t}}:=\ind{\evDirac{t}}\in \K.
				\end{equation*}
				For an $\alg A$-Fredholm quasicomplex $\complex{S}$, we similarly define 
				\begin{equation*}
					\ind{\complex{S}}:=\ind{\evDirac{S}}\in \K.
				\end{equation*}

				Now let $\complex{t}$ be any $\alg A$-Fredholm complex. By definition, the index of $\evDirac{t}$ and $\complex{t}$ are equal. Moreover, the index of $\Laplace{t}$ is zero because the operator is self-adjoint. We want to investigate the index of the remaining objects in \Cref{thm:fredholm-goal-I}. This is the content of the following theorem and its corollary.
												
				\begin{theorem}\label{thm:consistency-index-bounded-transform-hilbertized-complex}
					If $\complex{t}$ is an $\alg A$-Fredholm complex, then
					\begin{equation*}
						\ind{\complex{t}}=\ind{\complex{F_t}}=\ind{\complex[\md E_{\graph{t}}]{t}}
					\end{equation*}
					\begin{proof}
						By \Cref{thm:bounded-even-dirac-operator}, we know that $\evDirac{F_t}$ is the bounded transform of $\evDirac{t}$. Moreover, by \Cref{thm:dirac-operator-graph-complex}, there exist invertible operators $U$ and $V$ such that $\evDirac{F_t}=V\evDirac{t_\Gamma}U$. The result follows by definition and \Cref{thm:index-map-adjointable-operators}\ref{item5:thm:index-map-adjointable-operators}.
					\end{proof}
				\end{theorem}	
				
				\begin{corollary}\label{thm:index-adjoint}
					If $\complex{t}:=\{t_k\in \reg{\md E_k,\md E_{k+1}}\}_{k=0}^N$ is an $\alg A$-Fredholm complex, then
					\begin{equation*}
						\ind{\complex{t^\sharp}}=(-1)^{N+1}\ind{\complex{t}}.
					\end{equation*}
					\begin{proof}
						Since $\complex{F_t^\sharp}$ is the bounded transform of $\complex{t^\sharp}$, then the index of $\complex{t^\sharp}$ equals the index of $\evDirac{F_t^\sharp}$. The result follows by \Cref{thm:dirac-operator-adjoint-complex} and \Cref{thm:index-map-adjointable-operators}\ref{item1:thm:index-map-adjointable-operators}.
					\end{proof}
				\end{corollary}	

        \subsection{Fredholm complexes with weak Hodge decomposition}\label{subsec:fredholm_complex_with_WHD}
				
				As mentioned before, there exist $\alg A$-Hilbert complexes without weak Hodge decomposition. We mention here an example, showing that even C*-\emph{Fredholm} complexes can fail to have weak Hodge decomposition (let alone strong Hodge decomposition).
				
				\begin{example}\label{example:Fredholm-no-Hodge}
					Consider the C*-algebra of continuous functions $C([0,1])$. Observe that it can also be regarded as a trivial Hilbert $C([0,1])$-module. In particular, all adjointable operators on $C([0,1])$ are compact because the C*-algebra is unital. Consequently, every adjointable operator is $C([0,1])$-Fredholm. A well-known example of a C*-Fredholm operator without weak Hodge decomposition (polar decomposition) is the multiplication operator $m$, given by $(mf)(x) := x f(x)$, see, e.g., \cite[p.292]{Exel1993}. 
				\end{example}
				
				Motivated by the classical case, we consider the following alternative notion of a Fredholm index. 
				\begin{definition}                                                
					A \emph{weakly $\alg A$-Fredholm} complex is a finite-length $\alg A$-Hilbert complex $\complex{t}$ with weak Hodge decomposition such that, for all $k$, the Hilbert $\alg A$-module $\ker {\Laplace{t}_k}$ is a finite-rank module. The \emph{weak index} of $\complex{t}$ is then defined by
					\begin{equation*}
						\weakind{\complex{t}}:=\Kclass{\ker{\evDirac{t}}}-\Kclass{\ker{\oddDirac{t}}}\in \K.
					\end{equation*}
					Note that $\weakind{\complex{t}}$ is well defined by the assumption on $\ker {\Laplace{t}_k}$ and \Cref{thm:kernel-range-laplace-operator}\ref{item2:thm:kernel-range-laplace-operator}. 
				\end{definition}
				
				We note that not every $\alg A$-Fredholm complex is also weakly $\alg A$-Fredholm (indeed, it might not have weak Hodge decomposition, see \Cref{example:Fredholm-no-Hodge}). 
				Conversely, not every weakly $\alg A$-Fredholm complex is also $\alg A$-Fredholm (indeed, even on Hilbert spaces, a self-adjoint operator can have finite-dimensional kernel without being Fredholm). 

				The following results are straightforward consequences of \Cref{thm:isomorphic_regular_complexes} and \Cref{thm:kernel-range-laplace-operator}\ref{item1:thm:kernel-range-laplace-operator}.
				\begin{itemize}
					\item If $\complex{t}$ has weak Hodge decomposition, then the $K$-class of $\ker{ \Laplace{t}_k}$  is invariant under isomorphism of complexes with weak Hodge decomposition. Therefore, isomorphic weakly C*-Fredholm complexes have the same weak index.
					\item  A complex is weakly C*-Fredholm if and only if its bounded transform is weakly C*-Fredholm. Moreover, they have the same weak index. 
				\end{itemize}
				
				\begin{lemma}\label{thm:fredholm_complex_with_laplace_finite-rank-kernel}
					Let $\complex{t}$ be an $\alg A$-Hilbert complex. Assume that, for every $k$, the module  $\ker{\Laplace{t}_k}$ is complementable. If $\complex{t}$ is $\alg A$-Fredholm, then $\ker{\Laplace{t}_k}$ is a finite-rank module. 
                        \begin{proof}
						Without loss of generality, by \Cref{thm:kernel-range-laplace-operator} and \Cref{thm:fredholm-equivalence-complex-and-bounded-transform}, we prove the statement for the bounded transform complex $\complex{F_t}$. Fix $k$. Since $\Laplace{F_t}_k$ is $\alg A$-Fredholm and $\md E_k$ is orthogonally decomposable along $\ker{\Laplace{F_t}_k}$, we derive that each summand in
						\begin{equation*}
							\Laplace{F_t}_k=\Laplace{F_t}_k|_{\ker{\Laplace{F_t}_k}}\oplus \Laplace{F_t}_k|_{(\ker{\Laplace{F_T}_k})^\perp}
						\end{equation*}
						is an $\alg A$-Fredholm operator. Thus, $0=\Laplace{F_t}_k|_{\ker{ \Laplace{F_t}_k}}\colon  \ker {\Laplace{F_t}_k}\to \ker{\Laplace{F_t}_k}$ is also $\alg A$-Fredholm whence $1_{\ker{\Laplace{F_t}_k}}$ is compact.
					\end{proof}
				\end{lemma}

				\begin{theorem}\label{thm:index_weak_hodge_decomposition}
					If $\complex{t}$ is an $\alg A$-Fredholm complex with weak Hodge decomposition, then $\complex{t}$ is also weakly $\alg A$-Fredholm, and 
					\begin{equation*}
						\ind{\complex{t}}=\weakind{\complex{t}}\in \K.
					\end{equation*}
					\begin{proof}
						Without loss of generality, by \Cref{thm:kernel-range-laplace-operator}, \Cref{thm:weak-hodge-decomposition-derived-complexes}, and \Cref{thm:consistency-index-bounded-transform-hilbertized-complex}, it suffices to prove the statement for the bounded transform complex $\complex{F_t}$. By \Cref{thm:fredholm-equivalence-complex-and-dirac-operator}, $\evDirac{F_t}$ is $\alg A$-Fredholm. Moreover, the weak Hodge decomposition of $\complex{F_t}$ implies that $\evDirac{F_t}$ is polar decomposable (\Cref{thm:weak-hodge-decomposition-equivalences}). By the arguments in \cite[Lemma 17.1.4]{WeggeOlsen1993}, we can show that there exists an $\alg A$-Fredholm operator $G$ with closed range such that $G\sim \evDirac{F_t}$ and
						\begin{equation*}
							\ind{G}=\Kclass{\ker{\evDirac{F_t}}}-\Kclass{\ker{\oddDirac{F_t}}}\in \K.
						\end{equation*}
						The result follows by \Cref{thm:index-map-adjointable-operators}\ref{item3:thm:index-map-adjointable-operators}.
					\end{proof}
				\end{theorem}						
				
				Since complexes with strong Hodge decomposition clearly have weak Hodge decomposition, we can use the previous theorem to compute the index of C*-Fredholm complexes with strong Hodge decomposition similarly to the Hilbert-space case, cf. \eqref{eq3:subsec:C*-hilbert-complexes}. 
				
				\begin{corollary}\label{thm:index_topological_complex}
					Let $\complex{t}=\{t_k\in \reg{\md E_k,\md E_{k+1}}\}_{k=0}^N$  be an $\alg A$-Fredholm complex with strong Hodge decomposition. 
					Then
					\begin{equation}\label{eq1:thm:topological_fredholm_implies_boldfredholm}
						\ind{\complex{t}}=\sum_{k=0}^{N+1} (-1)^k\Kclass{\coH{k}{\complex{t}}}\in \K.
					\end{equation}
					\begin{proof}
						By \Cref{thm:index_weak_hodge_decomposition}, we have $\ind{\complex{t}} = \Kclass{\ker{\evDirac{t}}}-\Kclass{\ker{\oddDirac{t}}}$. 
						Using \Cref{thm:kernel-range-laplace-operator}\ref{item2:thm:kernel-range-laplace-operator}, we can rewrite this as  $\ind{\complex{t}} = \sum_k \Kclass{\ker{\Laplace{t}_{2k}}} - \sum_k \Kclass{\ker{\Laplace{t}_{2k+1}}}$. Finally, since $\complex{t}$ has strong Hodge decomposition, we have a Hilbert $\alg A$-module isomorphism $\ker {\Laplace{t}_k}\cong \coH{k}{\complex{t}}$ for each $k$ (\Cref{coro:cohomology-kernel}), and the statement follows.
					\end{proof}
				\end{corollary}
				
				We finally observe that our definition of the Fredholm property is consistent with the classical theory of Fredholm complexes.
				\begin{theorem}\label{thm:fredholm_complexes_of_hilbert_spaces}
					Assume that $\complex[\md H]{d}=\{d_k\}_{k=0}^N$ is a finite-length Hilbert complex. Then $\complex[\md H]{d}$ is $\mathbb C$-Fredholm if and only if each $k$-th cohomology group $\coH{k}{\complex[\md H]{d}}$ is finite-dimensional.  In this case,
					\begin{equation*}
						\ind{\complex[\md H]{d}}=\sum_{k=0}^{N+1} (-1)^k\dim \coH{k}{\complex[\md H]{d}}\in \mathbb Z,
					\end{equation*}
					and $\complex[\md H]{d}$ has strong Hodge decomposition.
					\begin{proof}
						Without loss of generality,  we might consider the bounded transform complex $\complex[\md H]{F_d}$. It suffices to show that $\complex[\md H]{F_d}$ has a parametrix if and only if the cohomology groups are finite dimensional, which is proven in \cite[Theorem 3, Sec.3.2.3.1]{RempelSchulze1982}. To compute the index, use Atkinson's theorem and \Cref{thm:index_topological_complex}.
					\end{proof}
				\end{theorem}

        \subsection{Putinar's construction of the index}\label{subsec:putinar_index}	
				For Banach spaces, Putinar defined in \cite[Section 2]{Putinar1982} an index map without recurring to cohomology groups, cf. \cite[Proposition 5.2]{Segal1970}. 
				For Hilbert C*-modules, Putinar's construction of the index is used in \cite{Bayramov2008}, which is the only place where we believe (adjointable) Fredholm complexes of (standard) Hilbert C*-modules (over unital C*-algebras) have been studied previously. Our work has been developed without being aware of \cite{Bayramov2008} (instead, our main references were  \cite{Segal1970,Vasilescu1979,BruningLesch1992,SchulzeTarkhanov1998}). 
				In this subsection, we will show that Putinar's construction of the index can be adapted to our setting, and that, for quasicomplexes, it agrees with our definition of the Fredholm index. 
				
				The following proof is in the preamble of \cite[Definition 2.4]{Putinar1982}.
				\begin{lemma}\label{thm:fredholm_equivalence_complex_TevTodd}
					Let $\complex{T}$ be a finite-length quasicomplex of Hilbert $\alg A$-modules. The quasicomplex $\complex{T}$ is $\alg A$-Fredholm if and only if there exists a quasicomplex $\complex P$ such that the operators
					\begin{align*}
						T_{\mathrm{ev}}:=\begin{pmatrix}
							T_0&P_1&0&0\\
							0&T_2&P_3&0\\
							0&0&T_4&\cdots\\
							0&0&0&\ddots
						\end{pmatrix}\quad\text{and}
						\quad
						T_{\mathrm{odd}}:=\begin{pmatrix}
							P_0&0&0&0\\
							T_1&P_2&0&0\\
							0&T_3&P_4&0\\
							0&0&\vdots&\ddots
						\end{pmatrix}
					\end{align*}									
					are $\alg A$-Fredholm in $\adj{\md E_{\mathrm{ev}},\md E_{\mathrm{odd}}}$ and $\adj{\md E_{\mathrm{odd}},\md E_{\mathrm{ev}}}$, respectively.  									
					\begin{proof}
						Let $\complex{T}:=\{T_k\in \adj{\md E_k,\md E_{k+1}}\}_{k=0}^N$ be an $\alg A$-Fredholm quasicomplex. By \Cref{thm:quasicomplex-parametrix}, there exists a parametrix $\complex P$ of $\complex{T}$ that is a quasicomplex. Use this parametrix to define $T_{\text{ev}}$ and $T_{\text{odd}}$. One can show that $T_{\text{odd}}$ is a parametrix of $T_{\text{ev}}$ because $\complex{P}$ is a quasicomplex parametrix. Conversely, assume there exists a quasicomplex (not necessarily a parametrix) $\complex P$ such that both $T_{\text{ev}}$ and $T_{\text{odd}}$ are $\alg A$-Fredholm. Thus, the products $T_{\text{ev}}T_{\text{odd}}$ and $T_{\text{odd}}T_{\text{ev}}$ induce, for every $k$, an $\alg A$-Fredholm operator $P_kT_k+T_{k-1}P_{k-1}$. The result follows by \Cref{thm:associated-quasicomplex}.
					\end{proof}
				\end{lemma}
				
				\begin{remark}\label{rmk:index-T_ev}
					The operator $T_{\text{ev}}$ in \Cref{thm:fredholm_equivalence_complex_TevTodd} arising from a quasicomplex para\-metrix can be found in \cite[Section 2]{Putinar1982} and \cite[Section 5]{SchulzeTarkhanov1998} under different settings, and the index of a quasicomplex is defined as the index of $T_{\text{ev}}$. One needs to show that the index of $T_{\text{ev}}$ does not depend on the choice of the quasicomplex parametrix. This can be done exactly as in the case of Banach and Fréchet spaces, so the proof is omitted and we refer the reader to \cite[Definition 2.4]{Putinar1982} or \cite[Theorem 5.1.9]{SchulzeTarkhanov1998}. 
				\end{remark}
				
				\begin{lemma}\label{thm:index_with_adjoint_of_parametrix}
					For any $\alg A$-Fredholm operator $T\in \adj{\md E,\md F}$ with parametrix $P\in \adj{\md F,\md E}$, we get the identity
					\begin{equation*}
						\ind{T}=\ind{T+P^*}\in \K.
					\end{equation*}
					\begin{proof}
						Since $P$ is a parametrix of $T$, observe that $P^*(1_{\md E}+T^*T)\sim T+P^*$. Moreover, we clearly have that $1_{\md E}+T^*T$ is invertible. The result follows by an application of \Cref{thm:index-map-adjointable-operators} items \ref{item1:thm:index-map-adjointable-operators} and \ref{item3:thm:index-map-adjointable-operators}-\ref{item5:thm:index-map-adjointable-operators}:
						\begin{equation*}
							\ind{T+P^*}=\ind{P^*(1_{\md E}+T^*T)}=\ind{P^*}=\ind{T}.\qedhere
						\end{equation*}
					\end{proof}
				\end{lemma}
				
				\begin{lemma}\label{thm:alternative_fredholm_index_putinar}
					Let $\complex{T}$ be a finite-length $\alg A$-Fredholm quasicomplex. 
					\begin{enumerate}[(i)]
						\item \label{item1:thm:alternative_fredholm_index_putinar}
						If $\complex{P}$ is a quasicomplex parametrix of $\complex{T}$, then 
						\begin{equation*}
							\evDirac{T}+\evDirac{P^\sharp}=T_{\text{ev}}+T_{\text{odd}}^*,
						\end{equation*}
						where $\evDirac{T}$ and $\evDirac{P^\sharp}$ are the even Dirac operators of $\complex{T}$ and $\complex{P^\sharp}$, respectively, and $T_{\mathrm{ev}}$ and $T_{\mathrm{odd}}$ are as in \Cref{thm:fredholm_equivalence_complex_TevTodd}. 
						\item \label{item2:thm:alternative_fredholm_index_putinar}
						There exists a quasicomplex parametrix $\complex{P}$ such that $\oddDirac{P^\sharp}$ is a parametrix of $\evDirac{T}$. 
					\end{enumerate}
					\begin{proof}
						The equality in \ref{item1:thm:alternative_fredholm_index_putinar} follows from a straightforward computation. 

						To prove \ref{item2:thm:alternative_fredholm_index_putinar}, we consider the $k$-th Laplace operator $\Delta_k$ of $\complex{T}$. By \Cref{thm:fredholm-goal-I_quasicomplexes}, the operator $\Delta_k$ is $\alg A$-Fredholm. Assume $R_k$ is a parametrix of $\Delta_k$. Since $\Delta_k$ is self-adjoint, we have that $R_k^*$ is also a parametrix of $\Delta_k$ whence $R_k\sim R^*_k$ for every $k$ (\Cref{thm:parmetrices-fredholm-operators}). Put $\Delta_{\text{ev}}=\oplus_{k} \Delta_{2k}$ and $\Delta_{\text{odd}}=\oplus_{k} \Delta_{2k+1}$, and define $R_{\text{ev}}$ and $R_{\text{odd}}$ in a similar fashion. Thus $R_{\text{ev}}$ and $R_{\text{odd}}$ are parametrices of $\Delta_{\text{ev}}$ and $\Delta_{\text{odd}}$, respectively. Moreover, by the proof of \Cref{thm:associated-quasicomplex}, the collection $\complex{P}:=\{P_k:=R_kT_k^*\}_{k}$ defines a parametrix of $\complex{T}$ and, for every $k$, we obtain $R_{k+1}T_{k}\sim T_{k}R_k$. All together, we have the following relations.
						\begin{align}
							\Delta_{\text{ev}}\sim \oddDirac{T}\evDirac{T}\quad & \text{and}\quad\Delta_{\text{odd}}\sim\evDirac{T}\oddDirac{T}, \label{eq1bis:thm:consistency_index_putinars}\\
							\Delta_{\text{ev}}R_{\text{ev}}\sim 1_{\md E_{\text{ev}}}\quad &\text{and}\quad R_{\text{odd}}\Delta_{\text{odd}}\sim 1_{\md E_{\text{odd}}},\label{eq2:thm:consistency_index_putinars}\\
							R_k&\sim R^*_k, \quad \forall k,\label{eq3:thm:consistency_index_putinars}\\
							R_{k+1}T_k&\sim T_kR_k,\quad \forall k.\label{eq4:thm:consistency_index_putinars}
						\end{align}
						We note that $\complex{P}$ is a quasicomplex: indeed, taking adjoints in \eqref{eq4:thm:consistency_index_putinars} and using  \eqref{eq3:thm:consistency_index_putinars}, we obtain $P_kP_{k+1}=R_kT_k^*R_{k+1}T_{k+1}^*\sim R_k^2T_{k}^*T_{k+1}^*\sim 0$.

						From \ref{item1:thm:alternative_fredholm_index_putinar} and \Cref{thm:index_with_adjoint_of_parametrix} (noting that $T_{\text{odd}}$ is a parametrix for $T_{\text{ev}}$), we know that $\evDirac{T}+\evDirac{P^\sharp}$ is $\alg A$-Fredholm, where $\evDirac{P^\sharp}$ is the even Dirac operator of the quasicomplex $\complex{P^\sharp}$. Once again, by \eqref{eq3:thm:consistency_index_putinars} and \eqref{eq4:thm:consistency_index_putinars}, we derive that
						\begin{align*}
							R_{\text{odd}}\evDirac{T}=
							\begin{pmatrix}
								R_1T_0&R_1T_1^*&0&\cdots\\
								0&R_3T_2&R_3T_3^*&\cdots\\
								0&0&R_5T_4&\cdots\\
								0&0&0&\ddots
							\end{pmatrix}\sim 
							\begin{pmatrix}
								P_0^*&P_1&0&\cdots\\
								0&P_2^*&P_3&\cdots\\
								0&0&P_4^*&\cdots\\
								0&0&0&\ddots
							\end{pmatrix}=\evDirac{P^\sharp},
						\end{align*}
						and, that $\evDirac{T}R_{\text{ev}}\sim \evDirac{P^\sharp}$. Multiply the first (second) equivalence in \eqref{eq1bis:thm:consistency_index_putinars} by $R_{\text{ev}}$ (by $R_{\text{odd}}$) from the right (from the left) and use \eqref{eq2:thm:consistency_index_putinars} to obtain:
						\begin{equation*}
							\evDirac{P^\sharp}\oddDirac{T}\sim 1_{\md E_{\text{odd}}}\quad\text{and}\quad \oddDirac{T}\evDirac{P^\sharp}\sim 1_{\md E_{\text{ev}}}.
						\end{equation*}
						We then observe (by taking adjoints) that $\oddDirac{P^\sharp}$ is a parametrix of $\evDirac{T}$. 
					\end{proof}
				\end{lemma}
				
				Finally, we show that the index of an  $\alg A$-Fredholm quasicomplex using the even Dirac operator, as suggested in \cite{Vasilescu1979,BruningLesch1992}, coincides with the index of a quasicomplex using the operator $T_{\text{ev}}$, as suggested in \cite{Putinar1982}. 
				
				\begin{theorem}\label{thm:consistency_index_putinars}
					Let $\complex{T}$ be an $\alg A$-Fredholm quasicomplex with quasicomplex para\-metrix $\complex{P}$, and put $T_{\mathrm{ev}}$ as constructed in \Cref{thm:fredholm_equivalence_complex_TevTodd}. Then
					\begin{equation*}
						\ind{\evDirac{T}}=\ind{T_{\mathrm{ev}}}\in \K,
					\end{equation*}
					where $\evDirac{T}$ is the even Dirac operator of $\complex{T}$.
					\begin{proof}										
						The index of $T_{\text{ev}}$ is independent of the parametrix; see \Cref{rmk:index-T_ev}. Thus, by \Cref{thm:alternative_fredholm_index_putinar}\ref{item2:thm:alternative_fredholm_index_putinar} we may choose the quasicomplex parametrix $\complex{P}$ such that $\oddDirac{P^\sharp}$ is a parametrix of $\evDirac{T}$. 
						From \Cref{thm:index_with_adjoint_of_parametrix} and \Cref{thm:alternative_fredholm_index_putinar}\ref{item1:thm:alternative_fredholm_index_putinar} we then obtain 
						\begin{equation*}
							\ind{\evDirac{T}} = \ind{\evDirac{T}+\evDirac{P^\sharp}} = \ind{T_{\text{ev}}+T_{\text{odd}}^*} = \ind{T_{\text{ev}}} .
							\qedhere 
						\end{equation*}
					\end{proof}
				\end{theorem}													

	\section{Stability of the index}\label{sec:stability}
		\subsection{Fredholm quasicomplexes}
				We now prove stability properties of the index for Fredholm quasicomplexes. Observe that \Cref{thm:properties_index_quasicomplex} below (and its corollary) holds for adjointable C*-Hilbert complexes if we restrict the statement to small and compact perturbations that are also adjointable complexes.

				\begin{theorem}\label{thm:properties_index_quasicomplex}	
					Let $\complex{T}:=\{T_k\in \adj{\md E_k,\md E_{k+1}}\}_{k=0}^N$ be an  $\alg A$-Fredholm quasicomplex, and let $\complex{S}$ be a finite-length $\alg A$-Hilbert quasicomplex.
					\begin{enumerate}[(i)]
						\item The adjoint quasicomplex $\complex{T^\sharp}$ is $\alg A$-Fredholm with index
						\begin{equation*}
							\ind{\complex{T^\sharp}}=(-1)^{N+1}\ind{\complex{T}}.
						\end{equation*}
						Moreover, the same is true for any of its quasicomplex parametrices.
						\item There exists an $\epsilon>0$ such that if $\complex{S}$ satisfies that $\sup_{k}\norm{T_k-S_k}<\epsilon$, then $\complex{S}$ is an $\alg A$-Fredholm quasicomplex.
						\item If $\complex{S}$ is a compact perturbation of $\complex{T}$, then $\complex{S}$ is an $\alg A$-Fredholm quasicomplex.
					\end{enumerate}
					If $\complex{S}$ satisfies (ii) or (iii), then
					\begin{equation*}
						\ind{\complex{S}}=\ind{\complex{T}}.
					\end{equation*}
					\begin{proof}
						Statements (ii) and (iii) follow by \Cref{thm:index-map-adjointable-operators} items \ref{item3:thm:index-map-adjointable-operators} and \ref{item7:thm:index-map-adjointable-operators}, and the definition of the Fredholm index; in particular, for small perturbations, notice that the distance between the even Dirac operators can be controlled by the distance between the differential maps. 
						
						For (i), to compute the index of a quasicomplex parametrix, we proceed as follows. If $\complex{T}$ is $\alg A$-Fredholm with quasicomplex parametrix $\complex{P}$, then $\complex{S}:=\complex{P^\sharp}$ is $\alg A$-Fredholm with quasicomplex parametrix $\complex{T^*}$. By definition (see \Cref{thm:fredholm_equivalence_complex_TevTodd}), $S_{\mathrm{ev}}=T^*_{\mathrm{odd}}$ and $S_{\mathrm{odd}}=T_{\mathrm{ev}}^*$, and consequently 
                        \(
                            T_{\mathrm{ev}}+T^*_{\mathrm{odd}}=S_{\mathrm{ev}}+S^*_{\mathrm{odd}}.
                        \)                  
                        Moreover, $T_{\mathrm{odd}}$ is a parametrix of $T_{\mathrm{ev}}$, and similarly $S_{\mathrm{odd}}$ is a parametrix of $S_{\mathrm{ev}}$.  Thus, by \Cref{thm:index_with_adjoint_of_parametrix} and \Cref{thm:consistency_index_putinars}, we have $\ind{\complex{T}} = \ind{T_{\mathrm{ev}}+T^*_{\mathrm{odd}}}$ and $\ind{\complex{P^\sharp}} = \ind{S_{\mathrm{ev}}+S^*_{\mathrm{odd}}}$. Hence we see that $\ind{\complex{P^\sharp}}=\ind{\complex{T}}$. 
                        The statement (i) then follows from \Cref{thm:dirac-operator-adjoint-complex} and \Cref{thm:index-map-adjointable-operators}\ref{item1:thm:index-map-adjointable-operators}.
					\end{proof}
				\end{theorem}			

				\begin{corollary}\label{thm:homotopical_Fredholm_quasicomplexes}
					Two finite-length $\alg A$-Hilbert quasicomplexes $\complex{T}$ and $\complex{S}$ are called \emph{homotopic} if there exists a continuous family (with respect to the operator norm) of $\alg A$-Fredholm quasicomplexes $\{\complex{T^{(\lambda)}}\}_{0\leq \lambda\leq 1}$ such that
					\begin{equation*}
						\complex{T}=\complex{T^{(0)}}\quad\text{and}\quad\complex{S}=\complex{T^{(1)}}.
					\end{equation*}
					In this case, $\complex{T}$ and $\complex{S}$ have the same index. 
				\end{corollary}							

        \subsection{``Small'' perturbations} 
				Since the graph of a regular operator in $\reg{\md E,\md F}$ is complementable in $\md E\oplus \md F$ \cite[Theorem 9.3]{Lance1995}, there is a natural metric on the set of regular operators, the so-called \emph{gap metric}. A priori, we define this metric as follows:
				\begin{equation*}
					\tilgapmet{s}{t}:=\norm{\mathscr P_{\graph{t}}-\mathscr P_{\graph{s}}},\quad \forall t,s\in \reg{\md E,\md F},
				\end{equation*}
				where $\mathscr P_{\graph{t}}$ and $\mathscr P_{\graph{s}}$ denote the orthogonal projection from $\md E\oplus \md F$ onto $\graph{t}$ and $\graph{s}$, respectively. The induced topology is called the \emph{gap topology}. One can define an equivalent metric for this topology using the bounded transform:
				\begin{multline*}
					\gapmet{s}{t} := \left( \norm{Q_t^2-Q_s^2}^2+\norm{Q_{t^*}^2-Q_{s^*}^2}^2+2\norm{F_tQ_t-F_sQ_s}^2 \right)^{1/2}, \\ \forall t,s\in \reg{\md E,\md F}.
				\end{multline*}
				This definition is not common in the literature, see, e.g., \cite{Joachim2003,Sharifi2011}. However, it will be convenient to study the stability of the index, cf. \cite[Equation 20]{CordesLabrousse1963}.		
				
				Observe that we have reduced the study of C*-Fredholm complexes to C*-Fredholm operators. Thus, following \cite[Definition 5.2]{CordesLabrousse1963},  let us define the operator
				\begin{equation*}
					V_{t,s}:=Q_tQ_s+F_t^*F_s, \quad \forall t,s\in \reg{\md E,\md F}.
				\end{equation*}
				In the next statement, we consider \cite[Lemmata 5.3 and 5.4]{CordesLabrousse1963}, where most of the arguments can be applied directly.

				\begin{lemma}\label{thm:the_operator_V}
					Consider regular operators $t,s\in \reg{\md E,\md F}$. Then
					\begin{equation*}
						V_{t,t}=1_{\md E}\quad\text{and}\quad V_{t,s}=V_{s,t}^*,
					\end{equation*}
					and
					\begin{equation*}
						|\norm{V_{t,s}x}^2-\norm{x}^2|\leq\norm{x}^2\gapmet{t}{s},\quad \forall x\in \md E.
					\end{equation*}
					In particular, if $\gapmet{t}{s}<1$, then $V_{t,s}$ is invertible and
					\begin{equation*}
						\norm{F_sV_{s,t}-F_t}\leq \gapmet{t}{s}.
					\end{equation*}
					\begin{proof}
						By \Cref{thm:properties-bounded-transform}, recall that $F_t=tQ_t$ and $Q_t^2=1_{\md E}-F_t^*F_t$. Thus,
						\begin{equation*}
							V_{t,t}=Q_t^2+F_t^*F_t=1_{\md E}-F_t^*F_t+F_t^*F_t=1_{\md E} .
						\end{equation*}
						To prove $V_{t,s}=V_{s,t}^*$, we use that $Q_t$ and $Q_s$ are self-adjoint. We now show the first inequality. Let $x\in \md E$. The $\alg A$-valued \airquotes{norm} of $x\in \md E$ is given by $|x|^2=\ev{x}{x}$, see, e.g., \cite[p.4]{Lance1995}. Following \cite[Lemma 5.3]{CordesLabrousse1963}, we show that
						\begin{align*}
							\norm[\alg A]{|V_{t,s}x|^2-|x|^2}\leq \norm[\md E]{x}^2\; \gapmet{t}{s}.
						\end{align*}
						The result is a consequence of the reverse triangle inequality and the identity $\norm[\alg A]{|x|^2}=\norm[\md E]{x}^2$ for every $x\in \md E$. Now, let us assume that $\gapmet{t}{s}<1$. Thus,
						\begin{equation*}
							\norm{V_{t,s}x}^2\geq (1-\gapmet{t}{s})\norm{x}^2,\quad \forall x\in \md E,
						\end{equation*}
						whence $V_{t,s}$ is injective with closed range. Moreover, since $\norm{V_{t,s}^*x}^2=\norm{V_{s,t}x}^2$, then both $V_{t,s}$ and $V_{t,s}^*=V_{s,t}$ are injective with closed range. By \Cref{thm:closed_range}, it follows that $V_{t,s}$ and $V_{s,t}$ are bijective. Hence, the operator $V_{t,s}^{-1}$ is adjointable with adjoint $V_{s,t}^{-1}$. For the remaining inequality, see \cite[Lemma 5.4]{CordesLabrousse1963}.	
					\end{proof}
				\end{lemma}
				
				Let us now consider small perturbations of regular operators with respect to the gap metric, the \emph{Riesz metric}, and \emph{displacements by unbounded operators}, cf. \cite[Corollary 3.5]{Vasilescu1979}. 

				\begin{proposition}\label{thm:fredholm_gap_topology}
					Let $t\in \reg{\md E,\md F}$ be $\alg A$-Fredholm. Then there exists an $\epsilon>0$ such that any $s\in \reg{\md E,\md F}$ satisfying $\gapmet{s}{t}<\epsilon$ is $\alg A$-Fredholm with the same index.
					\begin{proof}
						Since $t$ is $\alg A$-Fredholm, so is its bounded transform $F_t$. By stability for adjointable $\alg A$-Fredholm operators (\Cref{thm:index-map-adjointable-operators}\ref{item7:thm:index-map-adjointable-operators}), there exists $\epsilon_1>0$ such that $T\in \adj{\md E,\md F}$ is $\alg A$-Fredholm and $\ind{T}=\ind{F_t}$ whenever $\norm{T-F_t}<\epsilon_1$. Let $s\in \reg{\md E,\md F}$ be such that $\gapmet{s}{t}<\epsilon$, where $\epsilon:=\min\{\epsilon_1,1\}$. Hence, by \Cref{thm:the_operator_V}, the operator $F_sV_{s,t}$ is $\alg A$-Fredholm, and $V_{s,t}$ is invertible. Therefore, $F_s$ is $\alg A$-Fredholm (\Cref{thm:index-map-adjointable-operators}\ref{item5:thm:index-map-adjointable-operators}), and 
						\begin{equation*}
							\ind{F_s}=\ind{F_sV_{s,t}}=\ind{F_t}.
						\end{equation*}
						We conclude that $s$ is $\alg A$-Fredholm with the same index as $t$.
					\end{proof}
				\end{proposition}

				\begin{proposition}\label{thm:fredholm_riesz_topology}
					Define the \emph{Riesz metric} on $\reg{\md E,\md F}$ by setting
					\begin{equation*}
						\Rieszmet{t}{s}:=\norm{F_t-F_s},\quad \forall t,s\in \reg{\md E,\md F}.
					\end{equation*}
					Let $t\in \reg{\md E,\md F}$ be $\alg A$-Fredholm. Then there exists an $\epsilon>0$ such that any $s\in \reg{\md E,\md F}$ satisfying that $\Rieszmet{t}{s}<\epsilon$ is $\alg A$-Fredholm with the same index.
					\begin{proof}
						Since $s$ is Fredholm if and only if $F_s$ is Fredholm, the statement follows by applying \Cref{thm:index-map-adjointable-operators}\ref{item7:thm:index-map-adjointable-operators} to the bounded transform.
					\end{proof}
				\end{proposition}

				We now want to investigate displacements by unbounded operators. Let us first recall a well-known method, see, e.g., \cite[Lemma 2.3]{KaadLesch2012} and \cite[Lemma 1.15]{Sharifi2011}, to reduce some arguments to the case of self-adjoint, regular operators.
				
				\begin{lemma}\label{thm:reduction_to_sa_operators}
					Let $t$ be a closed and semi-regular operator. Put
					\begin{equation*}
						\hat t:=\begin{pmatrix}
							0&t^*\\
							t&0
						\end{pmatrix}.
					\end{equation*}
					The following statements hold.
					\begin{enumerate}[(i)]
						\item The operator $\hat t$ is closed and symmetric. Moreover, $t$ is regular if and only if $\hat t$ is self-adjoint and regular.
						\item If $t,s\in \reg{\md E,\md F}$, then $\gapmet{t}{s}=\gapmet{\hat t}{\hat s}$ and $\Rieszmet{t}{s}=\Rieszmet{\hat t}{\hat s}$. 
					\end{enumerate}
				\end{lemma}

				\begin{theorem}[Kato-Rellich]\label{thm:kato_rellich}
					Consider regular operators $t,s\in \reg{\md E,\md F}$ such that $\dom{t}\subset \dom{s}$ and $\dom{t^*}\subset \dom{s^*}$. Suppose that $\hat s$ is relative $\hat t$-bounded with relative bound less than $1$, that is, there exist $\alpha\in (0,1)$ and $\beta\geq 0$ such that
					\begin{equation*}
						\norm{\hat sx}\leq \alpha\norm{\hat tx}+\beta\norm{x},\quad \forall x\in \dom{\hat t}.
					\end{equation*}
					Then $t+s$ is regular in $\reg{\md E,\md F}$ with domain $\dom{t}$.
					\begin{proof}
						By \cite[Theorem 4.5]{KaadLesch2012} and \Cref{thm:reduction_to_sa_operators}, the operator $\hat t+\hat s$ is self-adjoint  and regular. Since $\hat t+\hat s\subset \widehat{t+s}$, then $\hat t+\hat s=\widehat{t+s}$ because $\widehat{t+s}$ is symmetric. 
					\end{proof}
				\end{theorem}
				
				\begin{proposition}\label{thm:fredholm_unbounded_displacement}
					Let $t\in \reg{\md E,\md F}$ be $\alg A$-Fredholm. Then there exists an $\epsilon>0$ such that if $s\in \reg{\md E,\md F}$ satisfies that $\dom{t}\subset \dom{s}$, $\dom{t^*}\subset \dom{s^*}$, and 
					\begin{equation}\label{eq1:thm:fredholm_unbounded_displacement}
						\norm{\hat s x}\leq \epsilon(\norm{\hat tx}+\norm{x}),\quad x\in \dom{\hat t},
					\end{equation}
					then $t+s\in \reg{\md E,\md F}$ is $\alg A$-Fredholm with $\ind{t+s}=\ind{t}$.
					\begin{proof}
					   Let $s\in \reg{\md E,\md F}$ satisfy the statement for some $0<\epsilon<1$. Then by \Cref{thm:kato_rellich}, the operator $t+s$ belongs to $\reg{\md E,\md F}$ with domain $\dom{t+s}=\dom{t}$. 
                          By \Cref{thm:fredholm_riesz_topology}, there exists an $\epsilon_0\in (0,1)$ such that any $r\in \reg{\md E,\md F}$ with $\Rieszmet{t}{r}<\epsilon_0$ is $\alg A$-Fredholm and it has the same index as $t$. 
                          By \Cref{thm:reduction_to_sa_operators}, it suffices to consider $\Rieszmet{\hat t}{\hat t+\hat s}$. 
                          By the assumption on $\hat s$, we have $\|\hat s(\hat t+i)^{-1}\| \leq 2\epsilon$ and $\|(\hat t+i)^{-1}\hat s\| \leq 2\epsilon$.
                          It then follows from (the proof of) \cite[Proposition A.4]{vdD25_Callias} that $\|F_{\hat t}-F_{\hat t+\hat s}\| \leq 8\epsilon$. 
                        Choosing $\epsilon$ sufficiently small, the statement then follows from \Cref{thm:fredholm_riesz_topology}. 
					\end{proof}
				\end{proposition}
				
				We can extend all these results to C*-Fredholm complexes. This is the content of the following theorem. Observe that, by definition, an $\alg A$-Fredholm complex has finite length. Thus, any other sequence of maps defined over the same collection of Hilbert $\alg A$-modules must have finite length.

				\begin{theorem}[Stability under small perturbations]\label{thm:small_stability_index_complex}
					Let $\complex{t}$ be an $\alg A$-Fredholm complex. There exists an $\epsilon>0$ such that the following holds.
					\begin{enumerate}[(i)]
						\item If $\complex{s}$ is an $\alg A$-Hilbert complex with $\sup_{k}\gapmet{s_k}{t_k}<\epsilon$, then $\complex{s}$ is an $\alg A$-Fredholm complex with index $\ind{\complex{s}}=\ind{\complex{t}}$.
						\item If $\complex{s}$ is an $\alg A$-Hilbert complex with $\sup_{k}\Rieszmet{s_k}{t_k}<\epsilon$, then $\complex{s}$ is an $\alg A$-Fredholm complex with index $\ind{\complex{s}}=\ind{\complex{t}}$.
						\item If $\{s_k\in \reg{\md E_k,\md E_{k+1}}\}_{k\in \geq 0}$ is a sequence of regular operators such that,  for every $k$, we have $\dom{t_k}\subset \dom{s_k}$, $\dom{t_k^*}\subset \dom{s_k^*}$,
						\begin{equation*}
							\norm{\hat s_kx_k}\leq \epsilon(\norm{\hat t_kx_k}+\norm{x_k}),\quad \forall x_k\in \dom{\hat t_k},
						\end{equation*}
						and $\complex{t+s}:=\{t_k+s_k\}_{k}$ is an $\alg A$-Hilbert complex, then $\complex{t+s}$ is an $\alg A$-Fredholm complex with index $\ind{\complex{t+s}}=\ind{\complex{t}}$.
						\item If $\{S_k\in \adj{\md E_k,\md E_{k+1}}\}_{k}$ are adjointable operators such that $\sup_{k}\norm{S_k}<\epsilon$ and $\complex{t+S}:=\{t_k+S_k\}_{k\in \mathbb Z}$ is an $\alg A$-Hilbert complex, then  $\complex{t+S}$ is an $\alg A$-Fredholm complex with index $\ind{\complex{t+S}}=\ind{\complex{t}}$.
					\end{enumerate}
					\begin{proof}
					By \Cref{thm:fredholm-goal-I} and the definition of the Fredholm index of a complex, it suffices to consider the even Dirac operator of the complex. 
                        We observe that the distance between the even Dirac operators is controlled by the distance between the differential maps, see \eqref{eq2:subsec:dirac-operator} and \Cref{thm:bounded-even-dirac-operator}. We then use \Cref{thm:fredholm_gap_topology} for (i), \Cref{thm:fredholm_riesz_topology} for (ii), and \Cref{thm:fredholm_unbounded_displacement} for (iii) and (iv).
					\end{proof}
				\end{theorem}

				We now consider the unbounded version of \Cref{thm:homotopical_Fredholm_quasicomplexes}. Clearly, we can also use the Riesz metric to define a continuous path of C*-Fredholm complexes.

				\begin{corollary}\label{thm:homotopy_fredholm_complexes}
					Two finite-length $\alg A$-Hilbert complexes $\complex{t}$ and $\complex{s}$ are called \emph{homotopic} whenever there exists a continuous family of $\alg A$-Fredholm complexes, with respect to the gap metric, $\complex{t^{(\lambda)}}_{0\leq \lambda\leq 1}$ such that
					\begin{equation*}
						\complex{t}=\complex{t^{(0)}}\quad\text{and}\quad\complex{s}=\complex{t^{(1)}}.
					\end{equation*}
					In this case, $\complex{t}$ and $\complex{s}$ have the same index.
				\end{corollary}
				
		\subsection{Relatively compact perturbations}													
				We now consider the stability under relatively compact perturbations. 
				Given $t\in \reg{\md E,\md F}$, we recall that a densely defined operator $r \colon \md E \to \md F$ is called \emph{relatively $t$-compact} if $\dom{t} \subset \dom{r}$ and $r (1+t^*t)^{-\frac12}$ is compact. 
				We have the following stability result for Fredholm operators. 
				
				\begin{proposition}[{cf.\ \cite[Proposition A.11]{vdD25_Callias}}]\label{prop:Fredholm_rel_cpt_pert}
				Let $t$ be a regular Fredholm operator on $\md E$, and let $r$ be a semi-regular operator on $\md E$ such that $r$ is relatively $t$-compact and $r^*$ is relatively $t^*$-compact. 
				Then: 
				\begin{enumerate}[(i)]
				\item 
				\label{item:rel_cpt_pert_reg_sa}
				$t+r$ is also regular and Fredholm, and any joint parametrix for $t$ is also a joint parame\-trix for $t+r$. 
				\item 
				$\ind{t+r}=\ind{t}$. 
				\end{enumerate}
				\end{proposition}
				
				\begin{theorem}[Stability under relatively compact perturbations]\label{thm:relative-compact-stability-index-complex}
					Let $\complex{t}$ be an $\alg A$-Fredholm complex, and let $\{r_k:\dom{r_k}\subset \md E_{k}\to \md E_{k+1}\}_{k}$ be a sequence of semi-regular operators such that each $r_k$ is relatively $t_k$-compact and each $r_k^*$ is relatively $t_k^*$-compact. 
					If $\complex{t+r}:=\{t_k+r_k\}_{k\in\mathbb Z}$ is an $\alg A$-Hilbert complex, then $\complex{t+r}$ is $\alg A$-Fredholm with index $\ind{\complex{t+r}}=\ind{\complex{t}}$.
					\begin{proof}
						Observe that the operator $\evDirac{r}$ is relatively $\evDirac{t}$-compact. (By abuse of notation, the semi-regular operator $\evDirac{r}$ is given by \eqref{eq3:subsec:dirac-operator}, even though it is not the even Dirac operator of a complex.) 
						Indeed, since $Q_{2k}$ commutes with $Q_{(2k-1)*}$ by \Cref{thm:properties-resolvent-sequence}\ref{item2:thm:properties-resolvent-sequence}, and using the equality (cf.\ \eqref{eq2:subsec:dirac-operator}) $1_{\md E_{\text{ev}}}+\oddDirac{t}\evDirac{t} = \bigoplus_k Q_{2k}^{-2} Q_{(2k-1)^*}^{-2}$, one computes that 
						\[                                                      
							\evDirac{r} Q_{\evDirac{t}} = \bigoplus_{k} r_{2k} Q_{2k} Q_{(2k-1)^*} + r_{2k-1}^* Q_{(2k-1)^*} Q_{2k} ,
						\]
						which is compact. Similarly also $(\evDirac{r})^*$ is relatively $\oddDirac{t}$-compact because $\oddDirac{r}\subset (\evDirac{r})^*$. 
						The result follows by \Cref{prop:Fredholm_rel_cpt_pert}.
					\end{proof}
				\end{theorem}
				
				\begin{corollary}[Stability under compact perturbations]\label{thm:compact-stability-index-complex}
					Let $\complex{t}$ be an $\alg A$-Fred\-holm complex. If $\{C_k\in \comp{\md E_k,\md E_k}\}_{k}$ is a sequence of compact operators such that $\complex{t+C}:=\{t_k+C_k\}_{k\in\mathbb Z}$ is an $\alg A$-Hilbert complex, then $\complex{t+C}$ is $\alg A$-Fredholm with index $\ind{\complex{t+C}}=\ind{\complex{t}}$.
				\end{corollary}

	\section{Direct sums, sequences, and tensor products}\label{sec:sums-seq-tensor}
				In this section, we consider three further constructions: direct sums of complexes, exact sequences of complexes, and tensor products of quasicomplexes. In particular, we investigate their indices.
	
		\subsection{Direct sums of complexes}										
				The \emph{direct sum} of two $\alg A$-Hilbert complexes $\complex{t}$ and $\complex[\md F]{s}$, denoted by $\complex{t}\oplus\complex[\md F]{s}$, is the $\alg A$-Hilbert complex with differential maps
				\begin{equation*}
					t_k\oplus s_k\colon \dom{t_k}\oplus \dom{s_k}\subset \md E_k\oplus \md F_k\to \md E_{k+1}\oplus\md F_{k+1}, \quad \forall k\in \mathbb Z.
				\end{equation*}
				Under this operation, the collection of $\alg A$-Fredholm complexes forms a semigroup.
				
				\begin{theorem}\label{thm:index_direct_sum}
					If $\complex{t}$ and $\complex[\md F]{s}$ are $\alg A$-Fredholm complexes, then its direct sum $\complex{t}\oplus\complex[\md F]{s}$ is an $\alg A$-Fredholm complex. Moreover,
					\begin{equation*}
						\ind{\complex{t}\oplus\complex[\md F]{s}}=\ind{\complex{t}}+\ind{\complex[\md F]{s}}.
					\end{equation*}
					\begin{proof}
						The map $\evDirac{t}\oplus\evDirac{s}$ is the even Dirac operator of $\complex{t}\oplus\complex[\md F]{s}$ with respect to the splitting $\md E_{\text{ev}}\oplus \md F_{\text{ev}}$. The result readily follows by \Cref{thm:index-map-adjointable-operators}\ref{item6:thm:index-map-adjointable-operators} and \Cref{thm:consistency-index-bounded-transform-hilbertized-complex}.
					\end{proof}
				\end{theorem}				
				
		\subsection{Exact sequences of Fredholm complexes}
				The definition of maps of complexes from \Cref{sec:maps_complexes} can be easily extended to quasicomplexes. In fact, if $\complex{T}$ and $\complex[\md F]{S}$ are $\alg A$-Hilbert quasicomplexes, we can say that a map of quasicomplexes $G\colon \complex{T}\to \complex[\md F]{S}$ is a collection of adjointable operators $\{G_k\in \adj{\md E_k,\md F_k}\}_{k\in\mathbb Z}$ such that $G_{k+1}T_k\sim S_kG_k$ for every $k\in \mathbb Z$. With this definition, we state the following theorem.
				
				\begin{theorem}\label{thm:fredholm_maps_complexes_sequence}
					Let $\complex{T}$, $\complex[\md F]{U}$, and $\complex[\md G]{V}$ be finite-length $\alg A$-Hilbert quasicomplexes, and assume there exists a sequence of maps of quasicomplexes
					\begin{equation}\label{eq2:subsec:isomorphic_fredholm_complexes}
						\begin{tikzcd}
							0 & {\complex{T}} & {\complex[\md F]{U}} & {\complex[\md G]{V}} & 0
							\arrow[from=1-1, to=1-2]
							\arrow["G", from=1-2, to=1-3]
							\arrow["H", from=1-3, to=1-4]
							\arrow[from=1-4, to=1-5]
						\end{tikzcd}
					\end{equation}
					such that
					\begin{equation}\label{eq3:subsec:isomorphic_fredholm_complexes}
						\begin{tikzcd}
							0 & {\md E_k} & {\md F_k} & {\md G_k} & 0
							\arrow[from=1-1, to=1-2]
							\arrow["{G_k}", from=1-2, to=1-3]
							\arrow["{H_k}", from=1-3, to=1-4]
							\arrow[from=1-4, to=1-5]
						\end{tikzcd}
					\end{equation}
					is an $\alg A$-Fredholm quasicomplex for every $k$. If at least two of the quasicomplexes $\complex{T}$, $\complex[\md F]{U}$, and $\complex[\md G]{V}$ are $\alg A$-Fredholm, then the third is $\alg A$-Fredholm as well.
					\begin{proof}
						Let $\Sigma$ be a Hilbert $\alg A$-module. Apply the functor $\phi_\Sigma$ to \eqref{eq2:subsec:isomorphic_fredholm_complexes}. By \Cref{thm:properties-putinar's-functor,thm:fredholm-equivalence-putinar-criterion}, the following is a short exact sequence of complexes:
						\[\begin{tikzcd}
							0 & {\complex[\phi_\Sigma(\md E)]{\phi_\Sigma(T)}} & {\complex[\phi_\Sigma(\md F)]{\phi_\Sigma(U)}} & {\complex[\phi_\Sigma(\md G)]{\phi_\Sigma(V)}} & 0.
							\arrow[from=1-1, to=1-2]
							\arrow["{\phi_\Sigma(G)}", from=1-2, to=1-3]
							\arrow["{\phi_\Sigma(H)}", from=1-3, to=1-4]
							\arrow[from=1-4, to=1-5]
						\end{tikzcd}\]
						Moreover, by assumption and \Cref{thm:fredholm-equivalence-putinar-criterion}, two out of the three complexes are exact. By \cite[Exercise 1.3.1]{Weibel1994}, then so is the third. Therefore, by \Cref{thm:fredholm-equivalence-putinar-criterion}, if at least two of the three quasicomplexes are $\alg A$-Fredholm, then so is the third. 
					\end{proof}
				\end{theorem}
				
				As the following theorem suggests, assuming that the sequence \eqref{eq3:subsec:isomorphic_fredholm_complexes} is an $\alg A$-Fredholm quasicomplex is a general version of a classical situation. Let us first consider the following result.
				
				\begin{theorem}\label{thm:exact_complexes_index}
					Let $\complex{T}$ be finite-length $\alg A$-Hilbert complex. If $\complex{T}$ is an exact sequence of Hilbert $\alg A$-modules, then $\complex{T}$ is an $\alg A$-Fredholm complex with index zero.
					\begin{proof}
						By \Cref{thm:fredholm-goal-I} and \Cref{thm:exact-complex-isomorphism-laplacians}, we know that $\complex{T}$ is a $\alg A$-Fredholm complex with strong Hodge decomposition. Moreover, by exactness of the sequence and \Cref{thm:index_topological_complex}, we deduce that $\ind{\complex{T}}=0$.
					\end{proof}
				\end{theorem}
				
				Observe that Putinar's functor $\phi_{\Sigma}$ transforms short exact sequences of Hilbert $\alg A$-modules into short exact sequences of Hilbert $\mathcal{C}(\Sigma)$-modules, for any Hilbert $\alg A$-module $\Sigma$; this an immediate application of \Cref{thm:fredholm-equivalence-putinar-criterion} and \Cref{thm:exact_complexes_index}. We would say that $\phi_{\Sigma}$ is an exact functor if the category of Hilbert C*-modules were an abelian category.

				\begin{corollary}\label{thm:exact-sequence-fredholm-quasicomplexes}
					Let $\complex{T}$, $\complex[\md F]{U}$, and $\complex[\md G]{V}$ be finite-length $\alg A$-Hilbert quasicomplexes, and assume there exists an exact sequence of maps of quasicomplexes
					\begin{equation*}
						\begin{tikzcd}
							0 & {\complex{T}} & {\complex[\md F]{U}} & {\complex[\md G]{V}} & 0.
							\arrow[from=1-1, to=1-2]
							\arrow["G", from=1-2, to=1-3]
							\arrow["H", from=1-3, to=1-4]
							\arrow[from=1-4, to=1-5]
						\end{tikzcd}
					\end{equation*}
					If at least two of the quasicomplexes $\complex{T}$, $\complex[\md F]{U}$, and $\complex[\md G]{V}$ are $\alg A$-Fredholm, then the third is $\alg A$-Fredholm as well.
					\begin{proof}
						It follows directly from \Cref{thm:fredholm_maps_complexes_sequence} and \Cref{thm:exact_complexes_index}.
					\end{proof}
				\end{corollary}
				
				In the Hilbert-space case, see, e.g., \cite[Corollary 1.8]{AmbrozieVasilescu1995}, one can show that for a short-exact sequence of maps of Fredholm complexes, say:
				\begin{equation*}
				\begin{tikzcd}
					0 & {\complex[\md H^{(1)}]{T}} & {\complex[\md H^{(2)}]{U}} & {\complex[\md H^{(3)}]{V}} & 0,
					\arrow[from=1-1, to=1-2]
					\arrow["G", from=1-2, to=1-3]
					\arrow["H", from=1-3, to=1-4]
					\arrow[from=1-4, to=1-5]
				\end{tikzcd}
				\end{equation*}
				their indices satisfy the formula:
				\begin{equation*}
					\ind{\complex[\md H^{(1)}]{T}}-\ind{\complex[\md H^{(2)}]{U}}+\ind{\complex[\md H^{(3)}]{V}}=0.
				\end{equation*}
				One can easily show, using \Cref{thm:exact_complexes_index}, that the same occurs in our case if we assume that the $\alg A$-Hilbert complexes have strong Hodge decomposition. However, we do not know if this holds in general. Let us consider a simpler case.
				
				\begin{theorem}\label{thm:index-unitarily-isomorphic-complexes}
					Let $\complex{T}$ and $\complex[\md F]{S}$ be adjointable $\alg A$-Fredholm complexes which are unitarily isomorphic (i.e., there is a map of complexes $G\colon \complex{T}\to \complex[\md F]{S}$ such that $G_k^*=G_k^{-1}$ for every $k$). Then $\ind{\complex{T}}=\ind{\complex{S}}$.
					\begin{proof}
						One can show that $G_{\mathrm{odd}}\evDirac{T}=\evDirac{S}G_{\mathrm{ev}}$, where $G_{\mathrm{ev}}:=\oplus_k G_{2k}$, and $G_{\mathrm{odd}}$ is defined similarly. The result follows by \Cref{thm:index-map-adjointable-operators}\ref{item5:thm:index-map-adjointable-operators}.
					\end{proof} 
				\end{theorem}													

		\subsection{Tensor product of complexes}\label{sec:tensor_product}
				Finally, we define the tensor product of adjointable complexes in analogy with the classical case, see, e.g., \cite[p.137]{SchulzeTarkhanov1998}. Let $\alg A$ and $\alg B$ be arbitrary C*-algebras. We denote by $\alg A\otimes \alg B$ the completion of $\alg A\otimes_{\mathrm{alg}} \alg B$ with respect to the minimal C*-norm.  For an introduction to the exterior tensor product of Hilbert C*-modules and adjointable operators, see \cite[Chapter 4]{Lance1995}.
				
				Let $\complex[\md L]{R}:=\{R_j\}_{j=0}^N$ be an adjointable finite-length $\alg A$-Hilbert complex, and let $\complex[\md M]{S}=\{S_i\}_{i=0}^M$ be an adjointable finite-length $\alg B$-Hilbert complex. For each $k$, consider
				\begin{equation*}
					\md E_k:=\bigoplus_{j+i=k} \md L_{j}\otimes \md M_{i}
				\end{equation*}
                as (finite) orthogonal sums of Hilbert $\alg A\otimes \alg B$-modules. Put $\sigma_j:=(-1)^j 1_{\md L_j}$ for every $j\geq 0$, and let us define an adjointable operator $T_k$ that is given by
				\begin{equation*}
					T_k:=\bigoplus_{j+i=k} (R_{j} \otimes 1_{\md M_i}+\sigma_j\otimes S_i).
				\end{equation*}		
				Notice that the maps $\{\sigma_j\}_{j\geq 0}$ satisfy the condition $\sigma_{j+1}R_{j}+R_j\sigma_j=0$ for all $j\geq 0$. For each $k$, the complex property of $\complex[\md L] {R}$ and $\complex[\md M]{S}$ yields the relation
				\begin{align*}
					T_{k+1}T_{k}&= \sigma_{j+1}R_j\otimes S_i+R_j\sigma_j\otimes S_i=0
				\end{align*}
				when it is restricted to $\md L_j\otimes_{alg}\md M_i$ for some $i+j=k$. Therefore, $\complex{T}$ is a finite-length $\alg A\otimes \alg B$-Hilbert complex. We write
				\begin{equation*}
					\complex[\md L]{R}\otimes \complex[\md M]{S}:=\complex{T}.
				\end{equation*}
				
				For adjointable complexes of length one (i.e., adjointable operators), the tensor product as complexes is not equal to its tensor product as operators. Indeed, let $R\colon L_0\to L_1$ and $S:M_0\to M_1$ be adjointable operators. 	Observe first that $\evDirac{R}=R$ and $\evDirac{S}=S$.  Hence, their tensor product (as complexes) has the form
				\[\begin{tikzcd}
					0 & {\md L_0\otimes \md M_0} & {(\md L_1\otimes\md M_0)\oplus(\md L_0\otimes \md M_1)} &   {\md  L_1\otimes\md M_1} & 0
					\arrow[from=1-1, to=1-2]
					\arrow["\begin{array}{c} \begin{pmatrix}R\otimes 1_{\md M_0}\\1_{\md L_0}\otimes  S\end{pmatrix}  \end{array}"'{yshift=8ex}, from=1-2, to=1-3]
					\arrow["{\begin{pmatrix}-1_{\md L_1}\otimes S,\; R\otimes 1_{\md M_1}\end{pmatrix}}"'{yshift=6ex}, from=1-3, to=1-4]
					\arrow[from=1-4, to=1-5]
				\end{tikzcd}\]
				with even Dirac operator
				\begin{equation*}
					\evDirac{R}\# \evDirac{S}:=\begin{pmatrix}
						\evDirac{R}\otimes 1_{\md M_{0}}&-1_{\md L_{1}}\otimes \oddDirac{S}\\
						1_{\md L_{0}}\otimes \evDirac{S}& \oddDirac{R}\otimes 1_{\md M_{1}}
					\end{pmatrix}.
				\end{equation*}																
				
				Let us now fix an adjointable finite-length $\alg A$-Hilbert complex $\complex[\md L]{R}:=\{R_j\in \adj{\md L_j,\md L_{j+1}}\}_{j=0}^N$ and an adjointable finite-length $\alg B$-Hilbert complex $\complex[\md M]{S}:=\{S_i\in \adj[\alg B]{\md M_i,\md M_{i+1}}\}_{i=0}^M$. In the following theorem, we show that the tensor product of adjointable C*-Fredholm complexes is again an adjointable C*-Fredholm complex.

            	\begin{theorem}\label{thm:fredholm_property_tensor_product}
					If $\complex[\md L]{R}$ and $\complex[\md M]{S}$ have the Fredholm property, then $\complex{T}:=\complex[\md L]{R}\otimes \complex[\md M]{S}$ is an $\alg A\otimes \alg B$-Fredholm complex.
					\begin{proof}
						Let $\complex[\md L]{\hat P}$ and $\complex[\md M]{\tilde P}$ be parametrices of $\complex[\md L]{R}$ and $\complex[\md M]{S}$, respectively, and consider the compact operators
                        \[
                        \hat C_j = 1 - \hat P_j R_j - R_{j-1} \hat P_{j-1} , 
                        \qquad 
                        \tilde C_i = 1 - \tilde P_i S_i - S_{i-1} \tilde P_{i-1} . 
                        \]
                        We will construct a parametrix as follows (the idea for this construction is borrowed from the proof of \cite[Theorem 2.6]{Eschmeier88}): for each $k$, define the adjointable operator
						\begin{equation*}
							P_k:=\bigoplus_{j+i=k}\hat P_j\otimes 1_{\md M_i}+\sigma_j \hat C_j \otimes \tilde P_i.
						\end{equation*}
					Using the complex property of $\complex[\md L]{R}$, we have $\hat C_{j+1} R_j = R_j - R_j \hat P_j R_j = R_j \hat C_j$. An explicit computation then shows that 
						\begin{equation*}
							P_kT_k+T_{k-1}P_{k-1} = 1_{\md E_k} - \bigoplus_{j+i=k} \hat C_j \otimes \tilde C_i.
						\end{equation*}
                        Indeed, observe that 
                        \begin{align*}
                             P_kT_k&=\hat P_jR_j\otimes 1_{\md M_i}+\sigma_{j+1}\hat C_{j+1}R_j\otimes \tilde P_{i-1}+\hat P_{j-1}\sigma_j\otimes S_i+\hat C_j\otimes \tilde P_iS_i,\\
                            T_{k-1}P_{k-1}&=
                            R_{j-1}\hat P_{j-1}\otimes 1_{\md M_i}+ \sigma_{j-1}\hat P_{j-1}\otimes S_i+R_j\sigma_j\hat C_j\otimes \tilde P_{i-1}+\hat C_j\otimes S_{i-1}\tilde P_{i-1}                        
                        \end{align*}
                        when they are restricted to $\md L_j\otimes_{\mathrm{alg}} \md M_i$ for some $i+j=k$. The remaining computations are strightforward. This proves that $\{P_k\}$ is a parametrix for $\complex{T}$. 
					\end{proof}
				\end{theorem}
                
				\begin{theorem}\label{thm:bold_operator_tensor_product}
					Consider the $\alg A\otimes \alg B$-Hilbert complex $\complex{T}:=\complex[\md L]{R}\otimes \complex[\md M]{S}$, and let $\evDirac{R},\evDirac{S}$, and $\evDirac{T}$ denote the even Dirac operators of the respective complexes. Then there exist unitaries $\mathcal U$ and $\mathcal V$ such that
					\begin{equation*}
						\mathcal V^*\evDirac{T}\mathcal U=
						\evDirac{R}\# \evDirac{S}:=\begin{pmatrix}
							\evDirac{R}\otimes 1_{\md M_{\text{ev}}}&-1_{\md L_{\text{odd}}}\otimes \oddDirac{S}\\
							1_{\md L_{\text{ev}}}\otimes \evDirac{S}&\oddDirac{R}\otimes 1_{\md M_{\text{odd}}}
						\end{pmatrix}.
					\end{equation*}
					\begin{proof}
						Observe that there are isometric surjective $\alg A\otimes \alg B$-linear maps
						\begin{align*}
							\mathcal U_0\colon  (\md L_{\text{ev}}\otimes_{\text{alg}} \md M_{\text{ev}})\oplus (\md L_{\text{odd}}\otimes_{\text{alg}} \md M_{\text{odd}})&\to \bigoplus_{k}\bigoplus_{j+i=2k} \md L_{j}\otimes_{\text{alg}} \md M_i,\\
							\mathcal V_0\colon  (\md L_{\text{odd}}\otimes_{\text{alg}} \md M_{\text{ev}})\oplus (\md L_{\text{ev}}\otimes_{\text{alg}} \md M_{\text{odd}})&\to \bigoplus_{k}\bigoplus_{j+i=2k+1} \md L_{j}\otimes_{\text{alg}} \md M_i.
						\end{align*}
						Thus, they can be extended to unitaries $\mathcal U$ and $\mathcal V$ by \cite[Theorem 3.5]{Lance1995}. Hence, it suffices to show that, for any $i,j,l,m\geq 0$, the operators $\evDirac{T}$ and $\evDirac{R}\# \evDirac{S}$ agree on
						\begin{equation*}
							\Gamma:=(\md L_{2j}\otimes_{\text{alg}} \md M_{2i})\oplus (\md L_{2l+1}\otimes_{\text{alg}}\md  M_{2m+1}).
						\end{equation*}
						Fix $i,j,m,l\geq 0$. A straightforward computation shows that
						\begin{equation*}
							(\evDirac{R}\# \evDirac{S})|_{\Gamma}:=\begin{pmatrix}
								(R_{2j}+R_{2j-1}^*)\otimes 1_{\md M_{2i}}&-1_{\md L_{2l+1}}\otimes (S^*_{2m}+S_{2m+1})\\
								1_{\md L_{2j}}\otimes (S_{2i}+S_{2i-1}^*)&(R_{2l+1}+R_{2l}^*)\otimes 1_{\md M_{2m+1}}
							\end{pmatrix}.
						\end{equation*}
						Moreover, for every $k$, and, up to a reorder of the orthogonal sums,
						\begin{align*}
							T_{2k}&=\bigoplus_{n\geq 0} (R_n\otimes 1_{\md M_{2k-n}}+\sigma_{n}\otimes S_{2k-n}),\\
							T_{2k-1}^*&=\bigoplus_{n\geq 0} (R_n^*\otimes 1_{\md M_{2k-1-n}}+\sigma_{n}\otimes S^*_{2k-1-n}).											 
						\end{align*}
						Thus,
						\begin{align*}
							T_{2(i+j)}|_{\md L_{2j}\otimes_{\text{alg}}\md M_{2i}}&=R_{2j}\otimes 1_{\md M_{2i}}+1_{\md L_{2j}}\otimes S_{2i},\\
							T_{2(j+i)-1}^*|_{\md L_{2j}\otimes_{\text{alg}}\md M_{2i}}&=R_{2j-1}^*\otimes 1_{\md M_{2i}}+1_{\md  L_{2j}}\otimes S_{2i-1}^*,\\
							T_{2(l+m)+2}|_{\md L_{2l+1}\otimes_{\text{alg}}\md M_{2m+1}}&=R_{2l+1}\otimes 1_{\md M_{2m+1}}-1_{\md L_{2l+1}}\otimes S_{2m+1},\\
							T_{2(l+m)+1}^*|_{\md L_{2l+1}\otimes_{\text{alg}}\md M_{2m+1}}&=R_{2l}^*\otimes 1_{\md M_{2m+1}}-1_{\md L_{2l+1}}\otimes S_{2m}^*.
						\end{align*}
						The result readily follows.
					\end{proof}
				\end{theorem}
																				
				\begin{corollary}\label{thm:index_tensor_quasicomplexes}
					If $\complex[\md L]{R}$ and $\complex[\md M]{S}$ have the Fredholm property, then $\complex[\md L]{R}\otimes \complex[\md M]{S}$ is $\alg A\otimes \alg B$-Fredholm and
					\begin{equation*}
						\ind{\complex[\md L]{R}\otimes \complex[\md M]{S}}=\ind{\evDirac{R}\# \evDirac{S}}\in \K[\alg A\otimes  \alg B].
					\end{equation*}
					\begin{proof}
						It follows by \Cref{thm:index-map-adjointable-operators}\ref{item5:thm:index-map-adjointable-operators}.
					\end{proof}
				\end{corollary}

	\section{Examples}\label{sec:applications}
				In this section, we consider two immediate examples of the theory of C*-Hilbert complexes, namely complexes over the C*-algebra of compact operators and C*-elliptic complexes.							
								
		\subsection{Complexes over the C*-algebras of compact operators}		
				Let us denote by $\alg K$ the C*-algebra of \emph{all} compact operators on some separable infinite-dimensional Hilbert space. We aim to prove the following statement.
				
				\begin{theorem}\label{thm:topological_fredholm_compact_operators}
					Let $\complex{t}:=\{t_k\in \reg{\md E_k,\md E_{k+1}}\}_{k=0}^N$ be a finite-length $\alg K$-Hilbert complex. The following conditions are equivalent:
					\begin{enumerate}[(i)]
						\item The complex $\complex{t}$ is $\alg K$-Fredholm.
						\item The complex $\complex{t}$ has strong Hodge decomposition and $\dim_{\alg K}\coH{k}{\complex{t}}<\infty$ for every $k$.
					\end{enumerate}	
					If one of these conditions is satisfied, we have
					\begin{equation*}
						\ind{\complex{t}}=\sum_{k=0}^{N+1}(-1)^k\dim_{\alg K} \coH{k}{\complex{t}}\in \mathbb Z\cong \K[\alg K].
					\end{equation*}
				\end{theorem}	
				
				Here $\dim_{\alg K}$ denotes the \emph{orthonormal dimension} of a Hilbert C*-module as presented in \cite{BakicGuljas2002}, which is defined as follows. Let $\md E$ be a Hilbert $\alg A$-module. A system $(x_{\lambda})_{\lambda\in\Lambda} $ in $\md E$ is orthonormal if each $x_\lambda$ is a basis vector, that is, if $\ev{x_\lambda}{x_\lambda}$ is a minimal projection in $\alg A$ and $\ev{x_\lambda}{x_\mu}=0$ for $\lambda\neq \mu$. An orthonormal system $(x_\lambda)$ is said to be an orthonormal basis if it generates a dense submodule of $\md E$. 
				
				\begin{lemma}[{\cite[Theorem 4 and Section 3]{BakicGuljas2002}}]
					Let $\md E$ be a Hilbert $\alg K$-module. The following statements holds.
					\begin{enumerate}[(i)]
						\item The module $\md E$ has an orthonormal basis. Moreover, any two orthonormal basis of $\md E$ have the same cardinal number, which is denoted by $\dim_{\alg K} \md E$.
						\item  Fix a minimal projection $e_0$ in $\alg K$ and put $\md E_{e_0}:=\{xe_0\mid  x\in \md E\}$. Then $\md E_{e_0}$ is a Hilbert space with respect to the inner product $(\cdot|\cdot):=\mathrm{Tr}(\ev{\cdot}{\cdot})$, where $\mathrm{Tr}$ denotes the (unbounded) trace on $\alg K$, and $\dim_{\alg K}\md E=\dim \md E_{e_0}$. 
					\end{enumerate}
				\end{lemma}
				
				In the following theorem, we show a straightforward application of the results in \cite[Theorems 5 and 6]{BakicGuljas2002} and \cite[Theorem 3.7]{NiknamSharifi2007}, where \Cref{thm:hodge-decomposition-equivalences}, \Cref{thm:fredholm-equivalence-laplace-dirac-operators} and \Cref{thm:C*-like-property-quotient} are useful to prove (ii).
				
				\begin{lemma}\label{thm:Psi_isomoorpism}
					Let $\md E$ and $\md F$ be Hilbert $\alg K$-modules, and let $e_0$ be an arbitrary minimal projection in $\alg K$. The following statements hold.
					\begin{enumerate}[(i)]
						\item The map $\Psi_{\md E,\md F}\colon \adj[\alg K]{\md E,\md F}\to\adj[\mathbb C]{\md E_{e_0},\md F_{e_0}}$, given by $\Psi_{\md E,\md F}(T):=T|_{\md E_{e_0}}$, is a *-homomorphism. Moreover, $\Psi_{\md E,\md E}$ is a *-isomorphism of C*-algebras.\label{item1:thm:Psi_isomoorpism}
						\item An operator $T$ is $\alg K$-Fredholm (is compact) in $\adj[\alg K]{\md E,\md F}$ if and only if the operator $\Psi_{\md E,\md F}(T)$ is $\mathbb C$-Fredholm (is compact) in $\adj[\mathbb C]{\md E_{e_0},\md F_{e_0}}$. Moreover, $T$ is $\alg K$-Fredholm in $\adj[\alg K]{\md E,\md F}$ if and only if $T$ has closed range and both $\dim_{\alg K}\ker T$ and $\dim_{\alg K} \ker{T^*}$ are finite.\label{item2:thm:Psi_isomoorpism}
					\end{enumerate}
				\end{lemma}
				
				We now apply these results to complexes of Hilbert $\alg K$-modules. 
				
				\begin{proof}[Proof of \Cref{thm:topological_fredholm_compact_operators}]
					Without loss of generality, it suffices to prove the statement for the bounded transform $\complex{F_t}$. Denote by $\Delta$ the Laplace operator of $\complex{F_t}$. Let $e_0$ be an arbitrary minimal projection in $\alg K$. 
					
					(i)$\Rightarrow$(ii). Assume $\complex{F_t}$ is $\alg K$-Fredholm. Using \Cref{thm:fredholm-goal-I}, we may apply \Cref{thm:Psi_isomoorpism}\ref{item2:thm:Psi_isomoorpism} to $\Delta$ and derive that it has closed range. By \Cref{thm:hodge-decomposition-equivalences}, each differential map has closed range. Hence, $\complex{F_t}$ is a $\alg K$-Fredholm complex with strong Hodge decomposition and one can compute its index by using \Cref{thm:index_topological_complex}. Recall from \cite[Corollary 6.4.2]{RordamLarsenLaustsen2000} that $\text{Tr}$ induces an isomorphism from $\K[\alg K]$ into $\mathbb Z$.  Therefore, $\dim_{\alg K} \coH{k}{\complex{F_t}}$ is finite for every $k$.
					
					(ii)$\Rightarrow$(i). By \Cref{thm:fredholm-goal-I}, it suffices to prove that $\Delta$ is $\alg K$-Fredholm. We write $\Psi_k$ instead of $\Psi_{\md E_k,\md E_{k+1}}$. Then, for every $k$, we get
					\begin{equation*}
						(\coH{k}{\complex{F_t}})_{e_0}=(\ker{F_{t_k}})_{e_0}/(\ran{F_{t_{k-1}}}) _{e_0}=\ker{\Psi_{k}(F_{t_k})}/\ran{\Psi_{k-1}(F_{t_{k-1}})}.
					\end{equation*}
					Since $\dim (\coH{k}{\complex{F_t}})_{e_0}=\dim_{\alg K}\coH{k}{\complex{F_t}}<\infty$, then $\{\Psi_k(F_{t_k})\}_{k=0}^N$ is a Fredholm Hilbert complex (\Cref{thm:fredholm_complexes_of_hilbert_spaces}). In particular, its Laplace operator is $\mathbb C$-Fredholm whence $\Delta$ is $\alg K$-Fredholm by \Cref{thm:Psi_isomoorpism}\ref{item2:thm:Psi_isomoorpism}. 
                    
                    To compute the index, we use the isomorphism from $\K[\alg K]$ into $\mathbb Z$ induced by the trace map and \Cref{thm:index_topological_complex}.
				\end{proof}
				
				\begin{corollary}
					Let $\complex{T}$, $\complex[\md F]{U}$, and $\complex[\md G]{V}$ be finite-length $\alg K$-Hilbert complexes, and assume there exists an exact sequence of maps of adjointable complexes
					\begin{equation*}
						\begin{tikzcd}
							0 & {\complex{T}} & {\complex[\md F]{U}} & {\complex[\md G]{V}} & 0.
							\arrow[from=1-1, to=1-2]
							\arrow["G", from=1-2, to=1-3]
							\arrow["H", from=1-3, to=1-4]
							\arrow[from=1-4, to=1-5]
						\end{tikzcd}
					\end{equation*}
					If at least two of the complexes $\complex{T}$, $\complex[\md F]{U}$, and $\complex[\md G]{V}$ are $\alg K$-Fredholm, then the third is $\alg K$-Fredholm as well, and
					\begin{equation*}
						\ind{\complex{T}}-\ind{\complex[\md F]{U}}+\ind{\complex[\md G]{V}}=0.
					\end{equation*}
					\begin{proof}
						As in the classical case, see, e.g., \cite[Corollary 1.8]{AmbrozieVasilescu1995}, we can show that the short exact sequence induces a long exact cohomology sequence. This sequence is $\alg K$-Fredholm with index zero (\Cref{thm:exact_complexes_index}), and its Hilbert $\alg K$-modules are finite-rank by \Cref{thm:index_topological_complex}. Hence, its index equals the \emph{Euler characteristic}.
					\end{proof}
				\end{corollary}
				
		\subsection{C*-elliptic complexes}
				Assume that $\alg A$ is a unital C*-algebra, and that $(M,g)$ is a (smooth) compact Riemannian manifold. Let $\xi$ be an $\alg A$-vector bundle over $M$, that is essentially, a locally trivial fiber bundle over $M$ whose fibers, $\xi_x$, are isomorphic to a finitely generated projective Hilbert $\alg A$-module, say $\md F$. As in the classical case, the space of smooth sections, $\Gamma^\infty(\xi)$, can be endowed with a pre-Hilbert $\alg A$-module structure and completed into a Sobolev space $H^s(\xi)$ for any $s\in \mathbb Z$, see, e.g.,  \cite[Section 3]{MiscenkoFomenko1979} or \cite[Section 3]{Krysl2014}. 
				
				\begin{definition}
						The set of $\alg A$-linear maps $T\colon \Gamma^\infty(\xi)\to \Gamma^\infty(\eta)$ such that for each $m\in \mathbb Z$ the map $T$ can be extended to a continuous map $T^{(m)}\colon  H^m(\xi)\to H^{m-r}(\eta)$ will be called the space of operators of order $r\in \mathbb Z$ and denoted by $\mathrm{OP}_r(\xi,\eta)$. 
				\end{definition}
				To ease the notation we simply write $T$ instead of $T^{(m)}$. Let $T^*M$ be the cotangent bundle and define the projection $\pi\colon  T^*M\setminus\{0\}\to M$ in the standard way. The pull-back bundle $\pi^*\xi$ is an $\alg A$-vector bundle whose fiber over the point $(v,x)\in ((TM_x)^*\setminus \{0\})\times M$ is isomorphic to $\xi_x$. For any two $\alg A$-vector bundles, say $\xi$ and $\eta$, we define $\mathrm{Hom}(\pi^*\xi,\pi^*\eta)$  as the space of all operator-valued functions $\sigma$ on $T^*M\setminus\{0\}$ such that $\sigma(v,x)\colon  \xi_x\to \eta_x$ is an adjointable operator, see \cite[p.74]{SolovyovTroitsky2001}. In addition, for each $r\in \mathbb Z$, we define
				\begin{equation*}
					\mathrm{Smbl}_r(\xi,\eta):=\{\sigma\in \mathrm{Hom}(\pi^*\xi,\pi^*\eta)\mid \sigma(\rho v,x)=\rho^r\sigma (v,x),\; \rho>0\}.
				\end{equation*}
				
				In \cite[Chapter 2]{SolovyovTroitsky2001}, the authors provided a generalization of the \emph{Seeley algebra} (see \cite{Seeley1965} or \cite[Chapter XI]{Palais1965}), constructing, for any $r\in \mathbb Z$, a subspace $\mathrm{Int}_r(\xi,\eta)$ of $\mathrm{OP}_r(\xi,\eta)$, called \emph{intregro-differential operators}, and a linear symbol map $\sigma_r\colon  \mathrm{Int}_r(\xi,\eta)\to \mathrm{Smbl}_r(\xi,\eta)$ with the following properties.
				
				\begin{theorem}[{\cite[Theorem 2.1.139]{SolovyovTroitsky2001}}]\label{thm:seeley_properties_I}
					The sequence
					\[\begin{tikzcd}
						0 & {\mathrm{OP}_{r-1}(\xi,\eta)} & {\mathrm{Int}_r(\xi,\eta)} & {\mathrm{Smbl}_r(\xi,\eta)} & 0
						\arrow[from=1-1, to=1-2]
						\arrow[from=1-2, to=1-3]
						\arrow["{\sigma_r}", from=1-3, to=1-4]
						\arrow[from=1-4, to=1-5]
					\end{tikzcd}\]
					is a short exact sequence of vector spaces.
				\end{theorem}
				\begin{theorem}[{\cite[Theorem 2.1.135]{SolovyovTroitsky2001}}]\label{thm:seeley_properties_II}
					Consider $T\in \mathrm{Int}_{r_1}(\xi,\eta)$ and $S\in \mathrm{Int}_{r_2}(\eta,\zeta)$. Then $T$ has a formal transpose, denoted by $T^\dagger$. Moreover, 
				 	\begin{align*}
					 	ST\in \mathrm{Int}_{r_1+r_2}(\xi,\zeta)\quad&\text{ and}\quad\sigma_{r_1+r_2}(ST)=\sigma_{r_2}(S)\sigma_{r_1}(T),\\
					 	T^\dagger\in \mathrm{Int}_{r_1}(\eta,\xi)\quad&\text{and}\quad\sigma_{r_1}(T^\dagger)=(-1)^{r_1}\sigma_{r_1}(T)^*.
					 \end{align*}
				\end{theorem}
				\begin{remark}
					The relation $\sigma_{r_1}(T^\dagger)=(-1)^{r_1}\sigma_{r_1}(T)^*$ can be replaced by $\sigma_{r_1}(T^\dagger)=\sigma_{r_1}(T)^*$. The constant $(-1)^{r_1}$ comes from the definition of the symbol map \cite[Definition 2.1.103]{SolovyovTroitsky2001}, where by following \cite{Seeley1965,Palais1965}, the authors aim to obtain a symbol map that is consistent with that for differential operators, see \cite[p.80]{SolovyovTroitsky2001}. As in the classical case, one can do the opposite, that is, modify the definition of the symbol map for differential operators to agree with the symbol of pseudo-differential operators, see, e.g. \cite[fn. p.115]{Wells2008}. As it is warned by \cite[Section XVI.8]{Palais1965}, the definition of the topological index is replete of sign conventions with arbitrary choices, and a change in any one of these choices calls for a compensating multiplication by a power of $-1$. For our discussion, without loss of generality, we assume that $\sigma_{r_1}(T^\dagger)=\sigma_{r_1}(T)^*$.
				\end{remark}				
				
				We say that $D\in \mathrm{Int}_r(\xi,\eta)$ is $\alg A$-elliptic if $\sigma_r(D)(v,x)\colon \xi_x\to \eta_x$ is an isomorphism for every $(v,x)\in T^*M\setminus\{0\}$. More generally, we consider the following definition.
				\begin{definition}
					Assume that $D_{\bullet}:=\{D_k\in \mathrm{Int}_{r}(\xi_k,\xi_{k+1})\}_{k=0}^N$ is a sequence of integro-differential operators such that $D_{k+1}D_k=0$. Then $D_\bullet$ is a complex, and it is an $\alg A$-elliptic complex if the symbol sequence
					\[\begin{tikzcd}
						0 & {\pi^*\xi_0} & {\pi^*\xi_1} & \cdots & {\pi^*\xi_{N+1}} & 0
						\arrow[from=1-1, to=1-2]
						\arrow["{\sigma_r(D_0)}", from=1-2, to=1-3]
						\arrow["{\sigma_r(D_{1})}", from=1-3, to=1-4]
						\arrow["{\sigma_r(D_{N})}", from=1-4, to=1-5]
						\arrow[from=1-5, to=1-6]
					\end{tikzcd}\]
					is an exact sequence of Hilbert $\alg A$-modules on each fiber (outside of the zero-section).
				\end{definition}
				Let $D_\bullet$ be a finite-length complex of integro-differential operators. As in the classical case, or by analogy with \Cref{subsec:dirac-operator,subsec:laplace-operator}, we can define the \emph{$k$-th Hodge-Laplacian} $\Delta_k:=D_k^\dagger D_k+D_{k-1}D_{k+1}^\dagger$ and the \emph{even Dirac operator} $\mathfrak{D}$, cf. \cite[Section 7]{AtiyahSinger1968}. Moreover, they satisfy the relation $\Delta=\mathfrak{D}^\dagger\mathfrak{D}\oplus \mathfrak{D}\mathfrak{D}^\dagger$. 
				\begin{theorem}
					$\Delta$ is $\alg A$-elliptic if and only if $\mathfrak D$ is $\alg A$-elliptic.
					\begin{proof}
						It is an immediate consequence of \Cref{thm:seeley_properties_II} and definition.
					\end{proof}
				\end{theorem}
				 
				We now prove a classical result in elliptic theory, cf. \cite[4.4.3]{Troitsky1987} and \cite[Corollary 10]{Krysl2014}.
				 
				\begin{theorem}\label{thm:ellipticity_equivalence_complex_Laplace}
					The complex $D_\bullet$ is $\alg A$-elliptic if and only if $\Delta$ is $\alg A$-elliptic.
				 	\begin{proof}
						By \Cref{thm:exact-complex-isomorphism-laplacians}, and \Cref{thm:seeley_properties_I,thm:seeley_properties_II}, 
				 		\begin{equation}\label{eq:symbol_elliptic_operators}
				 			\sigma_{2r}(\Delta_k)=\sigma_r(D_k)^*\sigma_r(D_k)+\sigma_r(D_{k-1})\sigma_r(D_{k-1})^*,\quad \forall k,
				 		\end{equation}
				 		is an isomorphism if and only if the symbol sequence is exact.
				 	\end{proof}
				\end{theorem}
				
				\cite[Theorem 2.1.146]{SolovyovTroitsky2001} showed that any $\alg A$-elliptic operator induces an $\alg A$-Fredholm operator for any of its continuous extensions. More generally, assume that $D_\bullet$ is an $\alg A$-elliptic complex. Then, for any $s\in \mathbb Z$, it induces a complex of Hilbert $\alg A$-modules
				\[\begin{tikzcd}
					0 & {H^{s}(\xi_0)} & \cdots & {H^{s-rN}(\xi_{N})} & {H^{s-r(N+1)}(\xi_{N+1})} & 0
					\arrow[from=1-1, to=1-2]
					\arrow["{D_0}", from=1-2, to=1-3]
					\arrow[from=1-3, to=1-4]
					\arrow["{D_{N}}", from=1-4, to=1-5]
					\arrow[from=1-5, to=1-6]
				\end{tikzcd}\]
				that we denote by $\complex[H^s(\xi)]{D}$; we refer to \cite[p.203]{Krysl2014} to infer how to compute the adjoint of each differential map. In particular, for any $0\leq k\leq N$, the adjoint of $D_k: H^{s-rk}(\xi_{k})\to H^{s-r(k+1)}(\xi_{k+1})$ is given by
                \[
                    D_k^*:=(1+\Delta_g)^{-s+rk}D_k^\dagger(1+\Delta_g)^{s-r(k+1)},
                \]
                where $\Delta_g$ is minus the Laplace-Beltrami operator. Hence, for every $s\in \mathbb Z$ and every $0\leq k\leq N$, the $k$-th Laplace operator of $\complex[H^s(\xi)]{D}$ is the continuous extension of an $\alg A$-elliptic operator because the $k$-th Hodge-Laplacian is $\alg A$-elliptic (by assumption and \Cref{thm:ellipticity_equivalence_complex_Laplace}). It follows from \cite[Theorem 2.1.146]{SolovyovTroitsky2001} that the $k$-th Laplace operator of $\complex[H^s(\xi)]{D}$ is $\alg A$-Fredholm. \Cref{thm:fredholm-goal-I} then implies the following result.
				
				\begin{theorem}
					If $D_\bullet$ is an $\alg A$-elliptic complex, then, for any $s\in \mathbb Z$, the induced complex $\complex[H^s(\xi)]{D}$ is an $\alg A$-Fredholm complex. 
				\end{theorem}
				
				\begin{remark}
					The compact operators defined in \cite{MiscenkoFomenko1979}  coincide with the compact operators in the sense of Kasparov (\Cref{sec:notation}) for unital C*-algebras; see \cite[Definitions 1.2.23 and 1.2.28]{SolovyovTroitsky2001}.
				\end{remark}

    
    \providecommand{\papertitle}[2]{\href{https://doi.org/#2}{#1}}
    \providecommand{\booktitle}[2]{\href{https://doi.org/#2}{\em #1}}

\end{document}